

\magnification=1200
\voffset=.7cm
\hoffset=0cm


\font\tenpc=cmcsc10

\font\eightrm=cmr8
\font\eighti=cmmi8
\font\eightsy=cmsy8
\font\eightbf=cmbx8
\font\eighttt=cmtt8
\font\eightit=cmti8
\font\eightsl=cmsl8
\font\sixrm=cmr6
\font\sixi=cmmi6
\font\sixsy=cmsy6
\font\sixbf=cmbx6

\skewchar\eighti='177 \skewchar\sixi='177
\skewchar\eightsy='60 \skewchar\sixsy='60


\font\tengoth=eufm10
\font\tenbboard=msbm10
\font\eightgoth=eufm7 at 8pt
\font\eightbboard=msbm7 at 8pt
\font\sevengoth=eufm7
\font\sevenbboard=msbm7
\font\sixgoth=eufm5 at 6pt
\font\fivegoth=eufm5

\newfam\gothfam
\newfam\bboardfam

\catcode`\@=11

\def\raggedbottom{\topskip 10pt plus 36pt
\r@ggedbottomtrue}

\def\pc#1#2|{{\bigf@ntpc #1\penalty
\@MM\hskip\z@skip\smallf@ntpc #2}}

\def\tenpoint{%
  \textfont0=\tenrm \scriptfont0=\sevenrm \scriptscriptfont0=\fiverm
  \def\rm{\fam\z@\tenrm}%
  \textfont1=\teni \scriptfont1=\seveni \scriptscriptfont1=\fivei
  \def\oldstyle{\fam\@ne\teni}%
  \textfont2=\tensy \scriptfont2=\sevensy \scriptscriptfont2=\fivesy
  \textfont\gothfam=\tengoth \scriptfont\gothfam=\sevengoth
  \scriptscriptfont\gothfam=\fivegoth
  \def\goth{\fam\gothfam\tengoth}%
  \textfont\bboardfam=\tenbboard \scriptfont\bboardfam=\sevenbboard
  \scriptscriptfont\bboardfam=\sevenbboard
  \def\bboard{\fam\bboardfam}%
  \textfont\itfam=\tenit
  \def\it{\fam\itfam\tenit}%
  \textfont\slfam=\tensl
  \def\sl{\fam\slfam\tensl}%
  \textfont\bffam=\tenbf \scriptfont\bffam=\sevenbf
  \scriptscriptfont\bffam=\fivebf
  \def\bf{\fam\bffam\tenbf}%
  \textfont\ttfam=\tentt
  \def\tt{\fam\ttfam\tentt}%
  \abovedisplayskip=12pt plus 3pt minus 9pt
  \abovedisplayshortskip=0pt plus 3pt
  \belowdisplayskip=12pt plus 3pt minus 9pt
  \belowdisplayshortskip=7pt plus 3pt minus 4pt
  \smallskipamount=3pt plus 1pt minus 1pt
  \medskipamount=6pt plus 2pt minus 2pt
  \bigskipamount=12pt plus 4pt minus 4pt
  \normalbaselineskip=12pt
  \setbox\strutbox=\hbox{\vrule height8.5pt depth3.5pt width0pt}%
  \let\bigf@ntpc=\tenrm \let\smallf@ntpc=\sevenrm
  \let\petcap=\tenpc
  \normalbaselines\rm}
\def\eightpoint{%
  \textfont0=\eightrm \scriptfont0=\sixrm \scriptscriptfont0=\fiverm
  \def\rm{\fam\z@\eightrm}%
  \textfont1=\eighti \scriptfont1=\sixi \scriptscriptfont1=\fivei
  \def\oldstyle{\fam\@ne\eighti}%
  \textfont2=\eightsy \scriptfont2=\sixsy \scriptscriptfont2=\fivesy
  \textfont\gothfam=\eightgoth \scriptfont\gothfam=\sixgoth
  \scriptscriptfont\gothfam=\fivegoth
  \def\goth{\fam\gothfam\eightgoth}%
  \textfont\bboardfam=\eightbboard \scriptfont\bboardfam=\sevenbboard
  \scriptscriptfont\bboardfam=\sevenbboard
  \def\bboard{\fam\bboardfam}%
  \textfont\itfam=\eightit
  \def\it{\fam\itfam\eightit}%
  \textfont\slfam=\eightsl
  \def\sl{\fam\slfam\eightsl}%
  \textfont\bffam=\eightbf \scriptfont\bffam=\sixbf
  \scriptscriptfont\bffam=\fivebf
  \def\bf{\fam\bffam\eightbf}%
  \textfont\ttfam=\eighttt
  \def\tt{\fam\ttfam\eighttt}%
  \abovedisplayskip=9pt plus 2pt minus 6pt
  \abovedisplayshortskip=0pt plus 2pt
  \belowdisplayskip=9pt plus 2pt minus 6pt
  \belowdisplayshortskip=5pt plus 2pt minus 3pt
  \smallskipamount=2pt plus 1pt minus 1pt
  \medskipamount=4pt plus 2pt minus 1pt
  \bigskipamount=9pt plus 3pt minus 3pt
  \normalbaselineskip=9pt
  \setbox\strutbox=\hbox{\vrule height7pt depth2pt width0pt}%
  \let\bigf@ntpc=\eightrm \let\smallf@ntpc=\sixrm
  \normalbaselines\rm}

\tenpoint


\catcode`\;=\active
\def;{\relax\ifhmode\ifdim\lastskip>\z@
\unskip\fi\kern\fontdimen2  -1.2 \fontdimen3 \string;}

\catcode`\:=\active
\def:{\relax\ifhmode\ifdim\lastskip>\z@\unskip\fi\penalty\@M\ \fi\string:}

\catcode`\!=\active
\def!{\relax\ifhmode\ifdim\lastskip>\z@
\unskip\fi\kern\fontdimen2  -1.1 \fontdimen3 \string!}

\catcode`\?=\active
\def?{\relax\ifhmode\ifdim\lastskip>\z@
\unskip\fi\kern\fontdimen2  -1.1 \fontdimen3 \string?}

\def\^#1{\if#1i{\accent"5E\i}\else{\accent"5E #1}\fi}
\def\"#1{\if#1i{\accent"7F\i}\else{\accent"7F #1}\fi}

\frenchspacing


\newif\ifpagetitre
\newtoks\auteurcourant \auteurcourant={\hfil}
\newtoks\titrecourant \titrecourant={\hfil}

\def\appeln@te{}
\def\vfootnote#1{\def\@parameter{#1}\insert\footins\bgroup\eightpoint
  \interlinepenalty\interfootnotelinepenalty
  \splittopskip\ht\strutbox 
  \splitmaxdepth\dp\strutbox \floatingpenalty\@MM
  \leftskip\z@skip \rightskip\z@skip
  \ifx\appeln@te\@parameter\indent \else{\noindent #1\ }\fi
  \footstrut\futurelet\next\fo@t}

\pretolerance=500 \tolerance=1000 \brokenpenalty=5000
\newdimen\hmargehaute \hmargehaute=0cm

\newdimen\lpage \lpage=13.3cm
\newdimen\hpage 
\hpage=20.9cm 
\newdimen\lmargeext \lmargeext=1cm
\hsize=11.67cm
\vsize=18.90cm 
\parskip 0pt
\parindent=12pt

\def\margehaute{\vbox to \hmargehaute{\vss}}%
\def\margebasse{\vss}

\footline={}

\output{\shipout\vbox to \hpage{\margehaute\nointerlineskip
  \corpsdepage\margebasse}
  \advancepageno \global\pagetitrefalse
  \ifnum\outputpenalty>-20000 \else\dosupereject\fi}

\titrecourant={\iftrue\botmark\fi}  

\def\newsectionmark{}

\def\shortsectionmark #1{\mark{\noexpand\eightrm\uppercase{#1}\noexpand\else\noexpand\eightrm
\uppercase{#1}} 
\def\newsectionmark {\noexpand\eightrm\uppercase{#1}}}

\def\corpsdepage{\hbox to \lpage{\hss\pagetexte\hskip\lmargeext}}
\def\pagetexte{\vbox{\makeheadline\kern.4cm\pagebody\makefootline}}
\headline={\ifpagetitre\titleheadline \else
  \ifodd\pageno\rightheadline \else\leftheadline\fi\fi}
\def\leftheadline{\eightpoint\rlap{\tenbf\folio}\hfil
{\eightrm\the\auteurcourant}\hfil}
\def\rightheadline{\eightpoint\hfil
{\eightrm\the\titrecourant}\hfil\llap{\tenbf\folio}}
\def\titleheadline{\hfill}
\pagetitretrue

\def\footnoterule{\kern-6\p@
  \hrule width 2truein \kern 5.6\p@} 

\def\pd#1#2 {\pc#1#2| }

\def\pointir{\discretionary{.}{}{.\kern.35em---\kern.7em}\nobreak
\hskip 0em plus .3em minus .4em }

\def\abstract#1{\vbox{\eightpoint \pc ABSTRACT|\pointir #1}}

\def\titre#1|{\message{#1}
              \par\vskip 30pt plus 24pt minus 3pt\penalty -1000
              \vskip 0pt plus -24pt minus 3pt\penalty -1000
              \centerline{\bf #1}
              \vskip 5pt
              \penalty 10000 }

\def\section #1|{\par\vskip .3cm \penalty -200
       \mark{\newsectionmark\noexpand\else \uppercase{#1}}
         {\bf #1}\shortsectionmark{#1}\unskip\pointir
            }

\def\sectiona #1|{\par\vskip .4cm \penalty -200
       \mark{\newsectionmark\noexpand\else \uppercase{#1}}
\centerline{\bf #1}\shortsectionmark{#1}\par\nobreak\vskip
5pt}

\def\ssection#1|{\par\vskip .2cm
                {\it #1}\pointir}

\long\def\th#1|#2\finth{\par\medskip
              {\petcap #1\pointir}{\it #2}\par\smallskip}

\long\def\tha#1|#2\fintha{\par\medskip
                    {\petcap #1}\par\nobreak{\it #2}\par\smallskip}
\def\cf{{\it cf}}

\def\rem#1|{
\medskip{{\it #1}\pointir}}

\def\rema#1|{\par\medskip
             {{\it #1}\par\nobreak }}

\def\ieme{\raise 1ex\hbox{\pc{}i\`eme|}}
\def\omini{\raise 1ex\hbox{\pc{}o|}}
\def\emini{\raise 1ex\hbox{\pc{}i\`eme|}}
\def\ermini{\raise 1ex\hbox{\pc{}er|}}
\def\remini{\raise 1ex\hbox{\pc{}re|}}

\def\article#1|#2|#3|#4|#5|#6|#7|
    {{\leftskip=7mm\noindent
     \hangindent=2mm\hangafter=1
     \llap{[#1]\hskip.35em}{#2}\pointir
     #3, {\sl #4}, t.\nobreak\ {\bf #5}, {\oldstyle #6},
     p.\nobreak\ #7.\par}}
\def\livre#1|#2|#3|#4|
    {{\leftskip=7mm\noindent
    \hangindent=2mm\hangafter=1
    \llap{[#1]\hskip.35em}{#2}\pointir
    {\sl #3}\pointir #4.\par}}
\def\divers#1|#2|#3|
    {{\leftskip=7mm\noindent
    \hangindent=2mm\hangafter=1
     \llap{[#1]\hskip.35em}{#2}\pointir
     #3.\par}}
\mathchardef\conj="0365
\def\dem{\par{\it D\'emonstration}\pointir}
\def\qed{\quad\raise -2pt\hbox{\vrule\vbox to 10pt{\hrule width 4pt
\vfill\hrule}\vrule}}

\def\virg{\raise 2pt\hbox{,}}   

\def\cqfd{\penalty 500 \hbox{\qed}\par\smallskip}

\long\def\entourer#1{\hbox{\vrule\vbox{\hrule\hbox{\kern15pt\vbox{\kern5pt
{#1}\kern5pt}\kern15pt}\hrule}\vrule}}
\def\\S {\vskip 5pt\hskip .5cm plus .1cm minus .1cm\relax}

\def\enonce#1|#2\finenonce{{\par\leftskip=36pt
\noindent\hbox to 0pt{\kern-\leftskip#1\hfill}{#2}\par}}

\def\senonce#1|#2\finsenonce{{\par\leftskip=36pt
\noindent\hbox to 0pt{\kern-24pt #1\hfill}{#2}\par}}

\def\ssenonce#1|#2\finssenonce{{\par\leftskip=48pt
\noindent\hbox to 0pt{\kern-24pt #1\hfill}{#2}\par}}

\def\decale#1|{\par\noindent\hskip 28pt\llap{#1}\kern 5pt}
\def\decaledecale#1|{\par\noindent\hskip 34pt\llap{#1}\kern 5pt}

\def\titrea#1|#2|{\message{#1 #2}
  \par\vskip.5cm plus .1cm minus .1cm\penalty -1000
  \centerline{\bf #1}
  \centerline{\bf #2}
  \vskip 5pt
  \penalty 10000 }

\def\ssectiona#1|{\par\vskip .2cm
  {\it #1}
  \par\nobreak\vskip 2pt }

\def\rest#1{\ifinner\setbox1=\hbox{$\textstyle{#1}$}
            \else\setbox1=\hbox{$\displaystyle{#1}$}\fi
            \dimen1=\ht1
            \advance\dimen1 by\dp1
            \divide\dimen1 by 2
            \box1\lower 2pt\hbox{$\left|\vbox to\dimen1{}\right.$}}

\def\displaylinesno#1{\displ@y\halign{
\hbox to\displaywidth{$\@lign\hfil\displaystyle##\hfil$}&
\llap{$##$}\crcr#1\crcr}}

\def\ldisplaylinesno#1{\displ@y\halign{ 
\hbox to\displaywidth{$\@lign\hfil\displaystyle##\hfil$}&
\kern-\displaywidth\rlap{$##$}\tabskip\displaywidth\crcr#1\crcr}}

\def\Eqalign#1{\null\,\vcenter{\openup\jot\m@th\ialign{
\strut\hfil$\displaystyle{##}$&$\displaystyle{{}##}$\hfil
&&\quad\strut\hfil$\displaystyle{##}$&$\displaystyle{{}##}$\hfil
\crcr#1\crcr}}\,}


\def\matrixd#1{\null \,\vcenter {\normalbaselines \m@th
\ialign {\hfil $##$&&\quad \hfil $##$\crcr
\mathstrut \crcr \noalign {\kern -\baselineskip } #1\crcr
\mathstrut \crcr \noalign {\kern -\baselineskip }}}\,}

\def\system#1{\left\{\null\,\vcenter{\openup1\jot\m@th
\ialign{\strut$##$     &\hfil$##$&$##$\hfil&&
       \enskip$##$\enskip&\hfil$##$&$##$\hfil\crcr
        #1\crcr}}\right.}

\def\lfq{\leavevmode\raise.3ex\hbox{$\scriptscriptstyle
\langle\!\langle$}\thinspace}
\def\rfq{\leavevmode\thinspace\raise.3ex\hbox{$\scriptscriptstyle
\rangle\!\rangle$}}

\catcode`\@=12


\def\og{\leavevmode\raise 1pt\hbox{$\scriptscriptstyle\langle\!\langle\,$}}
\def\fg{\ignorespaces\leavevmode\raise 1pt\hbox{$\,\scriptscriptstyle\rangle\!\rangle$}}

\def\decale#1|{\par\noindent\hskip 28pt\llap{#1}\kern
5pt } 
\def\decaledecale#1|{\par\noindent\hskip 34pt\llap{#1}\kern
5pt }

\def\card{\mathop{\rm card}\nolimits}

\def\tg{\mathop{\rm tg}}

\def\phiavec_#1{\varphi_{\hbox{\lower 3pt\hbox{$\scriptstyle #1$}}}}

\ifnum\day<0 \advance\day by 160 \fi

\def\fdate{\number\day\space
   \ifcase\month\or janvier\or f\'evrier\or 
mars\or avril\or mai\or 
   juin\or juillet\or ao\^ut\or septembre\or 
octobre\or novembre\or
   d\'ecembre\fi\space\number\year}

\def\illustration #1 by #2 (#3) scaled #4{\dimen1=#2
\divide\dimen1 by 1000
\multiply\dimen1 by #4
\vtop to \dimen1{\dimen1=#1
\divide\dimen1 by 1000
\multiply\dimen1 by #4
\hsize=\dimen1\vss
\noindent\special{illustration #3 scaled #4}}}

\def\Grille{\setbox13=\vbox to 5mm{\hrule width 110mm\vfill}
\setbox13=\vbox{\offinterlineskip
   \copy13\copy13\copy13\copy13\copy13\copy13\copy13\copy13
   \copy13\copy13\copy13\copy13\box13\hrule width 110mm}
\setbox14=\hbox to 5mm{\vrule height 65mm\hfill}
\setbox14=\hbox{\copy14\copy14\copy14\copy14\copy14\copy14
   \copy14\copy14\copy14\copy14\copy14\copy14\copy14\copy14
   \copy14\copy14\copy14\copy14\copy14\copy14\copy14\copy14\box14}
\ht14=0pt\dp14=0pt\wd14=0pt
\setbox13=\vbox to 0pt{\vss\box13\offinterlineskip\box14}
\wd13=0pt\box13}


\def\fleche(#1,#2)\dir(#3,#4)\long#5{%
\noalign{\nointerlineskip\leftput(#1,#2){\vector(#3,#4){#5}}\nointerlineskip}}


\def\hfl#1#2#3{\smash{\mathop{\hbox to#3{\rightarrowfill}}\limits
^{\scriptstyle#1}_{\scriptstyle#2}}}

\def\gfl#1#2#3{\smash{\mathop{\hbox to#3{\leftarrowfill}}\limits
^{\scriptstyle#1}_{\scriptstyle#2}}}


 \message{`lline' & `vector' macros from LaTeX}
 \catcode`@=11
\def\{{\relax\ifmmode\lbrace\else$\lbrace$\fi}
\def\}{\relax\ifmmode\rbrace\else$\rbrace$\fi}
\def\newcount{\alloc@0\count\countdef\insc@unt}
\def\newdimen{\alloc@1\dimen\dimendef\insc@unt}
\def\newwrite{\alloc@7\write\chardef\sixt@@n}

\newwrite\@unused
\newcount\@tempcnta
\newcount\@tempcntb
\newdimen\@tempdima
\newdimen\@tempdimb
\newbox\@tempboxa

\def\@spaces{\space\space\space\space}
\def\@whilenoop#1{}
\def\@whiledim#1\do #2{\ifdim #1\relax#2\@iwhiledim{#1\relax#2}\fi}
\def\@iwhiledim#1{\ifdim #1\let\@nextwhile=\@iwhiledim
        \else\let\@nextwhile=\@whilenoop\fi\@nextwhile{#1}}
\def\@badlinearg{\@latexerr{Bad \string\line\space or \string\vector
   \space argument}}
\def\@latexerr#1#2{\begingroup
\edef\@tempc{#2}\expandafter\errhelp\expandafter{\@tempc}%
\def\@eha{Your command was ignored.
^^JType \space I <command> <return> \space to replace it
  with another command,^^Jor \space <return> \space to continue without it.}
\def\@ehb{You've lost some text. \space \@ehc}
\def\@ehc{Try typing \space <return>
  \space to proceed.^^JIf that doesn't work, type \space X <return> \space to
  quit.}
\def\@ehd{You're in trouble here.  \space\@ehc}

\typeout{LaTeX error. \space See LaTeX manual for explanation.^^J
 \space\@spaces\@spaces\@spaces Type \space H <return> \space for
 immediate help.}\errmessage{#1}\endgroup}
\def\typeout#1{{\let\protect\string\immediate\write\@unused{#1}}}

\font\tenln    = line10
\font\tenlnw   = linew10

\newdimen\@wholewidth
\newdimen\@halfwidth
\newdimen\unitlength 

\unitlength =1pt


\def\thinlines{\let\@linefnt\tenln \let\@circlefnt\tencirc
  \@wholewidth\fontdimen8\tenln \@halfwidth .5\@wholewidth}
\def\thicklines{\let\@linefnt\tenlnw \let\@circlefnt\tencircw
  \@wholewidth\fontdimen8\tenlnw \@halfwidth .5\@wholewidth}

\def\linethickness#1{\@wholewidth #1\relax \@halfwidth .5\@wholewidth}

\newif\if@negarg

\def\lline(#1,#2)#3{\@xarg #1\relax \@yarg #2\relax
\@linelen=#3\unitlength
\ifnum\@xarg =0 \@vline
  \else \ifnum\@yarg =0 \@hline \else \@sline\fi
\fi}

\def\@sline{\ifnum\@xarg< 0 \@negargtrue \@xarg -\@xarg \@yyarg -\@yarg
  \else \@negargfalse \@yyarg \@yarg \fi
\ifnum \@yyarg >0 \@tempcnta\@yyarg \else \@tempcnta -\@yyarg \fi
\ifnum\@tempcnta>6 \@badlinearg\@tempcnta0 \fi
\setbox\@linechar\hbox{\@linefnt\@getlinechar(\@xarg,\@yyarg)}%
\ifnum \@yarg >0 \let\@upordown\raise \@clnht\z@
   \else\let\@upordown\lower \@clnht \ht\@linechar\fi
\@clnwd=\wd\@linechar
\if@negarg \hskip -\wd\@linechar \def\@tempa{\hskip -2\wd\@linechar}\else
     \let\@tempa\relax \fi
\@whiledim \@clnwd <\@linelen \do
  {\@upordown\@clnht\copy\@linechar
   \@tempa
   \advance\@clnht \ht\@linechar
   \advance\@clnwd \wd\@linechar}%
\advance\@clnht -\ht\@linechar
\advance\@clnwd -\wd\@linechar
\@tempdima\@linelen\advance\@tempdima -\@clnwd
\@tempdimb\@tempdima\advance\@tempdimb -\wd\@linechar
\if@negarg \hskip -\@tempdimb \else \hskip \@tempdimb \fi
\multiply\@tempdima \@m
\@tempcnta \@tempdima \@tempdima \wd\@linechar \divide\@tempcnta \@tempdima
\@tempdima \ht\@linechar \multiply\@tempdima \@tempcnta
\divide\@tempdima \@m
\advance\@clnht \@tempdima
\ifdim \@linelen <\wd\@linechar
   \hskip \wd\@linechar
  \else\@upordown\@clnht\copy\@linechar\fi}

\def\@hline{\ifnum \@xarg <0 \hskip -\@linelen \fi
\vrule height \@halfwidth depth \@halfwidth width \@linelen
\ifnum \@xarg <0 \hskip -\@linelen \fi}

\def\@getlinechar(#1,#2){\@tempcnta#1\relax\multiply\@tempcnta 8
\advance\@tempcnta -9 \ifnum #2>0 \advance\@tempcnta #2\relax\else
\advance\@tempcnta -#2\relax\advance\@tempcnta 64 \fi
\char\@tempcnta}

\def\vector(#1,#2)#3{\@xarg #1\relax \@yarg #2\relax
\@linelen=#3\unitlength
\ifnum\@xarg =0 \@vvector
  \else \ifnum\@yarg =0 \@hvector \else \@svector\fi
\fi}

\def\@hvector{\@hline\hbox to 0pt{\@linefnt
\ifnum \@xarg <0 \@getlarrow(1,0)\hss\else
    \hss\@getrarrow(1,0)\fi}}

\def\@vvector{\ifnum \@yarg <0 \@downvector \else \@upvector \fi}

\def\@svector{\@sline
\@tempcnta\@yarg \ifnum\@tempcnta <0 \@tempcnta=-\@tempcnta\fi
\ifnum\@tempcnta <5
  \hskip -\wd\@linechar
  \@upordown\@clnht \hbox{\@linefnt  \if@negarg
  \@getlarrow(\@xarg,\@yyarg) \else \@getrarrow(\@xarg,\@yyarg) \fi}%
\else\@badlinearg\fi}

\def\@getlarrow(#1,#2){\ifnum #2 =\z@ \@tempcnta='33\else
\@tempcnta=#1\relax\multiply\@tempcnta \sixt@@n \advance\@tempcnta
-9 \@tempcntb=#2\relax\multiply\@tempcntb \tw@
\ifnum \@tempcntb >0 \advance\@tempcnta \@tempcntb\relax
\else\advance\@tempcnta -\@tempcntb\advance\@tempcnta 64
\fi\fi\char\@tempcnta}

\def\@getrarrow(#1,#2){\@tempcntb=#2\relax
\ifnum\@tempcntb < 0 \@tempcntb=-\@tempcntb\relax\fi
\ifcase \@tempcntb\relax \@tempcnta='55 \or
\ifnum #1<3 \@tempcnta=#1\relax\multiply\@tempcnta
24 \advance\@tempcnta -6 \else \ifnum #1=3 \@tempcnta=49
\else\@tempcnta=58 \fi\fi\or
\ifnum #1<3 \@tempcnta=#1\relax\multiply\@tempcnta
24 \advance\@tempcnta -3 \else \@tempcnta=51\fi\or
\@tempcnta=#1\relax\multiply\@tempcnta
\sixt@@n \advance\@tempcnta -\tw@ \else
\@tempcnta=#1\relax\multiply\@tempcnta
\sixt@@n \advance\@tempcnta 7 \fi\ifnum #2<0 \advance\@tempcnta 64 \fi
\char\@tempcnta}

\def\@vline{\ifnum \@yarg <0 \@downline \else \@upline\fi}

\def\@upline{\hbox to \z@{\hskip -\@halfwidth \vrule
  width \@wholewidth height \@linelen depth \z@\hss}}

\def\@downline{\hbox to \z@{\hskip -\@halfwidth \vrule
  width \@wholewidth height \z@ depth \@linelen \hss}}

\def\@upvector{\@upline\setbox\@tempboxa\hbox{\@linefnt\char'66}\raise
     \@linelen \hbox to\z@{\lower \ht\@tempboxa\box\@tempboxa\hss}}

\def\@downvector{\@downline\lower \@linelen
      \hbox to \z@{\@linefnt\char'77\hss}}

\thinlines

\newcount\@xarg
\newcount\@yarg
\newcount\@yyarg
\newcount\@multicnt
\newdimen\@xdim
\newdimen\@ydim
\newbox\@linechar
\newdimen\@linelen
\newdimen\@clnwd
\newdimen\@clnht
\newdimen\@dashdim
\newbox\@dashbox
\newcount\@dashcnt
 \catcode`@=12


\newbox\tbox
\newbox\tboxa

\def\leftzer#1{\setbox\tbox=\hbox to 0pt{#1\hss}%
     \ht\tbox=0pt \dp\tbox=0pt \box\tbox}

\def\rightzer#1{\setbox\tbox=\hbox to 0pt{\hss #1}%
     \ht\tbox=0pt \dp\tbox=0pt \box\tbox}

\def\centerzer#1{\setbox\tbox=\hbox to 0pt{\hss #1\hss}%
     \ht\tbox=0pt \dp\tbox=0pt \box\tbox}

\def\image(#1,#2)#3{\vbox to #1{\offinterlineskip
    \vss #3 \vskip #2}}

\def\leftput(#1,#2)#3{\setbox\tboxa=\hbox{%
    \kern #1\unitlength
    \raise #2\unitlength\hbox{\leftzer{#3}}}%
    \ht\tboxa=0pt \wd\tboxa=0pt \dp\tboxa=0pt\box\tboxa}

\def\rightput(#1,#2)#3{\setbox\tboxa=\hbox{%
    \kern #1\unitlength
    \raise #2\unitlength\hbox{\rightzer{#3}}}%
    \ht\tboxa=0pt \wd\tboxa=0pt \dp\tboxa=0pt\box\tboxa}

\def\centerput(#1,#2)#3{\setbox\tboxa=\hbox{%
    \kern #1\unitlength
    \raise #2\unitlength\hbox{\centerzer{#3}}}%
    \ht\tboxa=0pt \wd\tboxa=0pt \dp\tboxa=0pt\box\tboxa}

\unitlength=1mm

\def\put(#1,#2)#3{\noalign{\nointerlineskip
                               \centerput(#1,#2){$#3$}
                                \nointerlineskip}}

{\nopagenumbers
\parindent=0pt
\pageno=0

\font\bighelvetica=cmss17 at 30 pt 
\font\midhelvetica=cmss17 at 25pt
\font\midmidhelvetica=cmss12 at 20pt 
\font\smallhelvetica=cmss12  

\noindent
{\bighelvetica
Lecture Notes in

\medskip
Mathematics}

\bigskip
{\smallhelvetica
A collection of informal reports and seminars

Edited by A. Dold, Heidelberg and B. Eckmann, Z\"urich

\bigskip
Series: Institut de Math\'ematique, Universit\'e de Strasbourg

Advisers: P. A. Meyer and M. Karoubi

}
\bigskip\bigskip

\vfill
{\midhelvetica 138
\medskip
\line{\hrulefill}
\bigskip
Dominique Foata}

\smallskip
{\smallhelvetica Universit\'e de Strasbourg}

\bigskip

{\midhelvetica
Marcel-P. Sch\"utzenberger}

\smallskip
{\smallhelvetica Universit\'e de Paris}

\vskip 2cm
{\midhelvetica Th\'eorie G\'eom\'etrique

\smallskip
des Polyn\^omes Eul\'eriens}

\medskip
\line{\hrulefill}
\vfill
\vfill
{\midmidhelvetica Springer-Verlag

\smallskip
Berlin $\cdot$ Heidelberg $\cdot$ New York 1970} 

\eject

}

{

\pagetitretrue

\def\leftheadline{\hfil}
\auteurcourant={}\titrecourant={}
\null

\nopagenumbers
\footnote{}{The present volume is a 2005 {\it revised} version
of the text that was originally published in the
Springer-Verlag Lecture Notes in Mathematics Series, back in
1970. We thank Anne Berstel for her careful proof-reading.}
\vfill\eject  

}

{

\overfullrule=0pt
\pagetitretrue
\titrecourant={\eightrm TABLE DES MATI\`ERES}
\auteurcourant={\eightrm TABLE DES MATI\`ERES}
\vglue 1.5cm
\centerline{{\bf TABLE DES MATI\`ERES}}

\pageno=-3
\def\flead{\leaders\hbox to 4pt{\hss.}\hfill}

\long\def\somr #1|#2|#3|{\hbox{\kern
23pt \vbox{\hsize=10.15cm
\parindent=-23pt
\strut{\pc CHAPITRE|} {#1}.\quad
{\bf #2}\flead\null}\hskip 1pt \hbox to 20pt{\hfill #3}}
\smallskip\goodbreak}

\long\def\somrr #1|#2|#3|{\hbox{\kern
23pt \vbox{\hsize=10.15cm
\strut {#1}.\quad
{#2}\flead\null}\hskip 1pt \hbox to 20pt{\hfill #3}}}

\vskip .5cm

{\parindent=0pt

\vskip 8pt minus 2pt
\somr 0|Introduction et historique des nombres d'Euler|1|
\somrr 1|Bref historique sur les nombres d'Euler|1|
\somrr 2|R\'esum\'e du m\'emoire|2|

\bigskip
\somr {\pc PREMIER|}|Propri\'et\'es g\'en\'erales des syst\`emes
d'exc\'e\-dances et de mont\'ees|5|

\somrr 1|Exc\'edances|5|

\somrr 2|Descentes et mont\'ees|6|

\somrr 3|La transformation fondamentale|7|

\somrr 4|Relations entre les exc\'edances et les descentes|8|

\somrr 5|Applications aux permutations altern\'ees|9|

\somrr 6|Relations entre les exc\'edances et les mont\'ees|10|

\somrr 7|Relations avec les permutations circulaires|11|

\somrr 8|Tableau des bijections utilis\'ees|12|

\somrr 9|Notations g\'en\'erales|12|

\bigskip
\somr 2|Les polyn\^omes eul\'eriens|14|

\somrr 1|Interpr\'etation des polyn\^omes eul\'eriens|14|

\somrr 2|Propri\'et\'es de sym\'etrie|15|

\somrr 3|Relations de r\'ecurrence|16|

\somrr 4|Relations avec le \og probl\`eme de Simon
Newcomb\fg|17|

\somrr 5|Relations avec les nombres de Stirling|18|

\somrr 6|Les identit\'es de Worpitzky|19|

\somrr 7|Table des polyn\^omes eul\'eriens|21|

\bigskip
\somr 3|La formule exponentielle|23|

\somrr 1|La formule de Hurwitz|23|

\somrr 2|Le compos\'e partitionnel|24|

\somrr 3|Une formule d'inversion pour les s\'eries
exponentielles|26|

\somrr 4|Le compos\'e partitionnel des applications|27|

\somrr 5|Applications|29|

\somrr 6|Une identit\'e entre d\'eterminants et permanents|30|

\bigskip
\somr 4|Fonctions g\'en\'eratrices des polyn\^omes
eul\'eriens|32|

\somrr 1|Fonction g\'en\'eratrice exponentielle de
${}^{0\kern-2pt}A_n(t)$, $A_n(t)$ et $B_n(t)$|32|

\somrr 2|Fonction g\'en\'eratrice exponentielle des
polyn\^omes
${}^{r\kern-3pt}A_n(t)$|34|

\somrr 3|Autres interpr\'etations des polyn\^omes eul\'eriens|36|

\bigskip\vfill\eject
\somr 5|Les sommes altern\'ees $A_n(-1)$ et $B_n(-1)$|38|

\somrr 1|Distribution du nombre des descentes sur
${\goth S}_n'$|38|

\somrr 2|Applications aux polyn\^omes eul\'eriens|40|

\somrr 3|Applications aux polyn\^omes $B_n(t)$|40|

\somrr 4|Les d\'eveloppements de $\tg u$ et de
$1/\cos u$|42|

\somrr 5|Table des nombres d'Euler|43|

\bigskip
\hbox{\kern
23pt \vbox{\hsize=10.15cm
\parindent=-23pt
\strut {\bf Bibliographie}\flead\null}\hskip 1pt \hbox to
20pt{\hfill 44}}

}}

\vfill\eject

\pageno=1
\pagetitretrue
\auteurcourant={ CHAPITRE 0 :
INTRODUCTION ET HISTORIQUE DES NOMBRES D'EULER}

\vglue 2cm
\centerline{{\eightrm CHAPITRE 0}}
\vskip 2mm
\centerline{{\bf INTRODUCTION ET HISTORIQUE}} 

\smallskip
\centerline{{\bf DES
NOMBRES D'EULER}}
\vskip 6mm plus 2mm
\sectiona 1. Bref historique sur les nombres d'Euler|
On sait depuis Euler que la relation
$$
\sum_{n\ge 0} A_n(t) {u^n\over n!}={1-t\over
-t+\exp(u(t-1))}\eqno(1)
$$ 
d\'efinit des polyn\^omes sym\'etriques
$$
A_0(t)=1\quad
{\rm et}\quad
A_n(t)=t^{n-1}A_n(t^{-1})=\sum_{0\le k\le n-1}
A_{n,k}t^k\quad (n\ge 1),
$$
de degr\'e $(n-1)$, dont les coefficents $A_{n,k}$ sont des
entiers positifs de somme $A_n(1)=n!$

Worpitzky [31] a donn\'e la formule
$$
x^n=\sum_{0\le k\le n-1}A_{n,k}{x+k\choose n}\eqno(2)
$$
et Frobenius [13], qui a appel\'e les $A_n(t)$ {\it
polyn\^omes eul\'eriens}, a indiqu\'e l'identit\'e
$$
A_n(t)=\sum_{0\le k\le n-1} (n-k)!\,(t-1)^kS(n,n-k),
\eqno(3)
$$
dans laquelle $S(n,j)$ d\'esigne le nombre de Stirling de
seconde esp\`ece, c'est-\`a-dire le nombre de partitions
en~$j$ classes d'un ensemble de~$n$ \'el\'ements. Une
bibliographie de cette question a \'et\'e rassembl\'ee par
Carlitz [4].

De par ailleurs, les polyn\^omes eul\'eriens apparaissent dans
divers probl\`emes d'\'enum\'eration concernant le groupe
sym\'etrique ${\goth S}_n$ sur l'ensemble totalement
ordonn\'e $[\,n\,]=\{1,2,\ldots, n\}$ ($=\emptyset$ pour
$n=0$). Ainsi, d'apr\`es (1) la somme altern\'ee
$(-1)^{p-1}A_{2p-1}(-1)$ est le coefficient de
$u^{2p-1}/(2p-1)!$ dans le d\'eveloppement de $\tg u$
(voir chapitre~V \S\kern2pt 4 du pr\'esent article) et
D\'esir\'e Andr\'e [1], [2] (voir aussi [21], chap.~4) a
d\'ecouvert que ce dernier nombre est celui des
permutations $\sigma\in {\goth S}_{2p-1}$ qui sont {\it
altern\'ees}, c'est-\`a-dire telles que pour chaque $j\in
[p-1]$, on ait \`a la fois $\sigma(2j)<\sigma(2j-1)$ et
$\sigma(2j)<\sigma(2j+1)$.

\goodbreak
MacMahon (19] a \'etudi\'e le nombre $A'_{n,k}$ des
permutations $\sigma\in {\goth S}_n$ telles que
$j<\sigma(j)$ pour exactement~$k$ \'el\'ements $j\in
[\,n\,]$ et en application de son \og Master Theorem\fg,
il a montr\'e que $A'_{n,k}$ est aussi le nombre des
$\sigma\in {\goth S}_n$ qui poss\`edent~$k$ {\it
mont\'ees}, c'est-\`a-dire qui satisfont \`a
$\sigma(j)<\sigma(j+1)$ pour exactement~$k$ valeurs
$j\in [\,n-1\,]$. Carlitz et Riordan~[7] ont reconnu que
les $A'_{n,k}$ sont pr\'ecis\'ement les coefficients des
polyn\^omes eul\'eriens et Riordan a g\'en\'eralis\'e ceux-ci au
moyen de sa th\'eorie des \og rook polynomials\fg. Appelons
$r$-{\it exc\'edance} de $\sigma\in {\goth S}_n$ tout
$j\in [\,n\,]$ tel que $j+r\le \sigma(j)$ $(0\le r$) et soit
${}^{r\kern-3pt}A_{n,k}$ le nombre des $\sigma\in {\goth S}_n$
ayant $k$ $r$-exc\'edances
$({}^{1\kern-3pt}A_{n,k}=A_{n,k})$ ; Riordan, dans le
dernier chapitre de son livre [24], consid\`ere avec des
notations un peu diff\'erentes les polyn\^omes
$$
{}^{r\kern-3pt}A_n(t)=\sum_{0\le k}
{}^{r\kern-3pt}A_{n,k}\,t^k\eqno(4)
$$
et \`a leur sujet \'etablit des extensions des formules (1) et
(3). D'autres g\'en\'era\-lisations sont dues \`a Carlitz et son
\'ecole ([5], [8], [10]).

Soit $|\Delta'M\sigma|$ le nombre des {\it mont\'ees} de 
$\sigma\in {\goth S}_n$. Roselle [25] a calcul\'e le
polyn\^ome $B_n(t)$ d\'efini comme la somme de
$t^{|\Delta'M\sigma|}$ \'etendue au sous-ensemble~${\cal
G}_n$ des $\sigma\in {\goth S}_n$ tels que $1\not
=\sigma(1)$ et $1+\sigma(j)\not=\sigma(j+1)$ pour
chaque $j\in [\,n-1\,]$. Il note que $B_n(1)$ est le
nombre de permutations sans points fixes de ${\goth
S}_n$ et il constate que la somme altern\'ee $(-1)^p
B_{2p}(-1)$ est \'egale au coefficient de $u^{2p}/(2p)!$
dans le d\'eveloppement de $1/\cos u$, c'est-\`a-dire au
$2p$\ieme~{\it nombre d'Euler}. Or, on sait depuis
longtemps que ce coefficient est le nombre des
$\sigma\in{\goth S}_{2p}$ tels que $\sigma(2p-1)<2p$
et $\sigma(2j-1)<\sigma(2j)>\sigma(2j+1)$ $(j\in [p-1])$
([1], [14]), ceci \'etant d'ailleurs un cas particulier d'une
formule plus g\'en\'erale due \`a Entringer [11] et dans une
direction assez diff\'erente de la th\'eorie des \og runs up
and down\fg\ d\'evelopp\'ee par David et Barton \`a des fins
statistiques [3].

Enfin, dans la th\'eorie dite du \og probl\`eme de
Newcomb\fg, on consid\`ere au lieu de l'ensemble {\it
ordonn\'e}
$[\,n\,]$ un ensemble {\it pr\'eordonn\'e} quelconque~$X$.
Les \'enonc\'es y d\'ependent donc de fa\c con cruciale de la
structure de~$X$, ce qui conduit \`a une probl\'ematique
sensiblement diff\'erente, sauf si l'on r\'eintroduit des
hypoth\`eses particuli\`eres sur~$X$ comme par exemple dans
le cas des polyn\^omes de Shanks [27], des polyn\^omes de
Poussin [22] ou dans celui de la \og sp\'ecification
$(1^r(n-r))$\fg\ qui fait appara\^\i tre directement les
entiers ${}^{r\kern-3pt}A_{n,k}$. Hormis ce dernier cas,
nous avons enti\`erement laiss\'e de c\^ot\'e le probl\`eme de
Newcomb qui nous e\^ut entra\^\i n\'e fort loin des polyn\^omes
eul\'eriens. Au demeurant, les m\^emes techniques de base ont
\'et\'e employ\'ees r\'ecemment par l'un de nous [9] pour traiter
le cas g\'en\'eral et certaines de ses applications.

\sectiona 2. R\'esum\'e du m\'emoire|

Les th\'eor\`emes qui viennent d'\^etre rappel\'es ont \'et\'e, en
r\`egle g\'en\'erale, \'etablis en utilisant conjointement quelques
propri\'et\'es g\'eom\'etriques (\og combi\-na\-toires\fg) des
permutations et les m\'ethodes plus exp\'editives du calcul
diff\'erentiel et int\'egral. En particulier, aucune connexion
sauf la co\"\i ncidence de l'aboutissement de deux s\'eries
de calculs ne semble avoir \'et\'e vue entre les sommes
altern\'ees $A_n(-1)$ et $B_n(-1)$ et la signification des
polyn\^omes $A_n(t)$ et $B_n(t)$ en termes d'exc\'edances
ou de mont\'ees.

Le but du pr\'esent m\'emoire est au contraire de d\'evelopper
la th\'eorie g\'eom\'etrique sous-jacente et c'est de fa\c con
subsidiaire que nous en d\'eduisons des identit\'es entre
s\'eries ou polyn\^omes. Nous nous sommes cependant attach\'es
\`a toujours retrouver les r\'esultats classiques. Dans de
nombreux cas, cette approche \'evacue pratiquement tous
les calculs : c'est ce qui se produit par exemple en ce qui
concerne les formules reliant les polyn\^omes eul\'eriens et
les nombres de Stirling. Dans d'autres cas, nous obtenons
des s\'eries d'identit\'es nouvelles (par exemple les \og
formules sommatoires\fg\ g\'en\'eralisant celle de Worpitzky
donn\'ees dans la section~6 du chapitre~II ou le
th\'eor\`eme~5.6). Les m\'ethodes du chapitre~III contiennent
implicitement l'\'enum\'eration du nombre des exc\'edances
pour les permutations dont les longueurs des cycles
satisfont \`a des conditions de divisibilit\'e donn\'ees.

Plus important nous semble la d\'emonstration du fait que
toutes les identit\'es classiques concernant les polyn\^omes
eul\'eriens sont seulement la traduction de propri\'et\'es tr\`es
simples des morphismes d'ensembles totalement ordonn\'es
finis. Pour l'essentiel, elles d\'erivent soit de m\'ethodes
\'el\'ementaires courantes comme l'inversion de M\"obius ou
la formule exponentielle, soit d'une op\'eration unique
nouvelle $\sigma\mapsto \hat\sigma$ appel\'ee ici {\it
transformation fondamentale}, d\'ej\`a introduite par l'un de
nous [12] dans le cadre g\'en\'eral du probl\`eme de Newcomb.
Simultan\'ement, les \'enonc\'es que nous proposons, expriment,
en r\`egle, des bijections entre ensembles. Ils sont donc plus
riches que les identit\'es \'enum\'eratives classiques
auxquelles ils se r\'eduisent quand, en fin de calcul, on
substitue \`a ces ensembles le nombre de leurs \'el\'ements.

Le chapitre I est consacr\'e \`a l'\'etude d\'etaill\'ee des
propri\'et\'es de la transformation fondamentale. Celle-ci est
une bijection de~${\goth S}_n$ sur lui-m\^eme ayant la
propri\'et\'e que l'ensemble des \og exc\'edances\fg\
de~$\sigma$ est envoy\'e, de fa\c con biunivoque, sur celui
des \og descentes\fg\ de~$\hat\sigma$. Elle permet ici
d'\'etablir que la distribution du nombre des exc\'edances sur
les permutations de~${\goth S}_n$ est la m\^eme que sur le
sous-ensemble ${\cal C}_{n+1}$ des {\it permutations
circulaires} de~${\goth S}_{n+1}$. Ce r\'esultat est \'etendu
dans le paragraphe~4, o\`u nous employons les {\it
bi-exc\'edances} (c'est-\`a-dire les $j\in [\,n\,]$ tels que~$j$
soit {\it \`a la fois} strictement plus petit
que~$\sigma(j)$ et que $\sigma^{-1}(j)$) et les {\it
creux} (c'est-\`a-dire les $j\in [\,n\,]$ tels que~$j$
soit {\it \`a la fois} strictement plus petit
que~$\sigma(j-1)$ et que $\sigma(j+1)$) pour \'etudier
les sommes altern\'ees. Cette dualit\'e pourrait \^etre \'etendue
\`a des constructions plus complexes sur lesquelles nous
reviendrons peut-\^etre dans un autre travail.

Dans le chapitre II, nous retrouvons et g\'en\'eralisons
diverses formules de r\'ecurrence de Riordan et l'identit\'e de
Worpitzky, en application des r\'esultats pr\'ec\'edents et de la
consid\'erations des {\it morphismes} $[\,n\,]\rightarrow
[m]$, c'est-\`a-dire, puisque $[\,n\,]$ et $[m]$ sont des
ensembles totalement ordonn\'es des applications $\phi
:[\,n\,]\rightarrow [m]$ telles que $i\le j$ implique
$\phi(i)\le \phi(j)$.

\vfill\eject
Dans le chapitre III, nous croyons utile de donner d'abord
une th\'eorie syst\'ematique de la formule exponentielle
classique de Cauchy exprimant le groupe sym\'etrique en
fonction des permutations circulaires. Cette formule est
un cas particulier d'une constructuon tr\`es g\'en\'erale
permettant de ramener divers probl\`emes d'\'enum\'eration \`a
un probl\`eme analogue sur une sous-famille \og
g\'en\'eratrice\fg\ constitu\'ee par des objets \og connexes\fg.
Afin de clarifier ces notions, nous donnons quelques
\'enonc\'es sous une forme qui permettrait de traiter les
\'enum\'erations d'arborescences. Retournant aux
permutations, une application de cette formule et des
r\'esultats du chapitre~I nous permet d'obtenir dans le
chapitre~IV la fonction g\'en\'eratrice exponentielle du
nombre des $r$-exc\'edances pour les permutations ayant
une composition en cycles donn\'ee.

Dans le chapitre V, nous \'etablissons un th\'eor\`eme sur la
distribution du nombre des mont\'ees pour les permutations
ayant un nombre de creux fix\'e. De fa\c con plus explicite,
soit ${\goth S}_{n,k}'$ l'ensemble des permutations
$\sigma\in {\goth S}_n$ ayant~$k$ creux et telles que
$\sigma(1)=n$ $(2\le 2k\le n)$ ; alors le nombre de
$\sigma\in {\goth S}_{n,k}'$ ayant~$j$ mont\'ees est
donn\'e par le coefficient binomial $n-2k\choose j$ $(0\le
j\le n-2k)$. Nous d\'eduisons de ce r\'esultat les expressions
des sommes altern\'ees $A_n(-1)$ et $B_n(-1)$ en
fonction du nombre des permutations altern\'ees. En
particulier, les d\'eveloppements de $\tg u$ et $1/\cos u$
sont obtenus sans calcul \`a partir de l'expression des
fonctions g\'en\'eratrices exponentielles des polyn\^omes
$A_n(t)$ et $B_n(t)$ donn\'ees dans le chapitre~IV.

Dans tout ce travail, \'etant donn\'es deux ensembles $A$
et $B$ finis et une application $\phi:A\rightarrow B$,
nous appellerons {\it ensemble pond\'er\'e} l'application
$\phi^{\#}:B\rightarrow{\bboard N}$ d\'efinie pour chaque
$b\in B$ par $\phi^{\#}(b)=\card \phi^{-1}(b)$. Par
abus de notation, on d\'esignera par $\phi A$ l'ensemble
pond\'er\'e~$\phi^{\#}$ et il sera commode
d'identifier~$\phi A$ et~$\phi^{\#}$ \`a l'\'el\'ement
$\sum\{\phi^{\#}(b):b\in B\}$ du $\bboard Q$-module
libre de base~$B$.

Nous sommes reconnaissants au professeur J. Riordan de
nous avoir fait b\'en\'eficier de ses conseils et de son
\'erudition. La dactylographie de ce m\'emoire est due \`a Mlle
Cler, du d\'epartement de math\'ematique de Strasbourg, que
nous tenons \`a remercier.

\vfill\eject

\pagetitretrue

\auteurcourant={{\eightrm 
CHAPITRE PREMIER : SYST\`EMES D'EXC\'EDANCES ET
DE MONT\'EES}}

\vglue 2cm
\centerline{{\eightrm CHAPITRE PREMIER}}
\vskip 2mm
\centerline{{\bf PROPRI\'ET\'ES G\'EN\'ERALES DES
SYST\`EMES}} 

\smallskip
\centerline{{\bf D'EXC\'EDANCES ET DE MONT\'EES}}
\vskip 6mm plus 2mm
\sectiona 1. Exc\'edances|
Dans tout ce chapitre, nous utilisons la notation
$x_+$ pour d\'esigner la {\it partie positive}
$x_+=\max\{0,x\}$ de tout $x\in {\bboard Z}$ et pour
$\sigma\in {\goth S}_n$, nous d\'esignons par $\sigma w$
le $n$-uple (ou vecteur) $(\sigma(1),\sigma(2),\ldots,
\sigma(n))\in {\bboard N}^n$. Par cons\'equent, $\sigma
w=\emptyset$ pour $n=0$.

\rem D\'efinition $1.1$|Pour $\sigma\in {\goth S}_n$, le
{\it syst\`eme des $0$-exc\'edances} de~$\sigma$ est le
$n$-uple $E\sigma=(E\sigma(1),
E\sigma(2),\ldots,E\sigma(n))\in {\bboard N}^n $, o\`u pour
chaque $k\in [\,n\,]$, on pose
$$
E\sigma(k)=\bigl(\sigma(k)-(k-1)\bigr)_+\;.
$$

Par exemple, avec $n=6$ et $\sigma w=(6,4,1,2,5,3)$, on a
\hfil\break $E\sigma=((6-0)_+, (4-1)_+, (1-2)_+,
(2-3)_+, (5-4)_+, (3-5)_+)=(6,3,0,0,1,0)$.

\rem D\'efinition $1.2$|Quelque soit l'entier $p\ge 1$, on
note $\Delta$, $\Delta'$ et $\Delta''$ les applications de
${\bboard N}^p$ dans ${\bboard N}^{p-1}$ envoyant
respectivement chaque vecteur $x=(x_1,x_2,\ldots, x_p)\in 
{\bboard N}^p$ sur
$$
\leqalignno{\Delta x&=\bigl((x_1-1)_+, (x_2-1)_+, \ldots,
(x_{p-1}-1)_+\bigr);\cr
\Delta' x&=(x_2,x_3,\ldots, x_p);\cr
\Delta'' x&=(x_1,x_2,\ldots, x_{p-1}).\cr
}
$$

Il est imm\'ediat que les trois op\'erateurs $\Delta$,
$\Delta'$, $\Delta''$ ainsi d\'efinis commutent deux \`a deux.
Prenant le m\^eme exemple que ci-dessus, on obtient :
$E\sigma=(6,3,0,0,1,0)\
(=\Delta^0E\sigma=\Delta'^0E\sigma=\Delta''^0E\sigma)$ ;
$\Delta E\sigma=(5,2,0,0,0)$ ;
$\Delta'E\sigma=(3,0,0,1,0)$ ;
$\Delta''E\sigma=(6,3,0,0,1)$ ;
$\Delta^2E\sigma=(4,1,0,0)$ ;
$\Delta\Delta'E\sigma=\Delta'\Delta E\sigma=(2,0,0,0)$ ;
$\Delta'^2E\sigma=(0,0,1,0)$, etc.

On notera que $\Delta E\sigma$ est simplement la suite 
$$\bigl((\sigma(1)-1)_+, (\sigma(2)-2)_+, \ldots,
(\sigma(n-1)-(n-1))_+\bigr)$$ et que plus g\'en\'eralement
le vecteur $\Delta^rE$ d\'ecrit les $r$-{\it exc\'edances}
de~$\sigma$ mentionn\'ees dans l'introduction. L'une des
raisons motivant l'introduction de $\Delta'$ est contenue
dans le lemme~1.4 ci-dessous. L'op\'erateur $\Delta''$
permettra dans le deuxi\`eme chapitre de formuler une
int\'eressante propri\'et\'e de sym\'etrie des polyn\^omes
eul\'eriens (propri\'et\'e~2.3).

\goodbreak
\rem Remarque $1.3$|Pour $\sigma\in {\goth S}_n$ et
$k\in [\,n\,]$, on a $E\sigma(k)=1$ si et seulement si
$(\sigma(k)-k+1)_+=1$, c'est-\`a-dire si $k=\sigma(k)$
est un {\it point fixe} de~$\sigma$.

\th Lemme 1.4|Soit $\zeta\in {\goth S}_n$ la
permutation circulaire envoyant~$n$ sur~$1$ et chaque
$k\le n-1$ sur $k+1$ ou encore la permutation d\'efinie
par $\zeta w=(2,3,\ldots, n,1)$. Pour chaque $r\ge 0$ et
chaque $\sigma\in {\goth S}_n$ on a
$$
\Delta'^rE\sigma=\Delta^rE\sigma\zeta^r.
$$
\finth

\dem
Posons $\sigma'=\sigma\zeta^r$. Pour chaque $k\in
[n-r]$, on a $\sigma'(k)=\sigma(r+k)$. Donc on obtient
$\Delta'^rE\sigma(k)=E\sigma(r+k)
=(\sigma(r+k)+1-r-k)_+=(\sigma'(k)+1-k-r)_+
=((\sigma'(k)+1-k)_+-r)_+=\Delta^rE\sigma'(k)$.\cqfd

\medskip
Ce simple r\'esultat a la cons\'equence imm\'ediate suivante
qui nous servira fr\'equemment par la suite.

\th Th\'eor\`eme 1.5|Quelque soit le mon\^ome $\Gamma$ de
degr\'e $r\ge 0$ en les applications~$\Delta$
et~$\Delta'$, les ensembles pond\'er\'es $\Delta^rE{\goth
S}_n$, $\Delta'^rE{\goth S}_n$ et $\Gamma E{\goth
S}_n$ sont \'egaux.
\finth

\dem
Puisque $\sigma\mapsto \sigma\zeta^r $ est une
bijection de ${\goth S}_n$ sur lui-m\^eme, l'\'egalit\'e
$\Delta^rE{\goth S}_n=\Delta'^r E{\goth S}_n$ r\'esulte
imm\'ediatement du lemme~1.4. Comme $\Delta$ et
$\Delta'$ commutent, on peut \'ecire
$\Gamma=\Delta^s\Delta'^{r-s}$, d'o\`u $\Gamma E{\goth
S}_n=\Delta^rE{\goth S}_n$.\cqfd

\sectiona 2. Descentes et mont\'ees|
En parall\`ele avec les exc\'edances, nous introduisons le
syst\`eme des 0-{\it descentes}, $D\sigma$, et des 0-{\it
mont\'ees}, $M\sigma$, de $\sigma\in {\goth S}_n$ par la
d\'efinition suivante, dans laquelle on convient que
$\sigma(0)=\sigma^{-1}(0)=\sigma(n+1)=0$.

\rem D\'efinition $1.6$|Pour $\sigma \in {\goth S}_n$, on
pose
$$
\leqalignno{D\sigma
&=\bigl(D\sigma(1),D\sigma(2),\ldots,D\sigma(n)\bigr)
\in {\bboard N}^n;\cr
M\sigma
&=\bigl(M\sigma(1),M\sigma(2),\ldots,M\sigma(n)\bigr)
\in {\bboard N}^n;\cr
\noalign{\hbox{o\`u pour chaque $k\in [\,n\,]$}}
D\sigma(k)
&=\bigl(\sigma(-1+\sigma^{-1}(k))-(k-1)\bigr)_+;\cr
M\sigma(k)
&=\bigl(\sigma(1+\sigma^{-1}(k-1))-(k-1)\bigr)_+.\cr
}
$$

\smallskip
Par exemple, prenant encore $\sigma w=(6,4,1,2,5,3)$, on
obtient\hfil\break
$D\sigma=\bigl((4-0)_+,(1-1)_+,(5-2)_+,
(6-3)_+,(2-4)_+,(0-5)_+\bigr)
=(4,0,3,3,0,0)$ 
$M\sigma=\bigl((6-0)_+,(2-1)_+,(5-2)_+,
(0-3)_+,(1-4)_+,(3-5)_+\bigr)
=(6,1,3,0,0,0)$.

\medskip
Par construction, tous les termes de $D\sigma$ sont nuls
ou sup\'erieurs ou \'egaux \`a~2. D'autre part, $D\sigma(n)$
est toujours nul. Par cons\'equent, $D\sigma$ et
$\Delta D\sigma$ ont le m\^eme nombre de termes
(strictement) positifs. Il est clair que $\Delta^rD\sigma$
et $\Delta^{r-1}M\sigma$ d\'ecrivent les diff\'erences
sup\'erieures ou \'egales \`a~$r$ $(r\ge 1)$ entre termes
cons\'ecutifs de~$\sigma w$, la connexion entre~$D$
et~$M$ \'etant explicit\'ee dans le lemme~1.7 ci-dessous. Il
est encore utile de noter que les termes positifs de
$\Delta'^rD\sigma$ (ou de $\Delta'^r\Delta D\sigma$)
correspondent aux paires $(j-1,j)$ $(0\le j-1)$ telles que
$\sigma(j-1)>\sigma(j)>r$.

\goodbreak
\th Lemme 1.7|Soit $\sigma\mapsto \tilde \sigma$ la
bijection de ${\goth S}_n$ sur lui-m\^eme d\'efinie par
l'identit\'e $\tilde\sigma(k)=\sigma(n+1-k)$ $(k\in [\,n\,])$.
On a $M\tilde\sigma(1)=\sigma(n)$ et
$M\tilde\sigma(k+1)=\Delta D\sigma(k)$ pour chaque $k\in
[n-1]$.
\finth

\dem
Par d\'efinition
$M\tilde\sigma(1)=(\tilde\sigma(1)-0)_+=\tilde\sigma(1)$
et $\tilde\sigma(1)=\sigma(n+1-1)=\sigma(n)$. Soit
$k=\tilde\sigma(j)$ avec $k\in [n-1]$ ; on a alors
$\tilde\sigma^{-1}(k)=j$ et $\sigma(n+1-j)=k$. D'o\`u il
vient
$$
\leqalignno{
M\tilde\sigma(k+1)&\!=\!(\tilde\sigma(1+j)-k)_+\!
=\!(\sigma(n-j)-k)_+
\!=\!(\sigma(-1+(n+1-j))-k)_+\cr
&\!=\!(\sigma(-1+\sigma^{-1}(k))-(k-1)-1)_+
\!=\!\Delta D\sigma(k).\qed\cr}
$$

\smallskip
Prenant $\sigma$ comme dans l'exemple ci-dessus, on
trouve
$\Delta D\sigma=(3,0,2,2,0)$,
$\tilde\sigma w=(3,5,2,1,4,6)$
et $M\tilde\sigma=(3,3,0,2,2,0)$. La construction d'une
bijection reliant~$E$ et~$M$ est l'objet des sections
suivantes.

\sectiona 3. La transformation fondamentale|

\'Etant donn\'ee une permutation $\sigma\in {\goth S}_n$,
l'ensemble $\sigma^*(k)\!=\!\{\sigma^p(k):p\in {\bboard
N}\}$ est l'{\it orbite} contenant~$k$ $(k\in [\,n\,]$). On
note $z(\sigma)$ le nombre des orbites de~$\sigma$ ou,
de fa\c con \'equivalente, le nombre des {\it cycles}
de~$\sigma$. A chaque $k\in [\,n\,]$, nous attachons la
paire $\Pi_\sigma(k)=(\overline k,q_k)$, o\`u $\overline k$
est l'\'el\'ement maximum de l'orbite $\sigma^*(k)$ et o\`u
$q_k=\min\{p\in {\bboard N}:\sigma^p(k)=\overline k\}$.

\rem D\'efinition $1.8$|Pour $\sigma\in {\goth S}_n$ on
d\'esigne par~$\hat\sigma$ la permutation telle que pour
chaque~$k$ la paire $\Pi_\sigma(\hat\sigma(k))$ est le
$k$\ieme\ terme de la suite
$\bigl(\Pi_\sigma(j)\bigr)_{(j\in [\,n\,])}$ ordonn\'ee par
ordre lexicographique.

\medskip
Par exemple, en prenant encore $\sigma
w\!=\!(6,4,1,2,5,3)$, on a $z(\sigma)\!=\!3$, les trois
orbites \'etant $\{4,2\}$, le point fixe~$\{5\}$ et
$\{6,3,1\}$. Comme
$4=\sigma^0(4)=\sigma^1(2)$, puis $5=\sigma^0(5)$ et
enfin $6=\sigma^0(6)=\sigma^1(1)=\sigma^2(3)$, la suite
des $\Pi_\sigma(j)$, ordonn\'ee suivant l'ordre
lexicographique, est\hfil\break \null\kern.5cm $\bigl((4,0),
(4,1), (5,0), (6,0),(6,1), (6,2)\bigr)$, d'o\`u $\hat\sigma
w=(4,2,5,6,1,3)$.

\medskip
Rappelons qu'un \'el\'ement $x_k$ d'une suite
$(x_1,x_2,\ldots, x_n)\in {\bboard N}^n$ est dit {\it
saillant} si et seulement s'il n'existe aucun \'el\'ement
d'indice inf\'erieur qui soit sup\'erieur ou \'egal \`a lui,
c'est-\`a-dire si l'on a $x_{k'}<x_k$ pour tout $k'<k$.
Donc par d\'efinition $x_1$ est toujours saillant.

\th Lemme 1.9|L'\'el\'ement $k$ de $[\,n\,]$ est maximum dans
son orbite $\sigma^*(k)$ $($c'est-\`a-dire
$\Pi_\sigma(k)=(k,0))$ si et seulement si~$k$ est
saillant dans $\hat\sigma w$. De plus,~$k$ est un point
fixe de~$\sigma$ si et seulement si $k=\hat\sigma(j)$
avec soit $j=n$, soit $j<n$ et $\sigma(j+1)$ un autre
\'el\'ement saillant de~$\hat\sigma w$.
\finth

\dem
Soit $\Pi_\sigma(k)=(\overline k,q)$. Si et seulement si
$q$ est diff\'erent de~0, l'entier~$k$ n'est pas l'\'el\'ement
maximum~$\overline k$ de son orbite et~$k$ n'est pas
saillant dans~$\hat\sigma w$ puisque
$\Pi_\sigma(\overline k)=(\overline k,0)$ pr\'ec\`ede
$\Pi_\sigma(k)$ dans l'ordre lexicographique,
donc~$\overline k$ pr\'ec\`ede~$k$ dans la
suite~$\hat\sigma w$.

R\'eciproquement, soient $q=0$ et $\hat\sigma(j)=k$. Ou
bien on a $j=1$ auquel cas~$k$ est saillant, ou bien
$j\ge 2$, auquel cas pour tout $j'<j$, l'\'el\'ement
$\Pi_\sigma(\hat\sigma(j'))$ est avant l'\'el\'ement
$\Pi_\sigma(\hat\sigma(j))$ ($=(k,0)$) pour l'ordre
lexicographique ; ce qui implique, en posant
$\hat\sigma(j')=k'$, que l'on a
$\hat\sigma(j')=k'<\overline {k'}<k=\hat\sigma(j)$ et
prouve l'\'equivalence entre les \'el\'ements maximaux des
orbites de~$\sigma$ et les \'el\'ements saillants
de~$\hat\sigma w$.

Supposons maintenant $k=\overline k=\hat\sigma(j)$. Si
et seulement si~$k$ n'est pas un point fixe de~$\sigma$,
il est imm\'ediatement suivi dans~$\hat\sigma w$ d'un
\'el\'ement~$k'$ tel que $\Pi_\sigma(k')=(k,1)$. Par
cons\'equent, on a $j<n$ et $\hat\sigma(j+1)=k'$ n'est pas
saillant.\cqfd

\medskip
Dans ce qui suit ${\goth S}_n'$ d\'esigne l'ensemble des
$\sigma\in {\goth S}_n$ telles que $\sigma(1)=n$
et~${\cal C}_n$ est le sous-ensemble des permutations
circulaires.

\th Propri\'et\'e 1.10|L'application $\sigma\mapsto
\hat\sigma$ est bijective et satisfait \`a
$\sigma(n)=\hat\sigma(n)$. Sa restriction \`a~${\cal C}_n$
est une bijection sur ${\goth S}_n'$.
\finth

\dem
\'Etant donn\'e $\tau\in {\goth S}_n$, il existe un et un
seul $\sigma\in {\goth S}_n$ tel que $\tau=\hat\sigma$,
car d'une part, les \'el\'ements saillants de $\tau w$ livrent
les \'el\'ements maximaux des orbites de~$\sigma$ d'apr\`es le
pr\'ec\'edent lemme, d'autre part, les restrictions
de~$\sigma$ \`a chacune de ses orbites sont d\'etermin\'ees
par la succession m\^eme des \'el\'ements de~$\tau w$. Ceci
\'etablit le caract\`ere bijectif de $\sigma\mapsto \hat
\sigma$. Comme $\tau(1)$ et~$n$ sont toujours des
\'el\'ements saillants de~$\tau w$, il est clair que
$z(\sigma)=1$ si et seulement si $\tau(1)=n$. Enfin, si
l'on a $k=\sigma(n)$, c'est-\`a-dire $n=\sigma^{-1}(k)$,
l'entier~$n$ est l'\'el\'ement maximum de~$\sigma^*(k)$ et
par cons\'equent, $\Pi_\sigma(k)=(n,q_k)$ est le dernier
terme de la suite
$(\Pi_\sigma(j))_{(j\in [\,n\,])}$ ordonn\'ee par ordre
lexicographique, c'est-\`a-dire $k=\hat\sigma(n)$.\cqfd

\sectiona 4. Relations entre les exc\'edances et les
descentes| Nous allons donner une premi\`ere application de
la transformation fondamentale
$\sigma\mapsto\hat\sigma$. Aux 1-exc\'edances
de~$\sigma$ correspondent les 1-descentes
de~$\hat\sigma$, mais les 0-descentes de~$\hat\sigma$
ne correspondent pas aux 0-exc\'edances de~$\sigma$.
Ceci am\`ene \`a introduire un vecteur $D'\hat \sigma$ tel
que $D\hat\sigma+D'\hat\sigma$ ($=(D+D')\hat\sigma$)
rende compte \`a la fois des descentes et des \'el\'ements
saillants dans la suite $\hat\sigma w$ et qu'en outre
l'identit\'e 
$E\sigma=(D+D')\hat\sigma$ soit v\'erifi\'ee.

\rem D\'efinition $1.11$|Pour $\tau\in {\goth S}_n$, le
vecteur $D'\tau$ est le $n$-uple
$$\bigl(D'\tau(1),\ldots, D'\tau(n)\bigr)\in {\bboard
N}^n,$$
o\`u pour $j\in [\,n\,]$, on pose $D'\tau(j)=+1$ si d'une
part,~$j$ est saillant et d'autre part, soit $j=n$, soit
$j\le n-1$ et $\tau(1+\tau^{-1}(j))$ est aussi saillant ;
$D'\tau(j)=0$ dans tous les autres cas.

\th Th\'eor\`eme 1.12|On a identiquement
$$E\sigma=(D+D')\hat\sigma\quad
{\it et}\quad\Delta E\sigma=\Delta D\hat\sigma.$$
\finth

\dem
Quelque soit $\tau\in {\goth S}_n$, remarquons d'abord
que $j\le n-1$ est saillant dans $\tau w$ seulement si
$\tau^{-1}(j)\le n-1$. Nous v\'erifions le r\'esultat pour
chaque
$k=\hat\sigma(j)\in [\,n\,]$ en posant
$\Pi_\sigma(k)=(\overline k,q)$ et en distinguant les cas
suivants :

\goodbreak
\decale (i)|$q\ge 1$ ; on a alors $\hat\sigma(j-1)=k'$ o\`u
$\Pi_\sigma(k')=(\overline k,q-1)$, c'est-\`a-dire
$k'=\sigma^{-(q-1)}(\overline
k)=\sigma^{-(q-1)}(\sigma^q(k))=\sigma(k)$. Donc par
d\'efinition, on a
$D\hat\sigma(k)
=(k'-(k-1))_+=(\sigma(k)-(k-1))_+=E\sigma(k)$ avec 
$D\hat\sigma(k)=(D+D')\hat\sigma(k)$ puisque~$k$ n'est
pas saillant.

\decale (ii)|$q=0$, c'est-\`a-dire $k=\overline k$. On a
toujours $D\hat\sigma(k)=0$ car soit $j-1=0$, soit
$j-1\ge 1$ avec $k'=\hat\sigma(j-1)$ appartenant \`a une
orbite dont l'\'el\'ement maximum est strictement plus petit
que $k=\overline k$. D'autre part, $E\sigma(k)=0$ ou~1
selon que $\sigma(k)\not=k$ ou non, car
$E\sigma(k)=(\sigma(k)-(k-1))_+$ o\`u $\sigma(k)\le k$
puisque $k=\overline k$. D'apr\`es la deuxi\`eme partie du
lemme~1.9 et la d\'efinition de $D+D'$, on a donc
$E\sigma(k)=1$ si et seulement si
$1=(D+D')\hat\sigma(k)\not=D\hat\sigma(k)=0$.\cqfd

\medskip
On remarquera que l'on a $(D+D')\tau(j)\not=D\tau(j)$ si
et seulement si $D\tau(j)=0$ et $(D+D')\tau(j)=1$.

\rem Exemple|On a vu dans les exemples
pr\'ec\'edents que si
$\sigma w=(6,4,1,2,5,3)$, on avait $E\sigma=(6,3,0,0,1,0)$
et $\hat\sigma w=(4,2,5,6,1,3)$. Comme les \'el\'ements
successifs $\hat\sigma(3)=5$ et $\hat\sigma(4)=6$ sont
saillants dans $\hat\sigma w$, on a
$(D+D')\hat\sigma(5)=1$. Il vient donc
$(D+D')\hat\sigma=(6,3,0,0,1,0)=E\sigma$ et aussi $\Delta
D\hat\sigma=(5,2,0,0,0)=\Delta E\sigma$.

\sectiona 5. Applications aux permutations altern\'ees|
Une permutation $\sigma\in {\goth S}_n$ est {\it
altern\'ee} si l'on a $\sigma(2j)<\sigma(2j-1)$,
$\sigma(2j+1)$ pour tout entier~$j$ tel que $2\le 2j\le
n-1$ et si l'on a encore $\sigma(n)<\sigma(n-1)$
lorsque~$n$ est pair. D'autre part, elle est {\it
biexc\'ed\'ee} si pour tout $j\in [\,n\,]$, on a $j<\sigma(j)$,
$\sigma^{-1}(j)$ ou $j>\sigma(j)$, $\sigma^{-1}(j)$.
Les ensembles des permutations {\it altern\'ees} et {\it
biexc\'ed\'ees} sont not\'es respectivement ${\cal T}_n$ et
${\cal B}_n$. Nous allons, dans cette section, utiliser la
transformation fondamentale pour construire entre ${\cal
T}_n$ et ${\cal B}_n$ une bijection qui servira au
chapitre~V.

\th Lemme 1.13|Soit $k=\hat\sigma(j)$ $(j\in [\,n\,])$ ; on a
$k<\sigma(k)$, $\sigma^{-1}(k)$ si et seulement si
$j\not=1$ et soit $j\le n-1$ et
$\hat\sigma(j)<\hat\sigma(j-1)$, $\hat\sigma(j+1)$,
soit $j=n$ et $\hat\sigma(j)<\hat\sigma(j-1)$.
\finth

\dem
L'hypoth\`ese $k<\sigma(k)$ entra\^\i ne $k\not=\overline k$
o\`u $\overline k$ est le maximum de l'orbite $\sigma^*(k)$.
Donc $j\not=1$ et $\hat\sigma(j-1)=\sigma(k)$
entra\^\i nant l'\'equivalence de $k<\sigma(k)$ et de
$\hat\sigma(j)<\hat\sigma(j-1)$, puisque l'on a
$\hat\sigma(j)>\hat\sigma(j-1)$ quand $k=\overline k$.
Distinguons deux cas selon que $j=n$ ou non.

Si $j=n$, on a $\sigma^{-1}(k)=\overline k=n$, donc
$k<\sigma^{-1}(k)$ et le r\'esultat est prouv\'e. Si
$j\not=n$, ou bien $\sigma^{-1}(k)=\overline k$, auquel
cas l'hypoth\`ese $k<\sigma^{-1}(k)$ \'equivaut \`a
$k=\hat\sigma(j)<\hat\sigma(j+1)$, puisque
$\hat\sigma(j+1)$ est le maximum d'une autre orbite ; ou
bien $\sigma^{-1}(k)\not=\overline k$, auquel cas
$\hat\sigma(j+1)=\sigma^{-1}(k)$ et
$k<\sigma^{-1}(k)$ \'equivaut \`a
$\hat\sigma(j)<\hat\sigma(j+1)$.\cqfd

\th Proposition 1.14|Toute permutation $\sigma\in {\cal
B}_n$ a tous ses cycles de longueur paire ; donc ${\cal
B}_n=\emptyset$ si~$n$ est impair. La transformation
fondamentale $\sigma\mapsto \hat\sigma$ \'etablit,
lorsque~$n$ est pair, une bijection de~${\cal B}_n$
sur~${\cal T}_n$.
\finth

\dem 
Dire que $\sigma$ est biexc\'ed\'ee \'equivaut \`a dire que dans
toute orbite de~$\sigma$, d'\'el\'ement maximum~$\overline
k$, on a identiquement $\sigma^{2p+1}(\overline
k)<\sigma^{2p}(\overline k)$,
$\sigma^{2p+2}(\overline k)$ $(0\le p)$. Ces in\'egalit\'es
impliquent d'abord que l'orbite n'est pas r\'eduite \`a un seul
\'el\'ement ; supposons maintenant qu'elle ait un nombre
impair d'\'el\'ements, disons $2q+1$ (avec $q\ge 1$). Ces
m\^emes in\'egalit\'es appliqu\'ees \`a $p=q$ entra\^\i nent
que l'on a
$\overline k=\sigma^{2q+1}(\overline
k)<\sigma^{2q}(\overline k)$, ce qui est impossible
puisque~$\overline k$ est l'\'el\'ement maximum dans son
orbite. La permutation~$\sigma$ n'a donc que des cycles
de longueur paire et par cons\'equent ${\cal
B}_n=\emptyset$ si~$n$ est impair.

La derni\`ere partie de la proposition r\'esulte du
lemme~1.13 et des deux \'equivalences suivantes, qui
d\'ecoulent imm\'ediatement de ce qui pr\'ec\`ede :

\decale (i)|$\sigma$ est dans ${\cal B}_{2p}$ si et
seulement si les in\'egalit\'es $k<\sigma(k)$,
$\sigma^{-1}(k)$ sont satisfaites pour exactement~$p$
indices~$k$ ;

\decale (ii)|$\tau$ est dans ${\cal T}_{2p}$ si et
seulement si, en posant $\tau(2p+1)=2p+1$, les in\'egalit\'es
$\tau(j)<\tau(j-1)$, $\tau(j+1)$ sont vraies pour
exactement~$p$ indices $j\ge 2$.\cqfd

\rem Exemple|Nous donnons ci-dessous le tableau des cinq
permutations biexc\'ed\'ees de ${\goth S}_4$ et en face de
chacune d'elles, la permutation altern\'ee qui lui correspond
par l'application fondamentale.

$$
\vbox{\halign{\vrule\thinspace \hfil$#$\ \vrule
&\strut\quad\hfil$#$\hfil
&\quad\hfil$#$\hfil
&\quad\hfil$#$\hfil
&\quad\hfil$#$\hfil\ \vrule
&\quad\hfil$#$\hfil
&\quad\hfil$#$\hfil
&\quad\hfil$#$\hfil
&\quad\hfil$#$\hfil\ \vrule
&\ \hfil$#$\hfil\  \vrule\cr
\noalign{\hrule}
i={}\ &1&2&3&4&1&2&3&4&=i\cr
\noalign{\hrule}
\sigma(i)=&2&1&4&3&2&1&4&3&=\hat\sigma(i)\cr
&3&4&1&2&3&1&4&2&\cr
{\cal B}_4&4&3&2&1&3&2&4&1&\kern-15pt{\cal T}_4\cr
&4&3&1&2&4&1&3&2&\cr
&3&4&2&1&4&2&3&1&\cr
\noalign{\hrule}
}}
$$

\sectiona 4. Relations entre les exc\'edances et les mont\'ees|
Toujours au moyen de la transformation fondamentale,
nous construisons une bijection $\sigma\mapsto
\overline\sigma$ telle que $E\sigma=M\overline\sigma$.
Soit en effet $\sigma\in {\goth S}_n$ et notons~$\zeta$
la permutation circulaire d\'efinie dans le lemme~1.4
($\zeta w=(2,3,\ldots, n,1)$). On pose successivement
$$
\leqalignno{\noalign{\vskip-8pt}
\sigma_1&=\sigma\,\zeta;\cr
\sigma_2&=\hat\sigma_1\quad\hbox{(transformation
fondamentale)};&\hbox{puis}\cr
\overline\sigma&=\tilde\sigma_2\quad
\hbox{(o\`u\kern5pt $\tilde{}$\kern5pt est d\'efini dans le
lemme 1.7)}.&\hbox{et enfin}\cr}
$$

\th Th\'eor\`eme 1.15|L'application $\sigma\mapsto
\overline\sigma$ est une bijection de ${\goth S}_n$ sur
lui-m\^eme satisfaisant \`a $E\sigma=M\overline \sigma$.
\finth

\dem
Le caract\`ere bijectif est \'evident d'apr\`es ce qui pr\'ec\`ede.
Ensuite, on a $E\sigma(1)=\sigma(1)=\sigma_1(n)$ et on
v\'erifie que $E\sigma(j+1)=\Delta E\sigma_1(j)$ pour
chaque $j\in [n-1]$. D'apr\`es la propri\'et\'e 1.10, on a
$\sigma_2(n)=\sigma_1(n)$ et d'apr\`es le th\'eor\`eme~1.12, il
vient $\Delta E\sigma_1=\Delta D\sigma_2$. Par
cons\'equent, en utilisant le lemme~1.7, on a bien 
$E\sigma=M\overline \sigma$.\cqfd

\medskip
Par exemple, partant de $\sigma w=(6,4,1,2,5,3)$, on a
$E\sigma=(6,3,0,0,1,0)$, puis
$\sigma_1w\!=\!(4,1,2,5,3,6)$ et
$\sigma_2 w\!=\!(5,4,1,2,3,6)$ ; enfin $\overline \sigma
w\!=\!(6,3,2,1,4,5)$. Comme on a $\Delta E\sigma_1=\Delta
D\sigma_2=(3,0,0,1,0)$, on v\'erifie bien que $M\overline
\sigma=E\sigma$.

\medskip
On a not\'e que $E\sigma(k)=1$ si et seulement si
$\sigma(k)=k$. Par cons\'equent, la restriction de
$\sigma\mapsto \overline \sigma$ au 
sous-ensemble~${\cal D}_n$ des permutations 
{\it sans points fixes} est une bijection sur le sous-ensemble
des $\overline \sigma\in {\goth S}_n$ telles que
$M\overline \sigma(j)\not=1$, c'est-\`a-dire
des~$\overline\sigma$ telles que $\overline \sigma(1)\not=1$
et $1+\overline\sigma(j)\not=\overline\sigma(j+1)$, pour
chaque $j\in [n-1]$. L'ensemble de ces
permutations~$\overline\sigma$ n'est autre que la classe,
que nous noterons~${\cal G}_n$, des permutations {\it
sans successions}. On a ainsi le corollaire suivant.

\th Corollaire 1.16|La restriction de $\sigma\mapsto
\overline \sigma$ \`a l'ensemble~${\cal D}_n$ des
permutations sans points fixes est une bijection sur
l'ensemble~${\cal G}_n$ des permutations sans
successions telle que $E\sigma=M\overline\sigma$.
\finth

\sectiona 7. Relations avec les permutations circulaires|
Pour terminer ce chapitre, il nous reste \`a \'etudier la
distribution des vecteurs-exc\'edances sur l'ensemble~${\cal
C}_n$. Combinant d'abord les r\'esultats de la
propri\'et\'e~1.10 et du th\'eor\`eme~1.12, nous avons d\'ej\`a la
propri\'et\'e suivante.

\th Propri\'et\'e 1.17|La restriction de la transformation
fondamentale $\sigma\mapsto\hat\sigma$ \`a
l'ensemble~${\cal C}_n$ est une bijection sur ${\goth
S}_n'$ telle que $\Delta E\sigma=\Delta D\hat\sigma$.
\finth

Nous construisons ensuite une bijection $\sigma\mapsto
\sigma'$ de 
$${\goth S}_n''=\{\sigma\in {\goth
S}_n:\sigma(n)=1\}$$ 
sur ${\goth S}_n'$ telle que $\Delta E\sigma=\Delta
D\sigma'$. Si~$\sigma$ est dans ${\goth S}_n''$, on pose
$i=\hat\sigma{}^{-1}(n)$ et~$\sigma'$ est d\'efini comme
l'unique permutation telle que
$$
\sigma'w=\bigl(\hat\sigma(i),\hat\sigma(i+1),\ldots,
\hat\sigma(n),\hat\sigma(1),\ldots,\hat\sigma(i-1)\bigr).
$$

\th Lemme 1.18|L'application $\sigma\mapsto \sigma'$
est une bijection de ${\goth S}_n''$ sur~${\goth S}_n'$
satisfaisant \`a $\Delta E\sigma=\Delta D\sigma'$.
\finth

\dem
D'apr\`es la propri\'et\'e 1.10, on a $\sigma(n)=\hat\sigma(n)$
et par cons\'equent~$\hat\sigma$ est dans ${\goth S}_n''$.
Il est clair que $\sigma\mapsto \sigma'$ est bijectif. En
outre, $\Delta E\sigma=\Delta D\hat\sigma$ d'apr\`es le
th\'eor\`eme~1.12. Il suffit donc de v\'erifier $\Delta
D\sigma'(k)=\Delta D\hat \sigma(k)$ pour chaque $k\in
[\,n\,]$.

Distinguons d'abord deux cas particuliers :

\decale (i)|$k=\hat\sigma(1)$. On a $k=\sigma'(n-i+2)$
avec $\sigma'(n-i+1)=\hat\sigma(n)=1$. Donc $\Delta
D\hat\sigma(k)=(0-k)_+$ et $\Delta
D\sigma'(k)=(1-k)_+$ sont nuls.

\decale (ii)|$k=\hat\sigma(i)=n$. On a encore $\Delta
D\hat\sigma(k)=(\hat\sigma(i-1)-n)_+=0$ et aussi,
puisque $n=\sigma'(1)$, l'\'egalit\'e $\Delta
D\sigma'(k)=(0-n)_+=0$.

Dans le cas g\'en\'eral o\`u $k=\hat\sigma(j+1)\not=n$ $(j\in
[n-1]$), on a $k=\sigma'(j'+1)$ avec $j'=j-i+1$ ou
$n+j-i+1$ selon que $j\ge 1$ ou non. Dans les deux cas,
on a $\hat\sigma(j)=\sigma'(j')$ et par cons\'equent,
$(\hat\sigma(j)-k)_+$ est la valeur commune de $\Delta
D\hat\sigma(k)$ et $\Delta D\sigma'(k)$.\cqfd

\goodbreak
\medskip
Ce lemme permet la construction suivante : soit
$\sigma\in {\goth S}_{n-1}$ ; on d\'efinit $\sigma_1\in
{\goth S}_n''$ en posant
$$
\eqalignno{
\sigma_1(n)&=1\ {\rm et}\ 
\sigma_1(k)=1+\sigma(k)\quad
\hbox{pour chaque $k\in [n-1]$ ; puis}\cr
\sigma_2&=\sigma'_1,\ 
\hbox{o\`u $\tau\mapsto\tau'$ est la bijection d\'efinie
dans le lemme 1.18 ; enfin}\cr
\sigma''&\hbox{ est la permutation d\'efinie per
$\hat\sigma''=\sigma_2$}.\cr}
$$

\th Th\'eor\`eme 1.19|L'application $\sigma\mapsto \sigma''$
est une bijection de ${\goth S}_{n-1}$ sur~${\cal C}_n$
telle que $E\sigma=\Delta E\sigma''$.
\finth

\dem
Tout d'abord la bijectivit\'e de $\sigma\mapsto \sigma''$
est \'evidente. Si~$\sigma$ est dans ${\goth S}_{n-1}$,
on a ensuite $\sigma_1\in {\goth S}_n''$ et
$E\sigma(k)=(\sigma(k)-(k-1))_+
=(\sigma_1(k)-k)_+=\Delta E\sigma_1(k)$ 
pour chaque $k\in [n-1]$, d'o\`u il r\'esulte
$E\sigma=\Delta E\sigma_1$. D'autre part, d'apr\`es le
lemme 1.18, on a $\sigma_2\in {\goth S}'_n$ et
$E\sigma=\Delta E\sigma_1=\Delta D\sigma_2$. Enfin, en
vertu de $\sigma_2\in {\goth S}_n'$, la propri\'et\'e~1.17
montre que l'on a $\sigma''\in {\cal C}_n$ et $\Delta
E\sigma''=\Delta D\sigma_2$.\cqfd

\sectiona 8. Tableau des bijection utilis\'ees|
Il para\^\i t int\'eressant de rappeler les propri\'et\'es des
bijections construites dans le premier chapitre et
d'indiquer leur r\'ef\'erence.

$$
\vbox{\offinterlineskip\halign{\vrule 
\hfil\thinspace #\thinspace\hfil\vrule
&\vrule height10pt depth 6pt width 0pt\hfil\thinspace
#\thinspace\hfil\vrule &\hfil\thinspace
#\thinspace\hfil\vrule &\hfil\thinspace
#\thinspace\hfil\vrule &\hfil\thinspace
#\thinspace\hfil\vrule\cr
\noalign{\hrule}
La bijection&envoie&sur&a la propri\'et\'e& et la r\'ef\'erence\cr
\noalign{\hrule}
$\sigma\mapsto \sigma\zeta^r$&${\goth S}_n$&
${\goth S}_n$&$\Delta'^rE\sigma=\Delta^rE\sigma \zeta^r
\ (r\ge 0)$&Lemme 1.4\cr
\noalign{\hrule}
$\sigma\mapsto \tilde\sigma$&${\goth S}_n$
&${\goth S}_n$&$M\tilde\sigma(1)=\sigma(n)$ et&Lemme
1.7\cr
&&&$M\tilde\sigma(k+1)\!=\!\Delta D\sigma(k)\;(k\!\in\!
[n-1])$&\cr
\noalign{\hrule} 
$\sigma\mapsto \hat\sigma$&${\goth S}_n$&
${\goth S}_n$&&D\'efinition 1.8\cr
&${\cal C}_n$&${\goth S}_n'$&&Proposition 1.10\cr
&${\goth S}_n$&
${\goth S}_n$&$E\sigma=(D+D')\hat\sigma$,
$\Delta E\sigma=\Delta D\hat\sigma$&Th\'eor\`eme
1.12\cr 
&${\cal B}_{2n}$&${\cal T}_{2n}$&&Proposition 1.14\cr
\noalign{\hrule} 
$\sigma\mapsto \overline \sigma$&${\goth S}_n$&
${\goth S}_n$&$E\sigma=M\overline \sigma$&Th\'eor\`eme
1.15\cr 
&${\cal D}_n$&${\cal G}_n$&&Corollaire 1.16\cr
\noalign{\hrule}
$\sigma\mapsto \sigma'$&${\goth S}_n''$&
${\goth S}_n'$&$\Delta E\sigma=\Delta D\sigma'$&Lemme
1.18\cr 
\noalign{\hrule}
$\sigma\mapsto \sigma''$&${\goth S}_{n-1}$&${\cal
C}_n$&$E\sigma=\Delta E\sigma''$&Th\'eor\`eme 1.19\cr
\noalign{\hrule}}}$$

\sectiona 9. Notations g\'en\'erales|
Nous r\'eunissons dans cette section toutes les notations
utilis\'ees pour les sur- et les sous-ensembles de~${\goth
S}_n$. Pour $n\ge 1$ et $1\le k,r\le n$, on consid\`ere les
sous-ensembles suivants de~${\goth S}_n$ :

\medskip
\itemitem{${\cal C}_n$}l'ensemble des permutations {\it
circulaires} ; (${\cal C}_0=\emptyset$) ;
\itemitem{${\cal G}_n$}l'ensemble des permutations {\it
sans successions}, c'est-\`a-dire des $\sigma \in {\goth
S}_n$ telles que $1\not=\sigma(1)$ et
$1+\sigma(j)\not=\sigma(j+1)$ pour $j\in [n-1]$ ;
\itemitem{${\cal D}_n$}l'ensemble des permutations {\it
sans points fixes} ;
\itemitem{${\cal T}_n$}l'ensemble des {\it permutations
altern\'ees}, c'est-\`a-dire des $\sigma\in {\goth S}_n$
telles que $\sigma(2j)<\sigma(2j-1)$, $\sigma(2j+1)$
pour tout entier~$j$ tel que $2\le 2j\le n-1$ et en plus
si~$n$ est pair, telles que $\sigma(n)<\sigma(n-1)$ ;
\itemitem{${\cal B}_n$}l'ensemble des permutations {\it
biexc\'ed\'ees}, c'est-\`a-dire des $\sigma\in {\goth S}_n$
telles que pour chaque $j\in [\,n\,]$, on ait $j<\sigma(j)$,
$\sigma^{-1}(j)$ ou $j>\sigma(j)$, $\sigma^{-1}(j)$ ;
\itemitem{${\goth S}_{n,k}$}l'ensemble des permutations
$\sigma\in {\goth S}_n$ telles que $\sigma(1)=k$ ; 
${\goth S}_n'={\goth S}_{n,n}$ ;
\itemitem{${\goth S}''_n$}l'ensemble des permutations
$\sigma\in {\goth S}_n$ telles que $\sigma(n)=1$ ;
\itemitem{${}_r{\goth S}_{n}$}l'ensemble des
permutations
$\sigma\in {\goth S}_n$ telles que\hfil\break
$\sigma^{-1}(n-r+1)<\sigma^{-1}(n-r+2)<\cdots<
\sigma^{-1}(n)$.

\medskip
On posera \'egalement :\quad
${\goth S}=\kern-5pt\bigcup\limits_{0\le n}
{\goth S}_n$ ; puis 
${\cal C}=\kern-5pt\bigcup\limits_{1\le n}{\cal C}_n$,\quad 
${\cal G}=\kern-5pt\bigcup\limits_{1\le n}{\cal G}_n$,\quad
${\cal T}=\kern-5pt\bigcup\limits_{1\le n}{\cal T}_n$,\quad 
${\cal B}=\kern-5pt\bigcup\limits_{1\le n}{\cal B}_n$\quad et
\quad ${\cal D}=\kern-5pt\bigcup\limits_{1\le n}{\cal D}_n$.
Enfin, on utilise les notations courantes suivantes :

\itemitem{$\bboard N$}l'ensemble des entiers naturels ;
\itemitem{$\bboard Z$}l'ensemble des entiers ;
\itemitem{$\bboard Q$}l'ensemble des nombres rationnels.

\vfill\eject

\pagetitretrue

\auteurcourant={{\eightrm CHAPITRE II :
LES POLYN\^OMES EUL\'ERIENS}}

\vglue 2cm
\centerline{{\eightrm CHAPITRE II}}
\vskip 2mm
\centerline{{\bf LES POLYN\^OMES EUL\'ERIENS}} 

\vskip 6mm plus 2mm
\sectiona 1. Interpr\'etation des polyn\^omes eul\'eriens|
Les th\'eor\`emes 1.12, 1.15, 1.19 et la propri\'et\'e 1.17
permettent d'\'etablir imm\'ediatement l'\'egalit\'e des cinq
ensembles pond\'er\'es
$E{\goth S}_n$, $(D+D'){\goth S}_n$, $M{\goth S}_n$,
$\Delta E {\cal C}_{n+1}$, $\Delta D {\goth S}'_{n+1}$ ;
et le th\'eor\`eme 1.12 donne encore : 
$\Delta E {\goth S}_n= \Delta D {\goth S}_n$.
La propri\'et\'e 2.1 suivante r\'esulte alors du th\'eor\`eme 1.5.

\th Propri\'et\'e 2.1|Soit $\Gamma$ un mon\^ome en~$\Delta$
et~$\Delta'$ de degr\'e $r\ge 0$. On a
$$\displaylines{
\Gamma E{\goth S}_n=\Gamma (D+D'){\goth S}_n
=\Gamma M {\goth S}_n=\Gamma \Delta E {\cal C}_{n+1}
=\Gamma \Delta D {\goth S}_{n+1}',\cr
\noalign{\hbox{o\`u en outre}}
\Gamma E{\goth S}_n=\Gamma D {\goth S}_n\cr}
$$
si et seulement si $\Gamma$ a au moins un $\Delta$
comme facteur.
\finth

Nous notons $|x|$ le nombre de termes positifs de tout 
$x\in {\bboard N}^p$ et introduisant une
ind\'etermin\'ee~$t$, nous posons $\theta x=t^{|x|}$.
Si~$K$ est une application dans ${\bboard N}^p$ d'une
partie~$\cal P$ de ${\goth S}_n$, l'ensemble pond\'er\'e
$$\theta K{\cal P}=\sum\{\theta K\sigma:\sigma\in {\cal
P}\}=\sum_{0\le k} t^k \cdot \card \{\sigma\in {\cal
P}:|K\sigma|=k\}$$
sera appel\'e, par abus de langage, {\it fonction
g\'en\'eratrice} de~$K$. Pour~${\cal P}$ fini, 
$\theta K{\cal P}$ sera donc un polyn\^ome en~$t$ \`a
coefficients dans~${\bboard N}$ et nous dirons que
$({\cal P},K)$ est une {\it interpr\'etation}.

\def \rA{{}^{r\kern-3pt}A}
Les polyn\^omes eul\'eriens $A_n(t)$ et leurs g\'en\'eralisations
$\rA_n(t)$ selon Riordan sont d\'efinis par :
$$
\rA_n(t)=\theta \Delta^rE{\goth S}_n
\quad{\rm pour}\ 0\le r\le n-1.
$$
On conviendra que $\rA_n(t)=n!$ pour $r\ge n$. On
posera $A_0(t)=1$ et $A_n(t)={}^{1\kern-3pt}A_n(t)$
pour $n\ge 1$. Par construction, les $\rA_n(t)$ sont des
polyn\^omes de degr\'e au plus \'egal \`a $n-r$. L'\'enonc\'e
suivant en donne plusieurs interpr\'etations par simple
application de la propri\'et\'e pr\'ec\'edente.

\th Propri\'et\'e 2.2|Soit $\Gamma$ un mon\^ome en $\Delta$
et $\Delta'$ de degr\'e~$r$. On a
$$
\rA_n(t)=\theta\Gamma E{\goth S}_n=\theta\Gamma(D+D')
{\goth S}_n=\theta\Gamma M{\goth S}_n=\theta\Gamma
\Delta E {\cal C}_{n+1}
=\theta \Gamma \Delta D {\goth S}_{n+1}',
$$
o\`u en outre $\rA_n(t)=\theta\Gamma D{\goth S}_n$ si et
seulement si $\Gamma$ a au moins un~$\Delta$ comme
facteur.
\finth

\goodbreak
De m\^eme, d'apr\`es le corollaire 1.16, on a l'\'egalit\'e entre
les ensembles pond\'er\'es $E{\cal D}_n$ et $M{\cal G}_n$,
d'o\`u encore
$$
\theta E{\cal D}_n=\theta M{\cal G}_n\quad
{\rm pour}\ n\ge 1.
$$
La valeur commune de ces deux derniers polyn\^omes sera
d\'esign\'ee par~$B_n(t)$. L'interpr\'etation $({\cal G}_n,M)$
de ces polyn\^omes est due \`a Roselle [25].
L'interpr\'e\-ta\-tion
$({\cal D}_n, E)$ servira \`a \'etablir au chapitre~IV
l'expression de la fonction g\'en\'eratrice exponentielle
des~$B_n(t)$.

\sectiona 2. Propri\'et\'es de sym\'etrie|
Les identit\'es (1) et (4) ci-dessous sont bien connues (\cf.
Riordan [24]). Nous en donnons ici des d\'emonstrations
\'el\'ementaires. Les polyn\^omes $\rA_n(t)$ pour $r\ge 2$
n'ont pas de propri\'et\'e de sym\'etrie \'evidente. En revanche,
si l'on forme les polyn\^omes r\'eciproques
$t^{n-r}\,\rA_n(t^{-1})$, on obtient plusieurs
interpr\'etations (\cf. les relations~(2) et~(3) ci-dessous)
qui nous serviront effectivement, dans les sections~2.4
et~2.6, pour \'etablir des connexions avec le probl\`eme de
Simon Newcomb et pour d\'emontrer de nouvelles identit\'es
sur les polyn\^omes eul\'eriens.

\th Propri\'et\'e 2.3|Pour $n\ge 1$ on a :
$${}^{0\kern-3pt}A_n(t)=t\,A_n(t).\eqno(1)$$
De plus,
si~$\Gamma$ est un mon\^ome de la forme
$\Delta'{}'^{r-1}\Delta$  ou $\Delta''{}^{r-1}\Delta'$ $(r\ge
1)$, on~a
$$\eqalignno{t^{n-r}\,\rA_n(t^{-1})
=\theta\Gamma E{\goth S}_n&=\theta\Gamma(D+D'){\goth
S}_n
=\theta\Gamma M{\goth S}_n.&(2)\cr
\noalign{\hbox{On a encore}}
t^{n-r}\,\rA_n(t^{-1})&=\theta \Delta''{}^{r-1}\Delta
D{\goth S}_n,&(3)\cr
\noalign{\hbox{d'o\`u en particulier pour $r=1$}}
t^{n-1}\,A_n(t^{-1})&=A_n(t).&(4)\cr}
$$
\finth

\dem 
Soit $\sigma\in {\goth S}_n$ ; d\'efinissant
$\check \sigma$ par l'identit\'e 
$\check\sigma(k)=n+1-\sigma(n+1-k)$, l'on a $E\check
\sigma(k)=1$ si et seulement si $E\sigma(n+1-k)=1$ et
$E\check \sigma(k)\ge 2$ si et seulement si
$E\sigma(n+1-k)=0$. Dans ces conditions, il vient 
$|E\check \sigma|+|\Delta E\sigma|=n$, ce qui \'etablit
$$
{}^{0\kern-3pt}A_n(t)=t^n\,A_n(t^{-1}).\eqno(5)
$$
D'autre part, la relation entre les vecteurs
$E\check\sigma$ et $E \sigma$ peut encore s'exprimer
par la condition
$$
\Delta E\check\sigma(k)\ge 1\ 
\hbox{si et seulement si } 
\Delta'E\sigma(n-k)=0
\hbox{ pour } 1\le k\le n-1.\eqno(6)
$$
Mais la condition (6) est encore \'equivalente \`a
$$
\displaylines{
|\Delta'{}^{r-1}\Delta E\check \sigma|+
|\Delta''{}^{r-1}\Delta' E\sigma|=n-r\cr
\noalign{\hbox{ou encore \`a}}
|\Delta''{}^{r-1}\Delta E\check \sigma|+
|\Delta'{}^{r} E\sigma|=n-r.\cr
}$$

\goodbreak
\noindent
On en d\'eduit 
$$
\eqalignno{\theta \Delta''{}^{r-1}\Delta E{\goth S}_n
&=t^{n-r}\sum_{0\le k}t^{-k}\card\{\sigma\in {\goth
S}_n : |\Delta'{}^{r} E\sigma|=k\}\cr
&=t^{n-r}\,\rA_n(t^{-1}),\cr}
$$
d'apr\`es la propri\'et\'e 2.2.

Les relations (2) et (3) r\'esultent alors de la propri\'et\'e~2.1
et faisant $r=1$ dans~(3), on obtient 
$t^{n-1}\,A_n(t^{-1})=\theta\Delta D {\goth
S}_n=A_n(t)$, c'est-\`a-dire la relation~(4). Enfin,
l'identit\'e~(2) r\'esulte \`a la fois de~(4) et de~(5).\cqfd

\rem Remarque $2.4$|D'apr\`es la d\'efinition~1.6, pour $1\le
r\le n$ et $\sigma\in {\goth S}_n$, il y a exactement
$|\Delta'\Delta''{}^{r-1}M\sigma|$ indices~$i$ tels que
$1\le i\le n-1$, $\sigma(i)<\sigma(i+1)$ et
$\sigma(i)<n-r$. Comme on a pos\'e d'autre part
$$\displaylines{
\rA_n(t)=\sum_{0\le k\le n-r}\rA_{n,k}\,t^k\cr
\noalign{\hbox{et que d'apr\`es la propri\'et\'e 2.3, on a}}
\theta\Delta'\Delta''{}^{r-1}M{\goth S}_n
=t^{n-r}\,\rA_n(t^{-1})
=\sum_{0\le s\le n-r}\rA_{n,n-r-s}\,t^s,\cr
\noalign{\vskip-8pt}
\noalign{\hbox{il vient :}}
\rA_{n,n-r-s}=\card\{\sigma\in {\goth S}_n :
|\Delta'\Delta''{}^{r-1}M\sigma|=s\}
\quad{\rm pour}\ 0\le s\le n-r.\cr}
$$
Cette remarque nous servira au paragraphe~6 du pr\'esent
chapitre.

\sectiona 3. Relations de r\'ecurrence|
Dans cette section, nous \'etablissons une relation de
r\'ecurrence sur les polyn\^omes eul\'eriens qui g\'en\'eralise
l'identit\'e~(1) ci-dessus et red\'emontrons la relation de
r\'ecurrence trouv\'ee par Riordan [24].

\th Propri\'et\'e 2.5|Pour $0\le r\le n$, on a l'identit\'e :
$$
t\cdot {}^{(r+1)\kern-3pt}A_n(t)=\rA_n(t)+r(t-1)\cdot
\rA_{n-1}(t).
$$
\finth

\dem
Pour $r=0$, l'identit\'e se r\'eduit \`a
$tA_n(t)={}^{0\kern-3pt}A_n(t)$. Pour $r=n$, elle est
encore vraie avec la convention que nous avons faite que
$\rA_n(t)=n!$ quand $r\ge n$. Nous supposons donc $1\le
r\le n-1$. Pour chaque $\sigma\in {\goth S}_n$, on a
$$
\Delta^rE\sigma=\bigl((\sigma(1)-r)_+,
(\sigma(2)-r-1)_+, \ldots, (\sigma(n-r)-(n-1))_+\bigr).
$$
D'une part, on a $\Delta^rE{\goth
S}_{n,s}=\Delta^rE{\goth S}_{n,r}$ pour $1\le s\le r$ ;
d'autre part,
$r\theta\Delta^r E{\goth S}_{n,1}
=r\theta\Delta^r\Delta'E{\goth S}_{n,1}=
r\,\rA_{n-1}(t)$ d'apr\`es la propri\'et\'e~2.2. Il en r\'esulte :
$\rA_n(t)-r\,\rA_{n-1}(t)=
\theta\Delta^r({\goth S}_n-({\goth S}_{n,1}
+\cdots+{\goth S}_{n,r}))
=\theta\Delta^rE({\goth S}_{n,r+1}
+\cdots+{\goth S}_{n,n})
=t\Delta'\Delta^rE({\goth S}_{n,r+1}
+\cdots+{\goth S}_{n,n})
=t\Delta'\Delta^rE({\goth S}_n-({\goth S}_{n,1}
+\cdots+{\goth S}_{n,r}))
=t(\Delta'\Delta^rE{\goth S}_{n,1} 
-r\Delta'\Delta^rE{\goth S}_{n,1})
=t({}^{(r+1)\kern-3pt}A_n(t)
-r\,\rA_{n-1}(t))$.\cqfd

\rem Remarque $2.6$|Riordan ([24], p. 214) a trouv\'e une
autre relation de r\'ecurrence, \`a savoir
$$
\rA_n(t)=(r+(n-r)t)\cdot \rA_{n-1}(t)+t(1-t)\cdot
\rA_{n-1}'(t)\ (0\le r\le n;\,1\le n),\eqno(7)
$$

\goodbreak
\noindent
o\`u $\rA_{n-1}'(t)$ d\'esigne la d\'eriv\'ee du polyn\^ome
$\rA_{n-1}(t)$. Comme on a pos\'e
$$
\rA_n(t)=\sum_{0\le k\le n-r}\rA_{n,k}t^k,
$$
cette relation de r\'ecurrence est encore \'equivalente aux
$(n-r+1)$ relations suivantes (\cf. [24] p. 215)
$$
\rA_{n,k}=(k+r)\cdot \rA_{n-1,k}+(n+1-k-r)\cdot
\rA_{n-1,k-1}\quad (0\le k\le n-r),\eqno(8)
$$
o\`u l'on a pos\'e $\rA_{n,k}=0$ si $k\le -1$ ou $k\ge
n-r+1$.

\medskip
Comme l'a not\'e Welschinger [30], on peut red\'emontrer
facilement (8) et par suite~(7) en prenant les polyn\^omes
eul\'eriens dans l'interpr\'etation
$$
\rA_n(t)=\theta \Delta'{}^{r-1}\Delta D{\goth S}_n.
$$
En effet, on v\'erifie tout d'abord que les relations (8) sont
vraies pour $n=r$, en notant que
${}^{n\kern-3pt}A_{n,0}=n!$ et 
${}^{n\kern-3pt}A_{n,k}=0$ pour $k\not=0$. On suppose
ensuite $0\le r\le n-1$ et l'on pose pour $i\in[\,n\,]$ et
$\sigma\in {\goth S}_{n-1}$
$$\displaylines{
\eta_i(\sigma)=\bigl(\sigma(1),\ldots,
\sigma(i-1),n,\sigma(i),\ldots, \sigma(n-1)\bigr).\cr
\noalign{\hbox{Il est clair que l'on a}}
{\goth S}_{n}=\{\eta_i(\sigma):i\in [\,n\,],\,\sigma\in
{\goth S}_{n-1}\}.\cr}
$$
Soit $\sigma\in {\goth S}_{n-1}$ ; alors
$|\Delta'{}^{r-1}\Delta D\sigma|=k$ [en abr\'eg\'e :
$\sigma\in{}^r{\goth A}_{n-1,k}$] si et seulement s'il
existe~$k$ couples $(j-1,j)$ tels que $1\le j-1$ et
$\sigma(j-1)>\sigma(j)\ge r$.
Prenons $\sigma$ dans ${}^r{\goth A}_{n-1,k}$ ; on
observe alors que $\eta_i(\sigma)$ appartient \`a
${}^r{\goth A}_{n,k}$ pour les seuls indices~$i$ suivants

\decale (i)|$1\le i-1$ et $\sigma(i-1)>\sigma(i)\ge r$ ;
\decale (ii)|$i\in [n-1]$ et $\sigma(k)\le r-1$ ;
\decale (iii)|$i=n$ ;

\noindent
c'est-\`a-dire pour exactement $k+(r-1)+1=k+r$ indices
$i\in [\,n\,]$. Pour les autres $n-(k+r)$ indices~$i$ ne
satisfaisant \`a aucune des conditions~(i), (ii), (iii), on a
$\eta_i(\sigma)\in {}^r{\goth A}_{n,k+1}$.
On constate donc que l'ensemble ${}^r{\goth A}_{n,k}$ est
contenu dans la r\'eunion $\bigcup\{\eta_i({}^r{\goth
A}_{n-1,k}\cup {}^r{\goth A}_{n-1,k-1}):i\in [\,n\,]\}$.
On voit ensuite que la relation $\eta_i(\sigma)\in
{}^r{\goth A}_{n,k}$ est v\'erifi\'ee pour exactement $(k+r)$
indice~$i$ si~$\sigma$ est dans
${}^r{\goth A}_{n-1,k}$ et pour exactement
$n-(k-1+r)=n+1-k-r$ indices~$i$ si~$\sigma$ est
dans ${}^r{\goth A}_{n-1,k-1}$. Les relations~(8) sont
ainsi d\'emontr\'ees.

\sectiona 4. Relations avec le \og probl\`eme de Simon
Newcomb\fg|
Nous allons expliciter maintenant le lien entre les
polyn\^omes $\rA_n(t)$ et les polyn\^omes g\'en\'erateurs que
l'on d\'efinit pour le \og probl\`eme de Simon Newcomb
avec une sp\'ecification $(1^r(n-r))$\fg\ (voir MacMahon
[20], vol.~1, chap.~4 et~5). Dans la d\'emonstration de la
propri\'et\'e qui suit, nous consid\'erons $\sigma w$ comme le
mot $\sigma(1)\sigma(2)\ldots \sigma(n)$ dont les lettres
sont des \'el\'ements de~$[\,n\,]$.

\th Propri\'et\'e 2.7|Soit $r\ge 2$ ; on a :
$t^{n-r}\,\rA_n(t^{-1})=r!\,\theta\Delta D\,{}_r{\goth
S}_n$, o\`u
${}_r{\goth S}_n=\{\sigma\in {\goth S}_n:
\sigma^{-1}(n-r+1)<\sigma^{-1}(n-r+2)<\cdots<
\sigma^{-1}(n)\}$.
\finth

\goodbreak
\dem
D'apr\`es la propri\'et\'e 2.3, il nous suffit d'\'etablir
$\theta\Delta''{}^{r-1}\Delta D{\goth S}_n=r!\,\theta
\Delta D\,{}_r{\goth S}_n$ ou encore de trouver une
surjection $\sigma\mapsto \sigma'$ de~${\goth S}_n$ sur
${}_r{\goth S}_n$ telle que

\decale (i)|$|\Delta''{}^{r-1}\Delta D\sigma|
=|\Delta D\sigma'|$ ;

\decale (ii)|l'application $\sigma\mapsto \sigma'$ soit
homog\`ene de degr\'e~$r!$, c'est-\`a-dire que l'image inverse
de chaque $\sigma'\in {}_r{\goth S}_n$ ait~$r!$
\'el\'ements.

Prenons $\sigma\in {\goth S}_n$ et soit
$\sigma w\!=\!g_1i_{n-r+1}g_2\ldots
g_ri_ng_{r+1}$, o\`u $\{i_{n-r+1},\ldots,
i_n\}\!\!=\{n-r+1,\ldots, n\}$. On pose
$\sigma'w=g_1(n-r+1)g_2\ldots g_rng_{r+1}$. Il est
clair que~$\sigma'$ est dans ${}_r{\goth S}_n$ et que
$\sigma\mapsto \sigma'$ est homog\`ene de degr\'e~$r!$
D'autre part, puisque les \'el\'ements $n-r+1$, $n-r+2$,
\dots~, $n$ se pr\'esentent dans cet ordre dans le
mot~$\sigma'w$, on a $\Delta D\sigma(j)=0$ pour
$j=n-r+1,n-r+2,\ldots, n-1$, ou encore 
$|\Delta D\sigma'|=|\Delta''{}^{r-1}\Delta D\sigma'|$.
Enfin, l'identit\'e $|\Delta''{}^{r-1}\Delta
D\sigma|=|\Delta''{}^{r-1}\Delta D\sigma'|$ r\'esulte du
fait que l'on a $\sigma'(j)=\sigma(j)$ si $\sigma(j)\le
n-r$ et que $\sigma'(j)\ge n-r+1$ si et seulement si
$\sigma(j)\ge n-r+1$.\cqfd

\medskip
Si l'on envoie tout
$\sigma'w=g_1(n-r+1)g_2(n-r+2)\ldots g_rng_{r+1}$
de~${}_r{\goth S}_n$ sur le mot
$f=g_1(n-r+1)g_2(n-r+1)\ldots g_r(n-r+1)g_{r+1}$,
on d\'efinit une bijection de~${}_r{\goth S}_n$ sur une
classe de mots de sp\'ecification $(1^r(n-r))$, 
c'est-\`a-dire des mots de longueur~$n$, qui contiennent
$n-r+1$ lettres distinctes
dont l'une d'entre elles est r\'ep\'et\'ee~$r$ fois. Si l'on
d\'efinit maintenant $|\Delta Df|$ comme le {\it nombre de
descentes}, c'est-\`a-dire le nombre de couples de lettres
successives dans~$f$ qui vont en d\'ecroissant, on voit que
l'on a $|\Delta D f|=|\Delta D\sigma'|$. Ainsi
$t^{n-r}\, \rA_n(t^{-1})$ est le {\it polyn\^ome g\'en\'erateur
du nombre des descentes pour un ensemble de
sp\'ecification} $(1^r(n-r))$.

\sectiona 5. Relations avec les nombres de Stirling|
Rappelons que pour $1\le q\le p$, le {\it nombre de
Stirling de seconde esp\`ece} $S(p,q)$ est le nombre de
partitions de~$[\,p\,]$ en~$q$ parties non vides. Le r\'esultat
suivant est obtenu par Riordan ([24], p.~213) au moyen de
calculs assez complexes.

\th Propri\'et\'e 2.8|L'entier $S(p,q)$ est le nombre de
parties $W\subset [\,p\,]\times [\,p\,]$ de $p-q$ \'el\'ements qui
satisfont aux conditions suivantes :
\decale {\rm (i)}|$W$ est une quasi-permutation,
c'est-\`a-dire qu'il existe au moins un $\sigma\in {\goth
S}_p$ tel que $W\subset \{(k,\sigma(k))\in [\,p\,]\times
[\,p\,] : k\in [\,p\,]\}$ ;

\decale {\rm (ii)}|$W$ est supra-diagonale, c'est-\`a-dire
que $(k,k')\in W$ implique $k<k'$.
\finth

\dem
Soit $\{E_1,E_2,\ldots, E_q\}$ une partition de $[\,p\,]$
que nous pouvons consid\'erer comme form\'ee des classes
d'une relation d'\'equivalence $E\subset [\,p\,]\times
[\,p\,]$. A chaque $E_j=\{i_1<i_2<\cdots <i_{n_j}\}$
comprenant $n_j\ge 1$ \'el\'ements, nous associons la
quasi-permutation supra-diagonale $E'_j=\{(i_1,i_2),
(i_2,i_3),\ldots,(i_{n_j-1},i_{n_j})\}$ contenant $n_j-1$
(\'eventuellement z\'ero) \'el\'ements de $[\,p\,]\times [\,p\,]$.
Posant $E'=\bigcup\limits _{1\le j\le q}E'_j$, on voit
que~$E'$ est une quasi-permutation supra-diagonale
ayant $\sum (n_j-1)=p-q$ \'el\'ements et que~$E$ est la
plus fine des \'equivalences sur $[\,p\,]$ qui
contienne~$E'$.

R\'eciproquement, \'etant donn\'ee une quasi-application
supra-diagonale~$W$ ayant $p-q$ \'el\'ements, soit $E$ la
plus fine de toutes les \'equivalences sur $[\,p\,]$ qui
contienne~$W$. On a $W=E'$ et la bijection d\'esir\'ee est
\'etablie.\cqfd

\medskip
La formule de Riordan ([24], p. 214)
$$
\rA_n(s+1)=\sum_{0\le k\le
n-r}s^k\,(n-k)!\,\,S(n+1-r,n+1-r-k)
$$
\'etant obtenue par celui-ci au moyen de la m\'ethode
alg\'ebrique classique d'inversion de M\"obius, nous pensons
pouvoir nous dispenser d'en reproduire ici la
d\'emonstration.

\sectiona 6. Les identit\'es de Worpitzky|
Soient $m$, $n$, $s$ trois entiers tels que
$0\le s<n\le m+s$ et $(i_1,i_2,\ldots, i_s)$ une suite
strictement d\'ecroissante d'entiers compris entre~1 et
$n-1$ que nous nommerons {\it indices distingu\'es}. Nous
allons d'abord d\'enombrer l'ensemble $\Phi_{m,n,s}$ de
tous les morphismes $\phi:[\,n\,]\rightarrow [\,m\,]$ tels que
$\phi(i)\not=\phi(i+1)$ si~$i$ n'est pas un indice
distingu\'e. Le d\'enombrement de $\Phi_{m,n,s}$ permet non
seulement d'obtenir l'identit\'e~(9) ci-dessous due \`a
Worpitzky, mais une g\'en\'eralisation de celle-ci au cas des
polyn\^omes $\rA_n(t)$ $(r\ge 2)$. La proposition~2.9
ci-dessous est bien connue.

\th Proposition 2.9|On a :
$\card \Phi_{m,n,s}={m+s\choose n}$.
\finth

\dem
Il suffit de faire correspondre, de fa\c con bijective, \`a tout
$\phi\in \Phi_{m,n,s}$ un morphisme injectif
$\psi:[\,n\,]\rightarrow [m+s]$ (c'est-\`a-dire une
application strictement croissante). Dans ce but,
d\'esignons pour tout entier $k\in[\,n\,]$, par $\theta(k)$
le nombre d'indices distingu\'es avant~$k$, \`a savoir le
nombre d'indices~$j$ tels que $i_j<k$. On a
$\theta(1)=0$ et $\theta(n)=s$. La bijection
$\phi\mapsto \psi$ est alors d\'efinie de la fa\c con suivante.
Pour tout $k\in[\,n\,]$, on pose $\psi(k)=\phi(k)+\theta(k)$.
On a $1\le \psi(n)\le m+s$ et ~$\psi$ est strictement
croissante, car si $k-1$ est distingu\'e, on a
$\theta(k-1)<\theta(k)$ et si $k-1$ ne l'est pas, on a
$\phi(k-1)<\phi(k)$. Dans les deux cas, il vient
$\psi(k-1)<\psi(k)$. Enfin l'application $\phi\mapsto
\psi$ est trivialement injective. Elle est aussi surjective,
puisque~$\phi$ est uniquement d\'etermin\'e par les
relations $\phi(k)=\psi(k)-\theta(k)$ $(1\le k\le n)$.\cqfd

\medskip
L'identit\'e de Worpitzky sur les nombres eul\'eriens s'obtient
par simple application de cette proposition. L'ensemble de
toutes les applications de~$[\,n\,]$ dans~$[\,m\,]$ \'etant
not\'e
$H_{n,m}$, soit $\phi\in H_{n,m}$ ; on d\'efinit
$\delta\phi$ comme l'unique $\sigma\in{\goth S}_n$
telle que la suite des paires $(\phi\sigma(1),\sigma(1))$,
$(\phi\sigma(2),\sigma(2))$,
\dots~, $(\phi\sigma(n),\sigma(n))$ soit croissante pour
l'ordre lexicographique. Par cons\'equent,
$\delta^{-1}\sigma$ est l'ensemble des applications
$\phi:[\,n\,]\rightarrow [m]$ telles que $\phi\sigma(i)\le
\phi\sigma(i+1)$ pour $i\in[n-1]$ et telles que l'\'egalit\'e
$\phi\sigma(i)=\phi\sigma(i+1)$ ne soit possible que si
$\sigma(i)<\sigma(i+1)$. D'apr\`es la remarque~2.4, il y a
exactement $s=|\Delta'M\sigma|$ indices~$i$ v\'erifiant
une telle in\'egalit\'e ; donc, d'apr\`es la pr\'ec\'edente
proposition $\card\delta^{-1}\sigma={m+s\choose n}$.
Comme le nombre de permutations $\sigma\in{\goth S}_n$
satisfaisant \`a $|\Delta M\sigma|=s$ est donn\'e par le
nombre eul\'erien $A_{n,s}$, il vient enfin
$$
m^n=\sum_{0\le s\le n-1}{m+s\choose n}A_{n,s}.\eqno(9)
$$
Cette identit\'e est un cas particulier de l'identit\'e (10)
ci-dessous.

\th Propri\'et\'e 2.10|Pour $r\in [\,n\,]$ et $r\le m$, on a :
$$
\sum_{0\le s\le n-r}\rA_{n,n-r-s}{m+s\choose
n}=m^{n-r}\,{m!\over (m-r)!}.\eqno(10)
$$
\finth

\dem 
Notons d'abord que pour $r=1$, on retrouve bien
l'identit\'e (9), puisque l'on a $A_{n,n-1-s}=A_{n,s}$
d'apr\`es la propri\'et\'e~2.3. D'autre part, l'identit\'e est vraie
pour $r=n$ avec les conventions que nous avons adopt\'ees.
On prendra donc $r\in [n-1]$. Soit $H_{m,n,r}$
l'ensemble des applications $\phi:[\,n\,]\rightarrow [m]$
dont la restriction \`a $\{n-r+1,n-r+2,\ldots,n\}$ est {\it
injective}. Il est imm\'ediat que l'on a $\card
H_{m,n,r}=m^{n-r}m!\,/(m-r)!$ 
Pour $\phi\in H_{m,n,r}$ on d\'efinit~$\delta\phi$ comme
\'etant l'unique $\sigma\in{\goth S}_n$ telle que la suite
$(\phi\sigma(1),\sigma(1))$, $(\phi\sigma(2),\sigma(2))$,
\dots~, $(\phi\sigma(n),\sigma(n))$ soit croissante pour
l'ordre lexicographique. Comme pr\'ec\'edemment, on a
$\phi\sigma(i)\le \phi\sigma(i+1)$ pour $i\in[n-1]$,
mais l'\'egalit\'e n'est possible que si l'on a
$\sigma(i)<\sigma(i+1)$ et $\sigma(i)\le n-r$, puisque
la restriction de~$\phi$ \`a l'ensemble
$\{n-r+1,n-r+2,\ldots,n\}$ est injective. Or, d'apr\`es la
remarque~2.4, il y a exactement
$|\Delta'\Delta''^{r-1}M\sigma|$ indices~$i$ tels que
$1\le i\le n-1$, $\sigma(i)<\sigma(i+1)$ et
$\sigma(i)\le n-r$. D'apr\`es la pr\'ec\'edente proposition, on
a $\card \delta^{-1}\sigma={m+s\choose n}$ avec 
$|\Delta'\Delta''^{r-1}M\sigma|=s$ et comme le nombre
de~$\sigma\in {\goth S}_n$ satisfaisant \`a
$|\Delta'\Delta''^{r-1}M\sigma|=s$ est \'egal \`a
$\rA_{n,n-r-s}$, on obtient l'identit\'e d\'esir\'ee.\cqfd

\medskip
Pour terminer cette section, nous donnons la formule
explicite des coefficients $\rA_{n,k}$ obtenue par un
simple calcul traduisant de nouveau l'inversion de
M\"obius (\cf. par exemple [26]) \`a partir des identit\'es~(9)
et~(10).

En effet, pour $1\le n+r$, on a
$\displaystyle{1\over (1-t)^{n+r}}=\sum_{0\le
k}t^k{n+r-1+k\choose n+r-1}$. Donc
$$
\eqalignno{
{\rA_{n-1+r}(t)\over (1-t)^{n+r}}&=
\sum_{0\le k}t^k{n+r-1+k\choose
n+r-1}
\sum_{0\le s\le n-1}\rA_{n+r-1,s}\,t^s\cr
&=\sum_{0\le j}t^j\kern-6pt\sum_{0\le s\le
\min(j,n-1)}\kern-15pt \rA_{n+r-1,s}
{n+r-1+j-s\choose n+r-1}\cr
&=\sum_{0\le j}t^j\sum_{0\le s\le
n-1}\rA_{n+r-1,n-1-s}{j+r+s\choose n+r-1}\cr
&=\sum_{0\le j}t^j(j+r)^{n-1}{(j+r)!\over j!}\cr }
$$
en utilisant le fait que 
${n+r-1+j-s\choose n+r-1}=0$ si $j\le n-2$ et
$s=j+1,\ldots, n-1$ et pour la derni\`ere \'etape en se
servant de l'identit\'e~(10). Il en r\'esulte
$$
\displaylines{
{\rA_{n-1+r}(t)\over  r!\,(1-t)^{n+r}}
=\sum _{0\le j}t^j(j+r)^{n-1}{j+r\choose r}.\cr
\noalign{\hbox{On a donc pour $0\le k\le n-1$ la formule
explicite :}}
\hfill
\rA_{n-1+r,k}=r!\sum_{0\le i\le k}
(-1)^i(k-i+r)^{n-1}{n+r\choose i}{k-i+r\choose r}.
\hfill\llap{(11)}\cr}
$$

\def\ra{{}^{r\kern-2pt}a}

\sectiona 7. Table des polyn\^omes eul\'eriens|
Le paragraphe~4 du pr\'esent chapitre ou l'identit\'e (11)
ci-dessus montre que tous les coefficients des polyn\^omes
$\rA_n(t)$ sont divisibles par~$r!$ (on verra une autre
d\'emonstration de ce r\'esultat \`a la fin du chaitre~IV).
Comme on a pos\'e
$$\displaylines{
\rA_n(t)=\sum_{0\le k\le n-r}\rA_{n,k}\, t^k\quad
(0\le r\le n),\cr
\noalign{\hbox{on peut \'ecrire}}
\rA_{n,k}=r!\,\ra_{n,k},\cr}
$$
o\`u $\ra_{n,k}$ est un entier $(0\le k\le n-r)$.
Le pr\'esent tableau donne les premi\`eres valeurs des
coefficients $\ra_{n,k}$ pour $r=1,2,3,4,5$, $r\le n\le 8$
et $0\le k\le n-r$.

\medskip
$r=1$ :

$$
\vbox{\halign{\vrule\thinspace \hfil$#$\ \vrule
&\strut\quad\hfil$#$\hfil
&\quad\hfil$#$\hfil
&\quad\hfil$#$\hfil
&\quad\hfil$#$\hfil
&\quad\hfil$#$\hfil
&\quad\hfil$#$\hfil
&\quad\hfil$#$\hfil
&\quad \hfil$#$\hfil\quad \vrule\cr
\noalign{\hrule}
k={}&0&1&2&3&4&5&6&7\cr
\noalign{\hrule}
n=1&1&&&&&&&\cr
2&1&1&&&&&&\cr
3&1&4&1&&&&&\cr
4&1&11&11&1&&&&\cr
5&1&26&66&26&1&&&\cr
6&1&57&302&302&57&1&&\cr
7&1&120&1191&2416&1191&120&1&\cr
8&1&247&4293&15619&15619&4293&247&1\cr
\noalign{\hrule}
}}
$$

\smallskip
$r=2$ :

$$
\vbox{\halign{\vrule\thinspace \hfil$#$\ \vrule
&\strut\quad\hfil$#$\hfil
&\quad\hfil$#$\hfil
&\quad\hfil$#$\hfil
&\quad\hfil$#$\hfil
&\quad\hfil$#$\hfil
&\quad\hfil$#$\hfil
&\quad \hfil$#$\hfil\quad \vrule\cr
\noalign{\hrule}
k={}&0&1&2&3&4&5&6\cr
\noalign{\hrule}
n=2&1&&&&&&\cr
3&2&1&&&&&\cr
4&4&7&1&&&&\cr
5&8&33&18&1&&&\cr
6&16&131&171&41&1&&\cr
7&32&473&1208&718&88&1&\cr
8&64&1611&7197&8422&2682&183&1\cr
\noalign{\hrule}
}}
$$

\smallskip
$r=3$ :
$$
\vbox{\halign{\vrule\thinspace \hfil$#$\ \vrule
&\strut\quad\hfil$#$\hfil
&\quad\hfil$#$\hfil
&\quad\hfil$#$\hfil
&\quad\hfil$#$\hfil
&\quad\hfil$#$\hfil
&\quad \hfil$#$\hfil\quad \vrule\cr
\noalign{\hrule}
k={}&0&1&2&3&4&5\cr
\noalign{\hrule}
n=3&1&&&&&\cr
4&3&1&&&&\cr
5&9&10&1&&&\cr
6&27&67&25&1&&\cr
7&81&376&326&56&1&\cr
8&243&1909&3134&1314&119&1\cr
\noalign{\hrule}
}}
$$

\smallskip
$r=4$ :
$$
\vtop{\halign{\vrule\thinspace \hfil$#$\ \vrule
&\strut\quad\hfil$#$\hfil
&\quad\hfil$#$\hfil
&\quad\hfil$#$\hfil
&\quad\hfil$#$\hfil
&\quad \hfil$#$\hfil\quad \vrule\cr
\noalign{\hrule}
k={}&0&1&2&3&4\cr
\noalign{\hrule}
n=4&1&&&&\cr
5&4&1&&&\cr
6&16&13&1&&\cr
7&64&113&32&1&\cr
8&256&821&531&71&1\cr
\noalign{\hrule}
}}$$

\smallskip

$r=5$ :
$$ 
\vtop{\halign{\vrule\thinspace \hfil$#$\ \vrule
&\strut\quad\hfil$#$\hfil
&\quad\hfil$#$\hfil
&\quad\hfil$#$\hfil
&\quad \hfil$#$\hfil\quad \vrule\cr
\noalign{\hrule}
k={}&0&1&2&3\cr
\noalign{\hrule}
n=5&1&&&\cr
6&5&1&&\cr
7&25&16&1&\cr
8&125&171&39&1\cr
\noalign{\hrule}
}}
$$

\vfill\eject

\pagetitretrue

\auteurcourant={{\eightrm CHAPITRE III :
LA FORMULE EXPONENTIELLE}}

\vglue 2cm
\centerline{{\eightrm CHAPITRE III}}
\vskip 2mm
\centerline{{\bf LA FORMULE EXPONENTIELLE}} 

\vskip 6mm plus 2mm
\vskip.4cm
Les trois premi\`eres sections de ce chapitre contiennent la
d\'efinition et quelques propri\'et\'es d'une construction tr\`es
g\'en\'erale que nous appelons \og compos\'e partitionnel\fg.
La motivation de cette notion appara\^\i t dans les sections
suivantes, ainsi qu'au chapitre~IV, o\`u nous appliquons
toutes ces techniques aux polyn\^omes eul\'eriens. Comme
nous l'avons d\'ej\`a mentionn\'e, ces r\'esultats ont \'et\'e
fr\'equemment \'etudi\'es en liaison avec divers probl\`emes
d'\'enum\'eration et tout particuli\`erement dans [29], [14] et
[15].

Dans ce chapitre, si $Z$ est un ensemble non vide, on note
$Z^*$ et $Z^+$ les mono\"\i des libre et ab\'elien libre
engendr\'es par~$Z$. On appellera {\it mots} les \'el\'ements
de~$Z^*$, qu'on pr\'esentera comme des suites
$g=z_1z_2\ldots z_r$, o\`u $z_1$, $z_2$, \dots~, $z_r$
appartiennent \`a~$Z$ ; l'entier~$r$ est la {\it longueur}
du mot~$g$. On appellera {\it mon\^omes} les \'el\'ements
de~$Z^+$. Si~$\alpha$ est le morphisme canonique
de~$Z^*$ sur~$Z^+$, on prendra dans chaque classe
$\alpha^{-1}(f)$, o\`u $f\in Z^+$, un mot canonique qu'on
identifiera \`a~$f$. Le mon\^ome~$f$ sera dit de {\it
degr\'e}~$r$ si le {\it mot}~$f$ est de longueur~$r$.

\sectiona 1. La formule de Hurwitz|
Dans cette section, $Y$ est un ensemble non vide, muni
d'une application $\lambda:Y\rightarrow {\bboard N}$. Le
m\^eme symbole d\'esignera les morphismes dans~${\bboard
N}$ \'etendant cette application aux mono\"\i des $Y^*$
et~$Y^+$. Le morphisme canonique de~$Y^*$
sur~$Y^+$ sera not\'e~$\alpha$.

Soit $\cal P$ l'ensemble des parties {\it finies} (y
compris la partie vide) de~$\bboard N$ ; on consid\`ere le
sous-ensemble~$Y_\lambda$ du produit cart\'esien
$Y\times {\cal P}$ compos\'e de tous les couples $(y,I)$
satisfaisant \`a la condition $\card I=\lambda y$.
On forme ensuite le mono\"\i de libre $Y_\lambda^*$
engendr\'e par~$Y_\lambda$. Soit
$$
\eqalignno{h&=(y_1,I_1)(y_2,I_2)\ldots (y_r,I_r)\cr
\noalign{\hbox{un mot de $Y_\lambda^*$ de longueur
$r\ge 1$. On pose}}
\beta h=g&=y_1y_2\dots y_r\in Y^*\quad {\rm et}\quad
\lambda h=\lambda g.&(1)\cr}
$$ 

\rem D\'efinition $3.1$|Pour chaque entier $r\ge 1$, le
{\it compos\'e partitionnel marqu\'e} de~$Y$ {\it de
degr\'e}~$r$ est le sous-ensemble $Y^{((r))}$
de~$Y_\lambda^*$ form\'e de tous les mots
$h=(y_1,I_1)(y_2,I_2)\ldots (y_r,I_r)$ de longueur~$r$
satisfaisant aux conditions
\decale (i)|$I_j\cap I_{j'}=\emptyset$ si $j\not=j'$ ;
\decale (ii)|$\bigcup \{I_j:j\in [r]\}=[\lambda h]$.

\medskip
On notera que si $\lambda y_j$ est strictement positif
pour tout $j\in [r]$, la {\it famille} $\{I_1,I_2,\ldots,
I_r\}$ est une {\it partition} de l'ensemble $[\lambda
h]$, aucun de ces sous-ensembles n'\'etant vide.

Par convention, on supposera l'existence d'un ensemble
$Y^{((r))}$ pour $r=0$ contenant un seul \'el\'ement, \`a
savoir l'\'el\'ement neutre de~$Y$. Ce dernier est envoy\'e
par~$\beta$ sur l'\'el\'ement neutre de~$Y^*$. On
identifiera, d'autre part, le compos\'e partitionnel marqu\'e
$Y^{((1))}$ de degr\'e~1 avec~$Y$.

\th Th\'eor\`eme 3.3 {\rm (Formule de Hurwitz)}|Soit
$E(Y)=\sum\{y/\lambda y! : y\in Y\}$ la fonction
g\'en\'eratrice exponentielle de~$Y$ $($par rapport
\`a~$\lambda)$. Pour tout $r\ge 0$, on a dans la $\bboard
Q$-alg\`ebre large de~$Y^*$ l'identit\'e
$$
\bigl(E(Y)\bigr)^r=\sum\{\beta h/\lambda h! :h\in
Y^{((r))}\}.\eqno(2)
$$
\finth

\dem
Avec nos conventions sur $\beta$, il n'y a rien \`a prouver
pour $r=0$. Pour chaque mot $g=y_1y_2\ldots y_r$ de
longueur~$r$ $(r\ge 1)$ de~$Y^*$, le nombre de mots
$h\in Y^{((r))}$ tels que $\beta h=g$ est \'egal au nombre
de suites $(I_1,I_2,\ldots, I_r)$ de parties de~$\bboard
N$ satisfaisant aux conditions~(i) et~(ii) de la
d\'efinition~3.1 ainsi qu'\`a la condition $\card I_j=\lambda
y_j$ pour chaque $j\in [r]$. Or, le nombre de telles suites
est \'evidemment donn\'e par le coefficient multinomial
$$\displaylines{
{\lambda\choose g}={\lambda g!\over \lambda
y_1!\,\lambda y_2!\,\ldots\,\lambda y_r!}.\cr
\noalign{\hbox{Donc}}
\sum\{\beta h/\lambda h! :h\in
Y^{((r))},\,\beta h=g\}
={1\over  \lambda
y_1!\,\lambda y_2!\,\ldots\,\lambda y_r!}\,g,\cr}
$$
puisque $\lambda h=\lambda g$ pour $\beta h=g$. Or, le
facteur $1/(\lambda
y_1!\,\lambda y_2!\,\ldots\,\lambda y_r!)$ est simplement
le coefficient de~$g$ dans le d\'eveloppement de 
$\bigl(E(Y)\bigr)^r$ et le r\'esultat s'en d\'eduit par
sommation sur tous les mots de longueur~$r$
de~$Y^*$.\cqfd

\sectiona 2. Le compos\'e partitionnel|
Nous conservons les m\^emes notations que dans la
section~1, mais nous supposons cette fois que $\lambda y
>0$ pour tout $y\in Y$. Soit $h=(y_1,I_1)(y_2,I_2)\ldots
(y_r,I_r)$ un mot de $Y^{((r))}$. L'hypoth\`ese ci-dessus
entra\^\i ne, puisque l'on a $\card I_j=\lambda y_j$ pour
chaque $j\in [r]$, que les sous-ensembles $I_1$, $I_2$,
\dots~, $I_r$ sont non vides, donc tous distincts. Il en
r\'esulte que~$h$ est {\it multilin\'eaire}, c'est-\`a-dire a
toutes ses lettres distinctes. En notant~$\delta$ le
morphisme canonique de~$Y_\lambda^*$
sur~$Y_\lambda^+$, on voit donc que la classe ab\'elienne
$\delta^{-1}\delta h$ contient exactement~$r!$ mots. De
plus, les conditions~(i) et~(ii) de la d\'efinition~3.1 ne
faisant pas intervenir l'ordre des lettres de~$h$, il
s'ensuit que $Y^{(r))}$ contient toute une classe
ab\'elienne 
$\delta^{-1}\delta h$ d\`es qu'il contient~$h$.

\rem D\'efinition $3.3$|On appelle {\it compos\'e
partitionnel} de $Y$, {\it de degr\'e}~$r$ $(r\ge 0)$,
l'ensemble $Y^{(r)}=\delta Y^{((r))}$ et l'union
$Y^{(+)}=\smash{\bigcup\limits_{0\le r}}Y^{(r)}$ est le
{\it compos\'e partitionnel} de~$Y$.

\medskip
Les \'el\'ements de $Y^{(r)}$ sont donc des {\it mon\^omes}
$f=(y_1,I_1)(y_2,I_2)\ldots (y_r,I_r)$
de~$Y_\lambda^+$. On d\'esigne par $\gamma f$ le {\it
mon\^ome} $m=y_1y_2\ldots y_s$ de~$Y^+$ et l'on pose
encore $\lambda f=\lambda m$. On a donc l'identit\'e
$$
\alpha \beta=\gamma\delta
$$
o\`u $\alpha$ est le morphisme $\alpha:Y^*\rightarrow
Y^+$ et o\`u $\beta$ a \'et\'e d\'efini en~(1). On peut
rassembler les remarques ainsi faites dans un lemme.

\th Lemme 3.4|Pour tout $f\in Y^{(r)}$, on a
$\card\{h\in Y^{((r))}:\delta h=f\}=r!$
et si $\delta h=f$, on a : $\gamma f=\alpha \beta h$.
\finth

Venons-en \`a la formule fondamentale de ce chapitre.

\th Th\'eor\`eme 3.5 {\rm (Formule exponentielle)}|Dans la
$\bboard Q$-alg\`ebre large de~$Y^+$, on a l'identit\'e
$$
\sum\{\gamma f/\lambda f!:f\in Y^{(+)}\}=\exp E(Y).
$$
\finth

\dem On a $\sum\{\gamma f/\lambda f!:f\in Y^{(0)}\}=1$.
D'autre part, pour $f\in Y^{(r)}$ $(r\ge 1)$, on a
$$
\displaylines{\noalign{\vskip -8pt}
\sum\{\alpha \beta h/\lambda h!: \delta h=f\}=r!\,
(\gamma f/\lambda f!)\cr
\noalign{\hbox{d'apr\`es le lemme pr\'ec\'edent. D'o\`u, il
r\'esulte}}
\sum\{\gamma f/\lambda f!:f\in Y^{(r)}\}
={1\over r!}\sum\{\alpha \beta h/\lambda h!: h\in
Y^{((r))}\}\cr}
$$
pour chaque $r\ge 1$. Utilisant le th\'eor\`eme 3.2, on
obtient donc par sommation sur tous les $r\in {\bboard
N}$
$$
\sum\{\gamma f/\lambda f! :f\in Y^{(+)}\}=\sum_{0\le r}
{1\over r!}\bigl(E(Y)\bigr)^r=\exp E(Y).\qed
$$

\rem Remarque $3.6$|Dans la $\bboard Q$-alg\`ebre large
de $Y^+$, on a pris la topologie des s\'eries formelles
induite par l'ordre~o suivant : si $a=\sum\{m\,a_m:m\in
Y^+\}$ est une s\'erie formelle, son ordre ${\rm o}(a)$ est
d\'efini par $${\rm o}(a)=\inf\{n\ge 1:\lambda m=n,
a_m\not=0\}.$$
Utilisons les notations abr\'eg\'ees
$$\displaylines{
\gamma\{Y^{(+)}\cap \lambda^{-1}n\}=\sum\{\lambda
f:f\in Y^{(+)},\,\lambda f=n\}\quad (n\ge 0);\cr
\{Y_n\}=\sum\{y:y\in Y,\,\lambda y=n\} \quad (n\ge
1).\cr}
$$
Les s\'eries formelles $\gamma\{Y^{(+)}\cap
\lambda^{-1}n\}$ et $\{Y_n\}$ sont d'ordre \'egal \`a~$n$,
ce qui permet d'\'ecrire la formule exponentielle sous la
forme :
$$
\sum_{0\le n}{1\over n!}\,\gamma\{Y^{(+)}\cap
\lambda^{-1}n\}=\exp\Bigl(\sum_{1\le n}{1\over
n!}\,\{Y_n\}\Bigr).\eqno(3)
$$

\sectiona 3. Une formule d'inversion pour les s\'eries
exponentielles|
Pour $f\in Y^{(+)}$ d\'esignons par $z(f)$ l'unique $r\in
{\bboard N}$ tel que $f\in Y^{(r)}$ ; l'entier~$z(f)$ n'est
autre que le {\it degr\'e} de~$f$. Posons
$$
\overline \gamma\{Y^{(+)}\cap \lambda^{-1}n\}
=\sum\{\lambda f\cdot (-1)^{z(f)+n}
:f\in Y^{(+)},\,\lambda f=n\}\quad (n\ge 0).
$$

\th Propri\'et\'e 3.7|Dans la $\bboard Q$-alg\`ebre large de
$Y^+$, on a l'identit\'e
$$
\Bigl(\sum_{0\le n}{1\over n!}\,
\gamma\{Y^{(+)}\cap \lambda^{-1}n\}\Bigr)^{-1}
=\sum_{0\le n}{(-1)^n\over n!}\,
\overline \gamma\{Y^{(+)}\cap \lambda^{-1}n\}.
\eqno(4)
$$
\finth

\dem
D'apr\`es le th\'eor\`eme 3.5, le membre de gauche de l'identit\'e
\`a \'etablir, soit~$U$, est \'egal \`a $(\exp(E(Y))^{-1}$,
c'est-\`a-dire \`a $\exp(-E(Y))$. Notons~$\phi$ le
morphisme envoyant sur $-y$ chaque $y\in Y$ ; ceci
\'equivaut \`a $U=\phi\exp E(Y)$, donc de nouveau d'apr\`es le
th\'eor\`eme~3.5, \`a
$$\eqalignno{
U=\phi \sum_{0\le n}{1\over n!}\,
\gamma\{Y^{(+)}\cap \lambda^{-1}n\}
&=\sum_{0\le n}{1\over n!}\sum_{0\le r}(-1)^r
\gamma\{Y^{(+)}\cap \lambda^{-1}n\cap
z^{-1}r\}\cr
&=\sum_{0\le n}{(-1)^n\over n!}\,
\overline \gamma\{Y^{(+)}\cap \lambda^{-1}n\}.\qed\cr}
$$

\medskip
Ceci termine l'\'etablissement des formules que nous
utiliserons par la suite. Le th\'eor\`eme~3.2 avc une
interpr\'etation ad\'equate des objets en cause exprime que la
{\it transformation de Borel}
$\sum\{y:y\in Y\}\mapsto \sum\{y/\lambda y!:y\in Y\}$
est un morphisme dans l'alg\`ebre large de~$Y^+$ de
l'alg\`ebre large de base~$Y$ par rapport au \og produit
d'intercalement\fg\ (\og shuffle\fg\ de Chen, Fox et
Lyndon). Le th\'eor\`eme~3.5 est appel\'e \og formule de
Cauchy\fg\ dans les probl\`emes concernant le groupe
sym\'etrique. Sous une forme ou sous une autre, elle a \'et\'e
retrouv\'ee et utilis\'ee souvent dans diverses questions
d'\'enum\'eration. Nous l'appellerons simplement {\it formule
exponentielle}. En prenant le logarithme, on obtiendrait
\'evidemment la fonction g\'en\'eratrice exponentielle~$E(Y)$
de~$Y$ en fonction de la fonction g\'en\'eratrice
exponentielle du compos\'e partitionnel $Y^{(+)}$ de~$Y$.

\rem D\'efinition $3.8$|Soient $Y^{(+)}$ un compos\'e
partitionnel et $A$ un mono\"\i de ab\'elien ; une
application
$\mu:Y^{(+)}\rightarrow A$ sera dite {\it multiplicative}
si et seulement s'il existe un morphisme
$\mu':=Y^+\rightarrow A$ tel que le diagramme suivant
soit commutatif.

\def\fleche(#1,#2)\dir(#3,#4)\long#5{%
{\leftput(#1,#2){\vector(#3,#4){#5}}}}

$$\vbox{\offinterlineskip
\centerput(0,0){$Y^{(+)}$}
\fleche(5,0)\dir(1,0)\long{20}
\centerput(28,0){$A$}
\centerput(15,2){$\mu$}
\centerput(-2,-6){$\gamma$}
\fleche(0,-2)\dir(0,-1)\long{7}
\fleche(5,-11)\dir(2,1)\long{18}
\centerput(15,-10){$\mu'$}
\centerput(0,-14){$Y^+$}
}\hskip3cm
$$

\vfill\eject
\goodbreak
\sectiona4. Le compos\'e partitionnel des applications|
Il existe de nombreuses familles de structures qui peuvent
\^etre consi\-d\'er\'ees comme le compos\'e partitionnel d'une de
leurs sous-familles. Nous examinerons ici, \`a titre
d'exemple, la famille des applications avec le but
d'introduire les notions n\'ecessaires pour traiter le cas
particulier des permutations.

\rem D\'efinition $3.9$|Soit $f:I\rightarrow I$ une
application d'un ensemble fini~$I$ dans lui-m\^eme. L'{\it
\'equivalence}~$f^*$ de~$f$ est la relation d'\'equivalence
dans $I\times I$ telle que deux \'el\'ements $i$ et $i'$
de~$I$ appartiennent \`a la m\^eme classe si et seulement
s'il existe des it\'er\'ees~$f^p$ et $f^{p'}$de~$f$
satisfaisant \`a
$f^p(i)=f^{p'}(i')$.

\medskip
Nous appellerons {\it sous-domaines} de $f$ les classes
de cette \'equivalence et leur nombre sera d\'esign\'e
par~$z(f)$. L'application~$f$ sera {\it connexe} si
$z(f)=1$. Ainsi les sous-domaines d'une permutation sont
les orbites de celle-ci ; les permutations circulaires sont
les permutations connexes.

Dans la suite, on notera $F_n$  l'ensemble des
applications de $[\,n\,]$ dans lui-m\^eme $(n\ge 0)$ et
l'on posera
$F=\bigcup_{0\le n} F_n$. D\'esignons par $I_1$, $I_2$,
\dots~, $I_r$ $(r=z(f))$ les sous-domaines d'une
application $f\in F_n$ $(n\ge 1)$. Pour tout $j\in [r]$,
on note~$\omega_j$ l'unique morphisme (d'ensembles
ordonn\'es) $\omega_j:[\card I_j]\rightarrow [\,n\,]$ qui a
pour image~$I_j$ et $f_j'$ la restriction de~$f$
\`a~$I_j$. Par d\'efinition de l'\'equivalence~$f^*$ on voit
que $f'_j(I_j)\subset I_j$ et il est licite de poser
$f_j=\omega_j^{-1}f'_j\,\omega_j$ $(j\in [r])$. Les
applications $f_j$ envoient $[\card I_j]$ dans lui-m\^eme
et sont toutes connexes $(j\in [r])$. Enfin, il est clair
que toute application $f\in F_n$ d\'etermine, de fa\c con
biunivoque le {\it mon\^ome} (appartenant au mono\"\i de
$(F\times{\cal P})^+$, o\`u $\cal P$ d\'esigne toujours
l'ensemble des parties finies de~$\bboard N$)
$(f_1,I_1)(f_2,I_2)\ldots (f_r,I_r)$, que l'on appellera sa
{\it factorisation canonique}, les~$f_j$ eux-m\^emes \'etant
les {\it facteurs} de~$f$. Par commodit\'e, on identifiera
tout $f\in F$ avec sa factorisation canonique et $f\in
F_0$ avec le mon\^ome unit\'e.

Le raccordement avec les trois premi\`eres sections se fait
de la fa\c con sui\-vante. Soit donn\'ee une famille $\cal F$
d'applications connexes dont les domaines sont des
ensembles de la forme $[\,n\,]$ $(n\in {\bboard N})$.
Posons $Y={\cal F}$ et prenons pour~$\lambda$
l'application qui envoie sur~$n$ chaque $f\in {\cal F}$
de domaine $[\,n\,]$ $(n\in {\bboard N})$. Formons
ensuite le compos\'e partitionnel ${\cal F}^{(+)}$. On
constate alors que la factorisation canonique d'une
application $f\in F$ appartient au compos\'e partitionnel 
${\cal F}^{(+)}$ si et seulement si les facteurs de~$f$
appartiennent \`a~${\cal F}$. Avec l'identification faite
ci-dessus, on a ainsi la proposition suivante.

\th Proposition 3.10|Soit ${\cal F}\subset F$ une famille
d'applications connexes. Le compos\'e partitionnel ${\cal
F}^{(+)}$ est l'ensemble des applications $f\in F$ dont
les facteurs appartiennent \`a~${\cal F}$.
\finth

La propri\'et\'e suivante d\'ecoule imm\'ediatement de la
d\'efinition du compos\'e partitionnel d'un ensemble
d'applications. Elle exprime le fait que la factorisation
canonique d'une application~$f$ conserve les exc\'edances
et les points fixes de~$f$. 

\th Propri\'et\'e 3.11|Soit $(f_1,I_1)(f_2,I_2)\ldots (f_r,I_r)$
$(r\ge 1)$ la factorisation canonique d'une
application~$f$. Pour tout $j\in [r]$, le morphisme
$\tau_j=\omega_j^{-1}$ est une bijection de~$I_j$ sur
$[\card I_j]$ telle que pour tout $i\in I_j$ on ait les
\'equivalences :\qquad
$i<f(i)\Leftrightarrow \tau_j(i)<f_j\,
\tau_j(i)$ ; \quad$i=f(i)\Leftrightarrow \tau_j(i)=f_j\,
\tau_j(i)$ ;
$i>f(i)\Leftrightarrow \tau_j(i)>f_j\,
\tau_j(i)$.
\finth

\dem
En effet, si l'entier $i$ est dans le sous-domaine $I_j$,
on~a
$f_j(i)=\omega_j^{-1}f'_j\,\omega_j(i)
=\tau_jf'_j\,\tau_j^{-1}(i)$
o\`u $f'_j$ est la restriction de~$f$ \`a~$I_j$. Les
\'equivalences ci-dessus r\'esultent alors du fait que
$\tau_j:I_j\rightarrow [\card I_j]$ est un morphisme
strictement croissant.\cqfd

\medskip
R\'ecrivons la formule exponentielle (3) et la formule
d'inversion~(4) dans ce cas particulier du compos\'e
partitionnel des applications. On a d'abord
${\cal F}^{(+)}\cap \lambda^{-1}n=F_n\cap {\cal
F}^{(+)}$ pour $n\ge 0$ et ${\cal F}\cap
\lambda^{-1}n=F_n\cap {\cal F}$ pour $n\ge 1$ et les
deux identit\'es (3) et (4) se pr\'esentent ainsi
$$
\eqalignno{
\sum_{0\le n}{1\over n!}\,\gamma\{F_n\cap
{\cal F}^{(+)}\}&=\exp\Bigl(\sum_{1\le n}{1\over
n!}\,\{F_n\cap {\cal F}\}\Bigr);&(5)\cr
\Bigl(\sum_{0\le n}{1\over n!}\,
\gamma\{F_n\cap
{\cal F}^{(+)}\}\Bigr)^{-1}
&=\sum_{0\le n}{(-1)^n\over n!}\,
\overline \gamma\{F_n\cap
{\cal F}^{(+)}\}.
&(6)\cr
}
$$

En fait, ces deux identit\'es seront appliqu\'ees ci-apr\`es
sous la forme sui\-vante. On suppose donn\'ee une
application multiplicative $\mu:{\cal F}^{(+)}\rightarrow
\Omega$ (\cf. d\'efinition~3.8). On forme ensuite l'alg\`ebre
sur~$\bboard Q$ du mono\"\i de~$\Omega$ ;
notons~$\overline
\Omega$ cette alg\`ebre. On consid\`ere enfin l'alg\`ebre
$\overline\Omega[[u]]$ des s\'eries formelles \`a
coef\-ficients dans~$\overline \Omega$ et \`a une
ind\'etermin\'ee~$u$.

\th Proposition 3.12|Soit $\mu : {\cal
F}^{(+)}\rightarrow \Omega$ une application
multiplicative. Dans l'alg\`ebre des s\'eries formelles
$\overline\Omega[[u]]$, on a les identit\'es :
$$
\eqalignno{
\sum_{0\le n}{u^n\over n!}\,\mu\{F_n\cap
{\cal F}^{(+)}\}&=\exp\Bigl(\sum_{1\le n}{u^n\over
n!}\,\mu\{F_n\cap {\cal F}\}\Bigr);&(7)\cr
\Bigl(\sum_{0\le n}{u^n\over n!}\,
\mu\{F_n\cap
{\cal F}^{(+)}\}\Bigr)^{-1}
&=\sum_{0\le n}{(-u)^n\over n!}\,
\overline \mu\{F_n\cap
{\cal F}^{(+)}\}.
&(8)\cr
}
$$
o\`u $\overline \mu f=(-1)^{z(f)+n}\mu f$ pour tout $f\in
F_n\cap {\cal F}^{(+)}$ $(n\ge 0)$.
\finth

\dem
Soit $\mu'$ le morphisme de ${\cal F}^+$ dans
$\Omega$ tel que $\mu=\mu'\gamma$. D\'esignons
par~$\phi$ l'application envoyant tout $f\in F_n\cap
{\cal F}$ sur le mon\^ome $u^n\mu'f$ $(n\ge 0)$.
Comme~$\mu$ est multiplicative, on a $\phi\gamma
f=u^n\mu f$ pour tout $f\in F_n\cap
{\cal F}^{(+)}$ $(n\ge 0)$. Tous les mon\^omes
$\phi\gamma f$ o\`u $f\in F_n\cap
{\cal F}^{(+)}$ sont donc de degr\'e~$n$ (en~$u$). On
peut donc prolonger~$\phi$ en un morphisme continu de
la $\bboard Q$-alg\`ebre large de~${\cal F}^+$ dans
$\overline\Omega[[u]]$. Appliquant ainsi~$\phi$ aux
deux membres des deux identit\'es~(5) et~(6), on obtient
les identit\'es~(7) et~(8).\cqfd

\sectiona 5. Applications|
Il est \'evident que si $\cal F$ est l'ensemble~$\cal C$
des {\it permutations circulaires}, le compos\'e
partitionnel ${\cal C}^{(+)}$ est exactement l'ensemble
${\goth S}=\bigcup_{0\le n}{\goth S}_n$. On a de plus
$F_n\cap
{\cal F}^{(+)}={\goth S}_n$ pour $n\ge 0$ et $F_n\cap
{\cal F}={\cal C}_n$ pour $n\ge 1$.
Enfin, si $\sigma$ est dans ${\goth S}_n$, se rappelant
que $z(\sigma)$ est le nombre des orbites de~$\sigma$,
on voit que le coefficient $(-1)^{z(\sigma)+n}$ est la
{\it signature} $\epsilon(\sigma)$ de~$\sigma$. Dans ces
conditions les deux identit\'es~(7) et~(8) s'\'ecrivent :
$$
\eqalignno{
\sum_{0\le n}{u^n\over
n!}\,\mu\{{\goth S}_n\}&=\exp\Bigl(\sum_{1\le
n}{u^n\over n!}\,\mu\{{\cal C}_n\}\Bigr);&(9)\cr
\Bigl(\sum_{0\le n}{u^n\over n!}\,
\mu\{{\goth S}_n\}\Bigr)^{-1}
&=\sum_{0\le n}{(-u)^n\over n!}\,
\overline \mu\{{\goth S}_n\}.
&(10)\cr
}
$$
o\`u $\overline \mu \sigma=\epsilon(\sigma)\mu\sigma$
pour tout $\sigma\in {\goth S}$. Donnons quelques
exemples d'application des formules~(9) et~(10).

\medskip
Soit $(x_n)_{n\ge 1}$ une suite d'ind\'etermin\'ees
commutatives. Si $\sigma$ est dans~${\goth S}_n$,
posons
$\mu\sigma=x_1^{m_1}x_2^{m_2}\ldots x_n^{m_n}$, o\`u,
pour tout $k\in [\,n\,]$, l'entier~$m_k$ est le nombre de
cycles de longueur~$k$ dans la $\sigma$. La
fonction~$\mu$ ainsi d\'efinie est multiplicative. Dans ce
cas, $\mu\{{\goth S}_n\}$ est le {\it polyn\^ome indicateur
de cycles} de~${\goth S}_n$ ou polyn\^ome de Bell (\cf.
[24] p.~68) et $(u^n/n!)\mu\{{\cal C}_n\}$ se r\'eduit \`a
$u^nx_n/n$ puisque $\card {\cal C}_n=(n-1)!$ La
formule exponentielle permet donc de retrouver
l'expression de la fonction g\'en\'eratrice exponentielle de
ces polyn\^omes, \`a savoir
$$\sum_{0\le n}{u^n\over
n!}\,\mu\{{\goth S}_n\}=\exp\Bigl(\sum_{1\le
n}{u^n\over n}\,x_n\Bigr).
$$

Un autre exemple de compos\'e partitionnel est donn\'e par
l'ensemble~$U$ des applications {\it ultimement
idempotentes}, c'est-\`a-dire, pour tout $n\ge 0$, des
applications $f\in F_n$ telles que $f^n =f^{n-1}$.
Consid\'erons, en effet, pour tout entier $n\ge 1$,
l'ensemble~$V_n$ des applications $f\in F_n$, dont
l'image de la $(n-1)$\ieme~it\'er\'ee~$f^{n-1}$ soit r\'eduite
\`a un seul point. Les \'el\'ements de~$V_n$ sont encore
appel\'es {\it arborescences}. Posons $V=\bigcup_{1\le
n}V_n$ ; il est alors clair que le compos\'e partitionnel
de~$V$ est l'ensemble~$U$. Posons $U_n=F_n\cap U$ pour
$n\ge 0$ ; on obtient donc pour toute application
multiplicative $\mu:U\rightarrow \Omega$, deux identit\'es
analogues \`a~(9) et~(10) en substituant~$U_n$ \`a~${\goth
S}_n$, $V_n$ \`a~${\cal C}_n$ et en posant $\overline
\mu f=(-1)^{z(f)+n}$ pour tout $f\in U_n$ $(n\ge 1)$.
On pose naturellement $\mu U_0=\overline\mu U_0=1$.

En particulier, prenons pour $\mu$ l'application qui
envoie sur~1 tout $f\in U$. Comme on a, de fa\c con
\'evidente $\card V_n=n\,\card U_{n-1}$ pour $n\ge 1$, il
vient
$$
\sum_{0\le n}{u^n\over n!}\,\mu\{V_n\}=u\sum_{0\le
n}{u^n\over n!}\card U_n
$$
et on retrouve, en appliquant la formule (7), ce r\'esultat
bien connu : que la s\'erie formelle
$w=\smash{\sum\limits_{0\le n}}(u^n/n!)\card U_n$ est
solution dans 
${\bboard Q}[[u]]$ de l'\'equation $w=\exp(uw)$.

\goodbreak
Enfin, prenons pour $\cal F$ l'ensemble $G$ de toutes les
applications connexes de~$F$. Dans ce cas, le compos\'e
partitionnel de~$G$ est~$F$ tout entier et l'on obtient
encore, pour toute application multiplicative 
$\mu:F\rightarrow \Omega$ donn\'ee, deux identit\'es
analogues \`a~(9) et \`a ~(10).

\def\per{\mathop{\rm per}\nolimits}

\sectiona 6. Une identit\'e entre d\'eterminants et
permanents|
Dans l'\'enonc\'e qui suit, $\Xi$ est une matrice infinie 
$\Xi=(\xi_{i,j})_{(i,j=1,2,\ldots\,)}$ \`a coefficients dans
un anneau commutatif~$\overline \Omega$ ; on d\'esigne
pour tout $n\ge 1$ par~$\Xi_n$ la matrice 
$(\xi_{i,j})_{(1\le i,j\le n)}$, par $\det \Xi_n$ son
d\'eterminant et par $\per \Xi_n$ son permanent.

\th Th\'eor\`eme 3.13|Soient $a$, $b$ et $c$ trois \'el\'ements
de $\overline \Omega$ et $\Xi$ une matrice infinie
ayant ses coefficients supradiagonaux $($resp. diagonaux,
resp. infradiagonaux$)$ \'egaux \`a~$a$
$($resp.~$b$, resp.~$c)$. On a l'identit\'e
$$
\Bigl(1+\sum_{1\le n}{u^n\over n!}\per \Xi_n\Bigr)^{-1}
=1+\sum_{1\le n}{(-u)^n\over n!}\det \Xi_n.\eqno(11)
$$
\finth

\dem
D\'esignons par $\xi_{i,j}$ les \'el\'ements de la
matrice~$\Xi$, puis posons
$\mu\sigma=\xi_{1,\sigma(1)}\xi_{2,\sigma(2)}\ldots
\xi_{n,\sigma(n)}$ pour tout $\sigma\in {\goth S}_n$.
Avec les notations de la propri\'et\'e 3.11 (en prenant
$f=\sigma$), on voit que pour un entier~$i$ appartenant
au sous-domaine (i.e. \`a l'orbite)~$I_j$, on a
$\xi_{i,\sigma(i)}=\xi_{\tau_j(i),f_j\,\tau_j(i)}$.
On a ainsi $\mu\sigma=\prod_j\prod_i\xi_{i,\sigma(i)}$,
o\`u $j$ varie dans $[r]$ et o\`u~$i$, pour~$j$ fix\'e, varie
dans~$I_j$. Si~$i$ est dans~$I_j$, l'\'el\'ement
$i'=\tau_j(i)$ est dans $[\card I_j]$ et l'on a, d'apr\`es ce
qui pr\'ec\`ede,
$\mu\sigma=\prod_j\prod_{i'}\xi_{i',f_j(i')}
=\prod_j\mu f_j$. L'application~$\mu$ est donc
multiplicative et l'on peut \'ecrire
$$\eqalignno{
\mu\{{\goth S}_n\}&=\sum_\sigma
\xi_{1,\sigma(1)}\cdots \xi_{n,\sigma(n)}=\per \Xi_n;\cr
\overline\mu\{{\goth S}_n\}&=\sum_\sigma
\epsilon(\sigma)
\xi_{1,\sigma(1)}\cdots \xi_{n,\sigma(n)}=\det \Xi_n.\cr
}
$$
Le th\'eor\`eme 3.13 r\'esulte alors de l'identit\'e (10).\cqfd

\medskip
\rem Remarque $3.14$|Notons que l'identit\'e (11) n'est pas
vraie pour toutes les matrices, mais comme l'a remarqu\'e
Kittel [18], on peut construire d'autres matrices
infinies~$\Xi$ que celles consid\'er\'ees dans l'\'enonc\'e du
th\'eor\`eme~3.13 pour lesquelles l'identit\'e (11) est v\'erifi\'ee.
Par exemple, si la premi\`ere colonne de la matrice~$\Xi$
n'a que des z\'eros, on a $\per \Xi_n=\det \Xi_n=0$ pour
tout $n\ge 1$ et l'identit\'e (11) est trivialement v\'erifi\'ee.

De m\^eme, consid\'erons la matrice~$\Xi$ d\'efinie par
$$
\xi_{i,j}=\cases{1,&si $1\le i\le j$ ;\cr
-(i-1),&si $1\le i-1=j$ ;\cr
0,&si $1\le j\le i-2$.\cr}
$$
On obtient facilement $\per\Xi_1=1$ et $\per\Xi_n=0$
pour tout $n\ge 2$, ainsi que $\det\Xi_n=n!$ pour tout
$n\ge 1$. L'identit\'e~(11) est encore v\'erifi\'ee ; on retrouve
en fait l'identit\'e
$(1+u)^{-1}=\sum_{0\le n}(-u)^n$.

\goodbreak
Notons enfin le r\'esultat \'el\'ementaire (pour $n\ge 1$)
$$
\det\Xi_n=\cases{\displaystyle{c(b-a)^n-a(b-c)^n\over
c-a},&si $c\not=a$ ;\cr
\noalign{\smallskip}
=(b-a)^{n-1}(b+(n-1)a),&si $c=a$.\cr}
$$
Portant ces valeurs dans la formule (11), on est conduit \`a
l'identit\'e
$$\displaylines{\quad
\Bigl(1+\sum_{1\le n}{u^n\over n!}\per\Xi_n\Bigr)^{-1}
\hfill\cr
\noalign{\vskip-6pt}
\hfill{}
=\cases{\displaystyle
{c\exp((a-b)u)-a\exp((c-b)u)\over c-a},\quad\hbox{si
$c\not=a$} ;\cr
\noalign{\medskip}
(1-au)\exp((a-b)u),\quad\hbox{si $a\not=b$ et
$c=a$}.\cr}\qquad\qquad(12)\cr}
$$

\vfill\eject

\pagetitretrue

\auteurcourant={\eightrm CHAPITRE IV :
FONCTIONS G\'EN\'ERATRICES}

\titrecourant={\eightrm 1. FONCTION G\'EN\'ERATRICE
EXPONENTIELLE}
\vglue 2cm
\centerline{{\eightrm CHAPITRE IV}}
\vskip 2mm
\centerline{\bf FONCTIONS G\'EN\'ERATRICES}
\smallskip
\centerline{{\bf DES POLYN\^OMES EUL\'ERIENS}} 

\def\zA{{}^{0\kern-3pt}A}

\vskip 6mm plus 2mm
\sectiona 1. Fonction g\'en\'eratrice exponentielle
de ${}^{0\kern-3pt}A_n(t)$, $A_n(t)$, $B_n(t)$|
Pour $\sigma\in {\goth S}_n$ $(n\ge 1)$, nous
d\'efinissons $E'\sigma\in \{0,1\}^n$ par la condition
$E'\sigma(k)=1$ ou~0 selon  que~$k$ est ou non un point
fixe de~$\sigma$. De par la d\'efinition de~$E$ et $\Delta
E$ on a donc imm\'ediatement :
$|E\sigma|=|E'\sigma|+|\Delta E\sigma|$.
Introduisant une nouvelle ind\'etermin\'ee~$t'$, nous posons
$\theta'\sigma=t'{}^{|E'\sigma|}t^{|\Delta E\sigma|}$
et $\overline A_n(t,t')=\theta'{\goth S}_n$ ($=1$ pour
$n=0$). Par cons\'equent, on a
$$
\displaylines{
\zA_n(t)=\overline A_n(t,t) ;\quad
A_n(t)=\overline A_n(t,1);\cr
\overline A_n(t,0)=\theta\Delta E\{\sigma\in {\goth S}_n
:|E'\sigma|=0\}
=\theta \Delta E{\cal D}_n=\theta E{\cal D}_n=B_n(t).\cr}
$$
(Voir la fin du paragraphe 1 du chapitre II.)

\th Th\'eor\`eme 4.1|On a
$$
\eqalignno{\noalign{\vskip-4pt}
\overline A(t,t',u)&=\sum_{0\le n}{u^n\over n!}\overline
A_n(t,t')
=\exp(ut'+C(t,u)),&(1)\cr
\noalign{\vskip-4pt}
\noalign{\hbox{o\`u}}
\noalign{\vskip-4pt}
C(t,u)&=\sum_{2\le n}{u^n\over n!} t\,A_{n-1}(t).&(2)\cr}
$$
\finth

\dem
Que $\theta'$ soit multiplicative d\'ecoule de la propri\'et\'e
3.11 et des d\'efinitions des vecteurs $E'\sigma$ et $\Delta
E\sigma$ $(\sigma\in {\goth S}_n)$. Par cons\'equent, le
membre de droite de l'identit\'e~(9) du chapitre~III devient
$\exp(\sum\limits_{1\le n}(u^n/n!)\,\theta'\{{\cal
C}_n\})$. Pour $n=1$, on a $\theta'\{{\cal
C}_n\}=t'$ et pour $n\ge 2$, d'apr\`es les propri\'et\'es
2.2 et~2.3, on a : $\theta'\{{\cal C}_n\}=\theta\Delta 
E\{{\cal C}_n\}=t\,A_{n-1}(t)$.\cqfd

\medskip
Le th\'eor\`eme 4.1 nous a donn\'e une identit\'e sur les
polyn\^omes $\overline A_n(t,t')$. Nous allons maintenant
trouver une formule explicite pour la fonction g\'en\'eratrice
$\overline A(t,t',u)$ en utilisant les r\'esultats de la
section~6 du chapitre~III.

\th Th\'eor\`eme 4.2|On a :
$$
\eqalignno{
\overline A(t,t',u)&=\sum_{0\le n}{u^n\over n!}
\overline A_n(t,t')
={1-t\over \exp((t-t')u)-t\exp((1-t')u)}.&(3)\cr
}
$$

\goodbreak
\noindent
En particulier,
$$
\eqalignno{\noalign{\vskip-8pt}
\overline A(t,t,u)&=\sum_{0\le n}{u^n\over n!}
\zA_n(t)
={1-t\over 1-t\exp((1-t)u)};&(4)\cr
\overline A(t,1,u)&=\sum_{0\le n}{u^n\over n!}
A_n(t)
={1-t\over -t+\exp((t-1)u)};&(5)\cr
\overline A(t,0,u)&=\sum_{0\le n}{u^n\over n!}
B_n(t)
={1-t\over \exp(ut)-t\exp(u)}.&(6)\cr
}
$$\finth

\dem
Avec les notations du th\'eor\`eme 3.13, si l'on pose $a=t$,
$b=t'$ et $c=1$, on a pour $\sigma\in{\goth S}_n$
$(n\ge 1)$ l'\'egalit\'e
$\xi_{1,\sigma(1)}\ldots\xi_{n,\sigma(n)}=\theta'\sigma$ ;
soit $\overline A_n(t,t')=\per \Xi_n$. La premi\`ere
identit\'e r\'esulte donc de la formule~(11) du chapitre~III.
En posant successivement $t'=t$, puis $t'=1$, enfin
$t'=0$, on obtient les trois suivantes.\cqfd

\rem Remarque $4.3$|Ces formules peuvent aussi s'obtenir
par le proc\'ed\'e suivant. D'apr\`es la propri\'et\'e~2.2, on a
$\zA_n(t)=t\,A_n(t)$ pour tout $n\ge 1$ ; on en tire
$$
\displaylines{\noalign{\vskip-8pt}
\overline A(t,t,u)=\sum_{0\le n}{u^n\over n!}
\zA_n(t)=1+t\sum_{1\le n}{u^n\over n!} A_n(t),\cr
\noalign{\vskip-4pt}
\noalign{\hbox{soit}}
\noalign{\vskip-4pt}
\hfill
\overline A(t,t,u)=1+t(\overline
A(t,1,u)-1).\hfill\llap{(7)}\cr
\noalign{\hbox{D'autre part, d'apr\`es le th\'eor\`eme 3.9 on a}}
\hfill
\exp(C(t,u))=\overline A(t,1,u)\exp(-u)\hfill\llap{(8)}\cr
\noalign{\hbox{ou encore}}
\hfill
\overline A(t,t,u)=\exp(ut-u)\overline
A(t,1,u).\hfill\llap{(9)}\cr}
$$
Du syst\`eme form\'e par les deux \'equations (7) et (9), on
d\'eduit imm\'ediatement les identit\'es~(4) et~(5). On calcule
ensuite $C(t,u)$ en utilisant la formule~(8) et l'on en tire
l'identit\'e~(3) en se servant du th\'eor\`eme~4.1.

\rem Remarque $4.4$|Les formules (4) et (5) sont connues
(\cf. Riordan [24], p.~215 et~39). La formule~(6) a \'et\'e
obtenue par Roselle [25], par les m\'ethodes traditionnelles
du calcul diff\'erentiel et int\'egral, dans le cas particulier
o\`u $B_n(t)=\theta M{\cal G}_n$ $(n\ge 1)$.

\medskip
Notons encore que du th\'eor\`eme 4.1 r\'esulte imm\'ediatement,
par simple d\'erivation, que la fonction g\'en\'eratrice ${\bf
A}=\overline A(t,1,u)$ est solution de l'\'equation
diff\'erentielle de Bernoulli :
$${\partial\over \partial u}{\bf A}={\bf A}(1+t({\bf
A}-1)).
$$
On peut aussi prouver ce r\'esultat directement et pour ce
faire, nous ferons la convention suivante que nous
utiliserons encore dans la section~2 : si $\sigma$ est
dans ${\goth S}_n$ $(n\ge 1)$, on consid\`ere~$\sigma w$
comme le mot $\sigma(1)\sigma(2)\ldots \sigma(n)$ dont
les lettres sont les \'el\'ements de~$[\,n\,]$ ; lorsque
$\sigma$ est l'\'el\'ement unique $\sigma_0\in{\goth S}_0$,
alors $\sigma w$ est le {\it mot vide} $\sigma_0 w$.
Si $f=y_1y_2\ldots y_m$ est un mot dont les lettres
$y_1$, $y_2$, \dots~, $y_m$ sont des entiers tous
distincts, on d\'esigne par $\omega$ l'unique morphisme
surjectif $\omega :\{y_1,y_2,\ldots, y_m\}\rightarrow
[\,m\,]$ et l'on note $\omega f$ le mot $\omega
y_1\,\omega y_2\,\ldots \omega\,y_m$. On pose encore
$\omega \sigma_0w=\sigma_0w$.

Prenons alors un mot $\sigma w\in {\goth S}_{n+1}$
$(n\ge 0)$ ; il s'\'ecrit univoquement $\sigma w=f(n+1)f'$.
Consid\'erons l'application $\sigma w\mapsto (\omega
f,\omega f')$ ; on a $\omega f\in {\goth S}_m$ et
$\omega f'\in {\goth S}_{n-m}$ pour un certain~$m$
tel que $0\le m\le n$ et puisque~$\omega$ est un
morphisme strictement croissant, on a encore
$$
|\Delta D\sigma|+\delta_{n,m}=|\Delta D \omega f|+|\Delta
D \omega f'|+1
$$
o\`u, comme d'usage, $\delta_{n,m}=1$ ou 0 selon que
$m=n$ ou $m\not=n$. D'autre part, l'image r\'eciproque par
l'application ci-dessus du couple $(\tau w, \tau'w)$ o\`u
$\tau\in {\goth S}_m$ et $\tau'\in {\goth S}_{n-m}$,
contient $n\choose m$ \'el\'ements. On en d\'eduit :
$$
A_{n+1}(t)=A_n(t)+t\sum_{0\le m\le n-1}{n\choose m}
A_m(t)A_{n-m}(t)\quad (n\ge 0).
$$
Il en r\'esulte que la fonction g\'en\'eratrice ${\bf
A}=\overline A(t,1,u)$ est bien solution de l'\'equation
diff\'erentielle pr\'ec\'edente. Ce r\'esultat a \'et\'e \'etabli pour la
premi\`ere fois par Riordan [23].

\titrecourant={\eightrm 2. FONCTION G\'EN\'ERATRICE
EXPONENTIELLE DES POLYN\^OMES}
\sectiona 2. Fonction g\'en\'eratrice exponentielle des
polyn\^omes $\rA_n(t)$|
Pour $r\ge 1$ posons
$$
\rA_n(t,u)=\sum_{r-1\le n}{u^{n-r+1}\over
(n-r+1)!}\,\rA_n(t).
$$
On a en particulier ${}^{1\kern-3pt}A(t,u)=\overline
A(t,1,u)$, dont on conna\^\i t d\'ej\`a la formule explicite (\cf.
(5)). Le but de la pr\'esente section est d'\'etablir l'identit\'e
remarquable suivante, due \`a Riordan ([24] p. 235)
$$
\rA(t,u)=(r-1)!\,\bigl({}^{1\kern-3pt}A(t,u)\bigr)^r.
$$

\rem Construction d'une bijection de
$\!\bigcup\limits_{0\le n}\!\!{\goth S}_{n+r-1}$ sur
${\goth S}_{r-1}\times\nobreak {\goth S}^{((r))}$|Il
s'agit d'une bijection sur le produit cart\'esien de
${\goth S}_{r-1}$ par le compos\'e partitionnel marqu\'e
${\goth S}^{((r))}$ $(r\ge 1)$. Pour d\'efinir ce dernier, il
faut munir l'ensemble~${\goth S}$ d'une
application~$\lambda$ ; nous prenons naturellement
l'application d\'efinie par
$$
\lambda\sigma=n\Leftrightarrow \sigma\in {\goth
S}_n\quad (n\in {\bboard N}).
$$
Notons que l'\'el\'ement unique $\sigma_0\in {\goth S}_0$
appartient \`a~$\goth S$, satisfait \`a $\lambda\sigma_0=0$
et est distinct de l'\'el\'ement neutre~$e$ du mono\"\i
de~${\goth S}^*$ pour lequel on a aussi $\lambda e=0$.

\goodbreak
Soit maintenant $\sigma\in {\goth S}_{n+r-1}$ $(0\le
n)$ ; le {\it mot} $\sigma w$ s'\'ecrit univoquement
$\sigma w=g_1i_1g_2i_2\ldots g_{r-1}i_{r-1}g_r$ o\`u
$\{i_1,i_2,\ldots,i_{r-1}\}=[r-1]$ et o\`u $g_1$, $g_2$,
\dots~, $g_r$ sont des mots (\'eventuellement vides) dont
les lettres sont des entiers. Soient $I_1$, $I_2$, \dots~,
$I_r$ les sous-ensembles de~$\bboard N$ dont les
\'el\'ements sont respectivement les lettres des mots $g_1$,
$g_2$,
\dots~, $g_r$. Posant $\sigma_j w=\omega g_j$ pour
chaque $j\in [r]$ (o\`u~$\omega$ est le morphisme d\'efini
\`a la fin de la section pr\'ec\'edente et o\`u l'on a
$\sigma_j=\sigma_0$ si le mot~$g_j$ est vide), on voit
imm\'ediatement que le mot
$$
h=(\sigma_1,I_1)(\sigma_2,I_2)\ldots (\sigma_r,I_r)
$$
est un \'el\'ement du compos\'e partitionnel marqu\'e ${\goth
S}^{((r))}$.

On note $\beta'\sigma\, (=\beta h$ dans les notations du
chapitre~III) le mot $\sigma_1\sigma_2\ldots
\sigma_r\in {\goth S}^*$ de longueur~$r$. Soit ensuite
$\overline\sigma w$ la permutation d\'efinie par
$\overline\sigma w=i_1i_2\ldots i_{r-1}$. On a
$\overline\sigma\in {\goth S}_{r-1}$ et comme
$\lambda h=\lambda
\beta'\sigma=\lambda\sigma_1+\cdots+\lambda\sigma_r
=\lambda\sigma-(r-1)=n$, on voit que l'application
$\sigma\mapsto (\overline\sigma,h)$ envoie ${\goth
S}_{n+r-1}$ dans 
${\goth S}_{r-1}\times {\goth S}^{((r))}\cap
\lambda^{-1}n$. Il est d'autre part imm\'ediat de v\'erifier
que cette application est bijective. Ceci ach\`eve la
construction de la bijection cherch\'ee.

\medskip
Maintenant, puisque $\card{\goth S}_{r-1}$ est \'egal \`a
$(r-1)!$, on peut \'ecrire
$$
\eqalignno{\sum\Bigl\{{\beta'\sigma\over
(\lambda\sigma-r+1)!} :\sigma\in \bigcup_{0\le n}{\goth
S}_{n+r-1}\Bigr\} 
&=(r-1)!\,\sum\Bigl\{{\beta h\over\lambda h!}
:h\in {\goth S}^{((r))}\Bigr\}\cr
&=(r-1)!\,\Bigl(\sum\Bigl\{{\sigma\over\lambda\sigma
!}:\sigma\in {\goth S}\Bigr\}\Bigr)^r,\cr}
$$
d'apr\`es le th\'eor\`eme 3.2. Notant $\varpi_r\sigma$ l'image
ab\'elienne de~$\beta'\sigma$, on obtient l'identit\'e
suivante valable dans l'alg\`ebre large sur~$\bboard Q$ du
mono\"\i de ab\'elien~${\goth S}^+$
$$
\sum\Bigl\{{\varpi_r(\sigma)\over(\lambda\sigma-r+1)!}
:\sigma\!\in\! \bigcup_{0\le n}{\goth S}_{n+r-1}\Bigr\}
=(r-1)!\,\Bigl(\sum\Bigl\{{\sigma\over\lambda\sigma
!}:\sigma\!\in\! {\goth S}\Bigr\}\Bigr)^r.\eqno(10)
$$
Soient enfin $\mu:{\goth S}^+\rightarrow \Omega$ un
morphisme dans un mono\"\i de ab\'elien~$\Omega$
et~$u$ une ind\'etermin\'ee. Comme d\'ej\`a vu au chapitre~III,
on v\'erifie que l'application $\sigma\mapsto u^{\lambda
\sigma}\mu\sigma$ $(\sigma\in {\goth S})$ peut \^etre
prolong\'ee en un morphisme continu~$\phi$ de la
${\bboard Q}$-alg\`ebre large de~${\goth S}^+$ dans
$\overline \Omega [[u]]$. Appliquant~$\phi$ aux deux
membres de l'identit\'e (10), on trouve
$$
\sum_{0\le n}{u^n\over n!}\,\mu\varpi_r\{{\goth
S}_{n+r-1}\}
=(r-1)!\,\Bigl(\sum_{0\le n}{u^n\over n!}\,\mu\{{\goth
S}_n\}\Bigr)^r.\eqno(11)
$$

\th Th\'eor\`eme 4.5|Pour $r\ge 1$ on a :
$\rA(t,u)=(r-1)!\,\bigl({}^{1\kern-3pt}A(t,u)\bigr)^r$.
\finth

\dem
Pour $r=1$, il n'y a rien \`a prouver. Supposons $r\ge 2$
et prenons pour~$\mu$ le morphisme prolongeant
l'application $\sigma\mapsto \theta\Delta
D\sigma=t^{|\Delta D\sigma |}$ au mono\"\i de~${\goth
S}^+$. Si l'on a $\sigma\in {\goth S}_{n+r-1}$ et
$\sigma w=g_1i_1g_2i_2\ldots g_{r-1}i_{r-1}g_r$
avec $\{i_1,i_2,\ldots, i_{r-1}\}=[r-1]$, il vient
$|\Delta D g_j|=|\Delta D \omega g_j|$ $(j\in [r])$,
puisque~$\omega$ est un morphisme injectif. D'autre
part, puisque le mot $g_1g_2,\ldots g_r$ contient toutes
les lettres du mot~$\sigma w$ sup\'erieures ou \'egales
\`a~$r$, on a $|\Delta'{}^{r-1}\Delta D\sigma|=\sum_j
|\Delta D\omega g_j|$, d'o\`u $\mu\varpi_r\sigma=\theta
\Delta'{}^{r-1}\Delta D\sigma$. On obtient alors, d'apr\`es
la propri\'et\'e~2.2, 
$\mu\varpi_r\{{\goth S}_{n+r-1}\}=\rA_{n+r-1}(t)$.
Comme on a d'autre part $\mu\{{\goth S}_n\}=A_n(t)$, le
th\'eor\`eme~4.5 r\'esulte de l'identit\'e~(11).\cqfd

\sectiona 3. Autres interpr\'etations des polyn\^omes
eul\'eriens|
Les techniques du chapitre pr\'ec\'edent pourraient \^etre
appliqu\'ees \`a d'autres probl\`emes d'\'enum\'eration. Par
exemple, au lieu d'introduire la fonction
multiplicative~$\theta'$ du paragraphe~1, on peut,
pour chaque $\sigma\in {\goth S}_n$ $(n\ge 1)$, poser
$\mu\sigma=\theta'\sigma\cdot r^{z(\sigma)}$, o\`u~$r$
est une ind\'etermin\'ee et o\`u $z(\sigma)$ est le nombre de
cycles de~$\sigma$. La fonction~$\mu$ est \'evidemment
multiplicative et l'on peut appliquer la proposition~3.12.
De plus, on a, comme dans la d\'emonstration du
th\'eor\`eme~4.1
$$
\eqalignno{\noalign{\vskip-8pt}
\mu\{{\cal C}_n\}&=t'\,r,&\hbox{pour $n=1$;}\cr
&=r\,\theta\Delta E{\cal C}_n=r\,t\,A_{n-1}(t),&\hbox{pour
$n\ge 2$.}\cr
\noalign{\hbox{La proposition 3.12 conduit donc \`a
l'identit\'e, dans laquelle 
$\mu\{{\goth S}_0\}=1$,}}
\sum_{0\le n}{u^n \over n!}\,\mu\{{\goth S}_n\}
&=\exp\Bigl(r(ut'+\sum_{2\le n}{u^n\over
n!}\,t\,A_{n-1}(t))\Bigr).&(12)\cr}
$$
Le membre de gauche de cette derni\`ere identit\'e est la
fonction g\'en\'eratrice exponentielle des permutations
class\'ees {\it \`a la fois par nombre de cycles et par nombre
d'exc\'edances}. Maintenant les identit\'es (3) et~(12), ainsi
que le th\'eor\`eme~4.1 permettent d'\'ecrire, lorsque~$r$ est
un entier positif,
$$
\sum_{0\le n}{u^n\over n!}\,\mu\{{\goth S}_n\}=
\bigl(\,\overline A(t,t',u)\bigr)^r.\eqno(13)
$$
On obtient donc d'apr\`es (3) la formule explicite de cette
fonction g\'en\'eratrice exponentielle.

\medskip
Si nous posons identiquement $t'=1$, nous obtenons
$$\mu\{{\goth S}_n\}=\sum\{t^{|\Delta
E\sigma|}r^{z(\sigma)}:\sigma\in {\goth S}_n\}$$ que
nous d\'esignons par $Q_n(t,r)$ $(n\ge 1)$. On posera
\'egalement $Q_0(t,r)=1$. Il r\'esulte de~(13) que l'on a,
lorsque~$r$ est un entier positif,
$$
\eqalignno{
(r-1)!\,\sum_{0\le n}{u^n\over n!}\,Q_n(t,r)
&=(r-1)!\,\bigl({}^{1\kern-3pt}A(t,u)\bigr)^r
\qquad\qquad\cr
\noalign{\vskip-5pt}
&=\rA(t,u)\qquad&\hbox{[d'apr\`es le th\'eor\`eme 4.5]}\cr
&=\sum_{0\le n}{u^n\over n!}\,\rA_{n+r-1}(t).\cr}
$$
On en d\'eduit une nouvelle interpr\'etation des polyn\^omes
$\rA_n(t)$, \`a savoir
$$
\rA_{n+r-1}(t)=(r-1)!\,Q_n(t,r)
\quad (r\ge 1).\eqno(14)
$$

\goodbreak
Enfin, d\'esignons par $s(\sigma)$ le nombre des \'el\'ements
saillants de la suite $\sigma w$ o\`u $\sigma\in {\goth
S}_n$ $(n\ge 1)$. Comme l'on a $|M\hat\sigma|+|\Delta
E\sigma|=n$ et $s(\hat\sigma)=z(\sigma)$, on obtient
$$
\eqalignno{\sum\{t^{|M\sigma|}r^{s(\sigma)}:\sigma\in
{\goth S}_n\}
&=t^n\sum\{t^{-|\Delta E\sigma|}r^{z(\sigma)}:\sigma\in
{\goth S}_n\}\cr
&=t^n\,Q_n(t^{-1},r).&(15)\cr}
$$
Le premier membre de l'identit\'e (15) est le polyn\^ome
g\'en\'erateur des permutations $\sigma\in {\goth S}_n$
class\'ees \`a la fois suivant {\it leur nombre d'\'el\'ements
saillants et leur nombre de mont\'ees}. Ce polyn\^ome
g\'en\'erateur a \'et\'e consid\'er\'e pour la premi\`ere fois par
Dillon et Roselle~[10] qui ont \`a son sujet prouv\'e un
certain nombre d'identit\'es, qu'on pourrait retrouver \`a
partir des formules~(14) et~(15) et des r\'esultats de ce
chapitre.

\medskip
Enfin, notons que la propri\'et\'e 2.6 fait appara\^\i tre que les
coefficents des polyn\^omes $\rA_n(t)$ sont tous divisibles
par~$r!$ (ce que ne fait pas appara\^\i tre le th\'eor\`eme~4.5).
Si donc on pose
$$
\rA_n(t)=r!\,{}^{r\kern-2pt}P_n(t),
$$
il semble int\'eressant d'obtenir une interpr\'etation pour les
polyn\^omes ${}^{r\kern-2pt}P_n(t)$
$(0\le r\le n)$.

D'abord, si $r=n$, on a ${}^{r\kern-2pt}P_n(t)=1$ et
pour $r=0$ et~1, on a ${}^{r\kern-2pt}P_n(t)=\rA_n(t)$.
On fait donc l'hypoth\`ese $2\le r\le n-1$. La restriction
de tout $\sigma\in {\goth S}_n$ \`a $[n-r]$ est une
injection de $[n-r]$ dans $[\,n\,]$ que nous noterons
$\phi\sigma$. L'application~$\phi$ est \'evidemment une
surjection de~${\goth S}_n$ sur l'ensemble ${\cal
I}_{n-r,n}$ des injections de $[n-r]$ dans~$[\,n\,]$ telle
que l'image inverse de tout $\tau \in {\cal I}_{n-r,n}$
a~$r!$ \'el\'ements. Introduisons l'application $\Lambda$ de
${\bboard N}^p$ $(p\ge 1)$ dans lui-m\^eme envoyant
chaque vecteur $(x_1x_2,\ldots, x_p)$ sur $((x_1-1)_+,
(x_2-1)_+,\ldots,(x_p-1)_+)$. On a ainsi
$\Delta=\Delta''\Lambda=\Lambda\Delta''$ et
$\Delta^r=\Delta''{}^r\Lambda^r$ pour $r\ge 1$.
D\'efinissant le vecteur-exc\'edance d'une injection de fa\c con
\'evidente, on a ainsi
$$
\displaylines{\noalign{\vskip-5pt}
\Delta^rE\sigma=\bigl((\sigma(1)-r)_+,\ldots,
(\sigma(n-r)-r)_+\bigr)
=\Lambda^rE\phi \sigma.\cr\noalign{\hbox{D'o\`u l'on
d\'eduit
$\Delta^rE{\goth S}_n=r!\,\Lambda^rE\,{\cal I}_{n-r,n}$
et par suite}}
\rA_n(t)\,{1\over r!}=
{}^{r\kern-2pt}P_n(t)=\theta\Lambda^rE\,{\cal
I}_{n-r,n}.\cr}
$$
Cette derni\`ere interpr\'etation des polyn\^omes $\rA_n(t)/r!$
est due \`a Strosser [28].

\vfill\eject

\pagetitretrue

\auteurcourant={\eightrm CHAPITRE V :
LES SOMMES ALTERN\'EES}

\def\mb{{\overline m}}
\def\db{{\overline d}}

\titrecourant={\eightrm 1. DISTRIBUTION DU NOMBRE DES
DESCENTES}
\vglue 2cm
\centerline{{\eightrm CHAPITRE V}}
\vskip 2mm
\centerline{\bf LES SOMMES ALTERN\'EES
$A_n(-1)$ ET $B_n(-1)$}

\vskip 6mm plus 2mm
\sectiona 1. Distribution du
nombre des descentes sur ${\goth S}_n'$|
Nous attachons \`a chaque $\sigma\in {\goth S}_n'\,(=
\{\sigma\in {\goth S}_n:\sigma(1)=n\}$) un mot not\'e
$V(\sigma)=v_1v_2\ldots v_n$ dans les lettres de
l'alphabet $X=\{m,\mb,d,\db\}$ par les r\`egles suivantes,
o\`u, par d\'efinition, $\sigma(n+1)=\sigma(1)\ (\,=n)$.

\decale (1)|Pour chaque $j\in [\,n\,]$, on a $v_j\in
\{d,\db\}$ ou $v_j\in \{m,\mb\}$ selon que
$\sigma(j)>\sigma(j+1)$ ou $\sigma(j)<\sigma(j+1)$ ;

\decale (2)|si $v_j\in \{d,\db\}$, $j\in [n-1]$, on a
$v_j=d$ ou~$\db$ selon que $v_{j+1}$ est dans
$\{d,\db\}$ ou dans $\{m,\mb\}$ ;

\decale (3)|si $v_j\in \{m,\mb\}$, $(2\le j\le n)$, on a
$v_j=m$ ou~$\mb$ selon que $v_{j-1}$ est dans
$\{m,\mb\}$ ou dans $\{d,\db\}$.

\medskip
En raison de $\sigma(1)=\sigma(n+1)=n$, on a toujours
$v_1\in \{d,\db\}$, $v_n\in \{m,\mb\}$ et les seules
occurrences des lettres~$\db$ et~$\mb$ se rencontrent
dans les facteurs $v_jv_{j+1}=\db\,\mb$ correspondant
aux indices $j\in [n-1]$ tels que
$\sigma(j)>\sigma(j+1)<\sigma(j+2)$. Par exemple, pour
$\sigma (w)=(7,1,4,6,3,2,5)$, on aurait
$V(\sigma)=\db\,\mb\,m\,d\,\db\,\mb\,m$.

\medskip
Introduisons maintenant pour tout lettre~$x$ et tout
mot~$g$ la {\it d\'erivation} $g\,(\partial/\partial x)$
envoyant chaque mot $f=x_1x_2\ldots x_p$ sur l'ensemble
pond\'er\'e form\'e de tous les mots obtenus en rempla\c cant
dans~$f$ chaque occurrence de la lettre~$x$ par le
mot~$g$. Formellement, $g\,(\partial/\partial x)$ est
l'op\'erateur lin\'eaire d\'efini par sa restriction \`a~$X$, \`a
savoir
$$\eqalignno{
\Bigl(g{\partial\over \partial x}\Bigr)x'
&=\cases{ g,&si $x'=x$ ;\cr
x',&si $x'\in X\setminus\{x\}$ ;\cr}\cr
\noalign{\hbox{et par l'identit\'e}}
\Bigl(g{\partial\over \partial x}\Bigr)ff'
&=
\Bigl(g{\partial\over \partial x}\Bigr)f\cdot f'+
f\cdot \Bigl(g{\partial\over \partial x}\Bigr)f'.\cr}
$$
Donc si $f=f_1xf_2x\ldots f_{r-1}xf_r$, o\`u les~$f_i$ ne
contiennent pas la lettre~$x$, l'on aura :
$$\displaylines{
\Bigl(g{\partial\over \partial x}\Bigr)f
=f_1gf_2x\ldots f_{r-1}xf_r
+f_1xf_2g\ldots f_{r-1}xf_r+\cdots+
f_1xf_2x\ldots f_{r-1}gf_r.\cr
\noalign{\hbox{Par exemple, on a :}}
\Bigl(\db\,\mb{\partial\over \partial m}\Bigr)
(\db\,\mb\,m\,d\,\db\,\mb\,m)
=\db\,\mb\,\db\,\mb\,d\,\db\,\mb\,m+
\db\,\mb\,m\,d\,\db\,\mb\,\db\,\mb.\cr}
$$

\goodbreak
\th Lemme 5.1|Soit
$$
\nabla=\Bigl(d\,\db\,{\partial \over \partial \db}\Bigr)
+\Bigl(\mb\,m\,{\partial \over \partial\, \mb}\Bigr)
+\Bigl(\db\,\mb{\partial \over \partial d}\Bigr)
+\Bigl(\db\,\mb\,{\partial \over \partial m}\Bigr).$$
On a identiquement :
$V{\goth S}_n'=\nabla\,V{\goth S}_{n-1}'$ $(n\ge 2)$.
\finth

\dem
Il existe une bijection de ${\goth S}_{n-1}'\times [\,n\,]$
sur ${\goth S}_{n}'$ envoyant chaque $(\sigma',k)\in
{\goth S}_{n-1}'\times [\,n\,]$ sur la permutation
$\sigma\in {\goth S}_n'$ telle que $\sigma w$ soit
obtenue en ajoutant~1 \`a tous les chiffres de~$\sigma' w$
et en ins\'erant 1 entre le $k$\ieme\ et le $(k+1)$\ieme\
terme de~$\sigma'w$. Soit $V(\sigma')=v_1'v_2'\ldots
v'_{n-1}$ et supposons $v'_k\in \{d,\db\}$ c'est-\`a-dire
$1\le k\le n-2$ et $\sigma'(k)>\sigma'(k+1)$. On a
$\sigma(k)=1+\sigma'(k)>\sigma(k+1)
=1<\sigma(k+2)=1+\sigma'(k+1)$,
donnant dans $V(\sigma)$ le facteur
$v_kv_{k+1}=\db\,\mb$. Maintenant :

\decale (i)|si $v'_k=\db$, c'est-\`a-dire si
$v'_{k+1}=\mb$ et $\sigma'(k+1)<\sigma'(k+2)$, on a
$v_{k+2}=m$ puisque
$\sigma(k+2)=1+\sigma'(k+1)<\sigma(k+3)
=1+\sigma'(k+2)$
et toute l'op\'eration \'equivaut au remplacement de
$v'_{k+1}=\mb$ par $v_{k+1}v_{k+2}=\mb\,m$.
Remarquons qu'avec nos conventions, si $k=n-2$, on a
$\sigma'(k+2)=\sigma'(n)=\sigma'(1)=n-1$.

\decale(ii)|si $v'_k=d$, c'est-\`a-dire si
$\sigma'(k+1)<\sigma'(k+2)$, on a encore $v_{k+2}\in
\{d,\db\}$ et $V(\sigma)$ est d\'eduit de $V(\sigma')$ en
rempla\c cant $v'_k=d$ par $v_kv_{k+1}=\db\,\mb$.
Un raisonnement analogue s'applique si $v'_k=m$ ou
$\mb$.\cqfd

\medskip
Notons maintenant $\alpha$ le morphisme canonique
envoyant le mono\"\i de libre engendr\'e par
$\{m,\mb,d,\db\}$ sur le mono\"\i de commutatif libre de
m\^eme base.

\th Th\'eor\`eme 5.2|Il existe des entiers positifs $c_{n,k}$
tels que
$$
\alpha V{\goth S}_n'=\sum_{2\le 2k\le n}
c_{n,k}\,(\db\,\mb)^k\,(d+m)^{n-2k}
\quad (n\ge 2).
$$
\finth

\dem Pour $n=2$, on a $V{\goth S}'_n=\db\,\mb$ et le
r\'esultat s'en d\'eduit par induction sur~$n$
puisque~$\nabla$ commute avec~$\alpha$.\cqfd

\rem Remarque $5.3$|Les coefficients $c_{n,k}$ des
polyn\^omes $\alpha V{\goth S}_n'$ ob\'eissent \`a des
relations de r\'ecurrence qu'il est facile d'\'etablir. Posons,
par convention, $c_{n,k}=0$ si $k\le 0$ ou si $2k\ge
n+1$ ; on  alors les deux relations :
$$
\eqalignno{c_{2,1}&=1\quad \hbox{ et pour $n\ge 3$,
$k\ge 1$,}\cr
c_{n,k}&=kc_{n-1,k}+2(n+1-2k)c_{n-1,k-1}.\cr}
$$

\rem Remarque $5.4$|La fonction g\'en\'eratrice des nombres
$c_{n,k}$ est donn\'ee par Barton \& David ([3] p. 180,
voir aussi [17]). La th\'eorie de ces auteurs se rattache aux
consid\'erations pr\'esentes en utilisant l'observation suivante
dont la d\'emonstration est laiss\'ee au lecteur.

Pour $\sigma\in{\goth S}_{n-1}$, soit $\sigma'\in
{\goth S}_{n}'$ d\'efinie par $\sigma'(1)=n$,
$\sigma'(1+j)=n+1-\sigma(n-j)$. Le nombre des
facteurs~$d$ ou~$\db$ de~$V(\sigma')$ surpasse de~1 le
nombre des $j\in [n-1]$ tels que
$\sigma(j)>\sigma(j-1)$ et le nombre de facteurs
$\db\,\mb$ de $V(\sigma')$ est \'egal au nombre des $j\in
[n-2]$ tels que $\sigma(j)>\sigma(j+1)<\sigma(j+2)$,
augment\'e d'une unit\'e.

\goodbreak
\sectiona 2. Applications aux polyn\^omes eul\'eriens|
Le th\'eor\`eme 5.2 va nous permettre de donner une
interpr\'etation combinatoire aux nombres $A_n(-1)$. Il
est commode, tout d'abord, de noter la relation suivante
sur les cardinaux des ensembles~${\cal T}_n$ des
permutations {\it altern\'ees} (\cf. chap.~I, \S\kern2pt 9).

\th Propri\'et\'e 5.5|Pour $p\ge 1$, on a :
$\card {\cal T}_{2p-1}=\card {\cal T}_{2p}\cap
{\goth S}_{2p}'$.
\finth

\dem
En effet, l'application qui envoie chaque $\sigma'\in
{\goth S}_n'$ telle que $V(\sigma')=(\db\,\mb)^p$
$(n=2p\ge 2$) sur l'\'el\'ement $\sigma\in {\goth S}_{n-1}$
d\'efini par $\sigma(j)=\sigma'(j+1)$ $(j\in [n-1])$, est
une bijection sur ${\cal T}_{n-1}$. D'autre part, il est
clair que l'on a :
${\cal T}_{2p}\cap {\goth S}_{2p}'
=\{\sigma'\in {\goth S}_{2p}':
V(\sigma')=(\db\,\mb)^p\}$.\cqfd

\th Th\'eor\`eme 5.6|Pour $n\ge 2$ on a l'identit\'e :
$$\displaylines{\hfill
t\,A_{n-1}(t)=\sum_{2\le 2k\le
n}c_{n,k}\,t^k(1+t)^{n-2k}.\hfill\llap{\rm (1)}\cr
\noalign{\hbox{De plus, pour $p\ge 1$, on a :}}
A_{2p}(-1)=0\quad{\it et}\quad
(-1)^{p-1}\,A_{2p-1}(-1)=\card {\cal T}_{2p-1}.\cr}$$
\finth

\dem
Pour $\sigma'\in {\goth S}'_n$, le nombre des occurrences
des lettres~$d$ ou~$\db$ dans $V(\sigma')$ est \'egal \`a
$|\Delta D\sigma'|$, puisque l'on a $v_n\in \{m,\mb\}$.
Or d'apr\`es les propri\'et\'es~2.2 et~2.3, on a $\theta \Delta D
{\goth S}_n'={}^{0\kern-3pt}A_{n-1}(t)= t\,A_{n-1}(t)$.
On en d\'eduit que $t\,A_{n-1}(t)$ est obtenu en faisant
$d=\db=t$ et $m=\mb=1$ dans $\theta V{\goth S}'_n$.
La formule~(1) r\'esulte alros du th\'eor\`eme~5.2.

D'autre part, le second membre de la formule~(1) admet le
facteur $(1+t)$ si~$n$ est impair. On en conclut que
$A_{2p}(-1)$ est nul pour $p\ge 1$. Au contraire, pour
$n=2p\ge 2$, on voit que
$c_{2p,p}=(-1)^{p-1}A_{2p-1}(-1)$. Or
$$
c_{2p,p}=\card\{\sigma'\in {\goth
S}'_{2p}:V(\sigma')=(\mb\,\db)^p\}
=\card {\cal T}_{2p}\cap {\goth S}'_{2p}=\card {\cal
T}_{2p-1},
$$
d'apr\`es la propri\'et\'e 5.5. Le th\'eor\`eme 5.6 en r\'esulte.\cqfd

\titrecourant={\eightrm 3. APPLICATIONS AUX
POLYN\^OMES}

\sectiona 3. Applications aux polyn\^omes $B_n(t)$|
Nous donnons enfin des identit\'es analogues \`a celles du
th\'eor\`eme~4.6, concernant les polyn\^omes $B_n(t)=\theta
E{\cal D}_n=\theta M{\cal G}_n$, o\`u comme pr\'ec\'edemment
$$
{\cal D}_n=\{\sigma\in {\goth S}_n:\sigma(j)\not=j\};\quad
{\cal G}_n=\{\sigma\in {\goth S}_n:
1\not=\sigma(1),\,1+\sigma(j)\not=\sigma(j+1)\}.
$$
Pour d\'emontrer le th\'eor\`eme 5.9 ci-dessous, nous allons de
nouveau appliquer la formule exponentielle et utiliser les
propri\'et\'es \'el\'ementaires des permutations et de la
transformation fondamentale du chapitre~I.
Pour $n=2p\ge 2$ et $\sigma\in {\goth S}_n$, nous
posons $\mu\sigma=1$ si et seulement si~$\sigma$ est
{\it biexc\'ed\'ee}, c'est-\`a-dire si $\sigma\in {\cal B}$ (\cf.
chap.~1, \S\kern2pt9) et $\mu\sigma=0$ dans les autres
cas.

\th Lemme 5.7|L'application $\mu$ est multiplicative. En
d'autres termes,~$\sigma$ est une permutation biexc\'ed\'ee,
si tous les termes de sa factorisation cano\-nique sont
aussi des permutations biexc\'edes.
\finth

\dem
Ce lemme r\'esulte encore de la propri\'et\'e 3.11. Soit
$\sigma_1\sigma_2\ldots \sigma_r$ la d\'ecomposition en
produit de cycles disjoints d'une permutation $\sigma\in
{\cal B}$ et $(f_1,I_1)(f_2I_2)\cdots (f_r,I_r)$ sa
factorisation canonique. Avec les m\^emes notations que
dans la propri\'et\'e~3.11, on a
$f_j=\tau_j\sigma_j\tau_j^{-1}$ $(j\in [r])$. Par suite,
$i<\sigma(i)\Leftrightarrow \tau_j(i)<f_j\tau_j(i)$ et
$i<\sigma^{-1}(i)\Leftrightarrow
\sigma\sigma^{-1}(i)<\sigma^{-1}(i)
\Leftrightarrow f_j\tau_j\sigma^{-1}(i)
<\tau_j\sigma^{-1}(i)
\Leftrightarrow\tau_j(i)<f_j^{-1}\tau_j(i)$,
puisque les entiers $i$ et $\sigma^{-1}(i)$
appartiennent \`a la m\^eme orbite. On a les m\^emes
\'equivalences en rempla\c cant le symbole ``$<$" par
``$>$".\cqfd

\th Lemme 5.8|Pour $p\ge 1$ on a :
$$
\displaylines{
\hfill
\mu\{{\goth S}_{2p}\}=\card {\cal B}_{2p}\quad
{\it et}\quad
\mu\{{\goth S}_{2p-1}\}=\mu \{{\cal C}_{2p-1}\}=0;
\hfill
\llap{\rm (2)}\cr
\hfill
\mu\{{\cal C}_{2p}\}=\card {\cal B}_{2p}\cap {\cal
C}_{2p}=\card {\cal T}_{2p}\cap {\goth S}_{2p}'.\hfill
\llap{\rm (3)}\cr}
$$
\finth

\dem
Les relations (2) r\'esultent de la d\'efinition de $\mu$ et de
la proposition 1.14. D'apr\`es la propri\'et\'e~1.10 et la
proposition~1.14, la transformation fondamentale
$\sigma\mapsto \hat\sigma$ est une bijection de ${\cal B}_{2p}\cap {\cal
C}_{2p}$ sur ${\cal T}_{2p}\cap{\goth
S}_{2p}'$. La relation (3) est ainsi v\'erifi\'ee.\cqfd

\medskip
Pour la d\'emonstration du th\'eor\`eme ci-dessous,
l'utilisation des nombres complexes est une simple
commodit\'e d'\'ecriture \'evitant de recourir au produit
d'Hadamard.

\th Th\'eor\`eme 5.9|Pour $p\ge 1$ on a:
$$B_{2p-1}(-1)=0\quad{\it et}\quad
(-1)^pB_{2p}(-1)=\card {\cal T}_{2p}.$$
\finth

\dem
On applique la formule (9) du chapitre~III avec
l'application multiplicative~$\mu$ du lemme~5.7. 
D'apr\`es~(2), le
premier membre de cette formule s'\'ecrit :
$1+\sum\limits_{1\le n}(u^n/n!)\,\card{\cal B}_n$.
Maintenant pour $p\ge 1$ on~a
$$
\eqalignno{
\mu\{{\cal C}_{2p}\}&=\card {\cal T}_{2p}\cap {\goth
S}'_{2p}&\hbox{[d'apr\`es (3)]}\cr
&=\card {\cal T}_{2p-1}&\hbox{[d'apr\`es la propri\'et\'e
5.5]}\cr
&=(-1)^{p-1}A_{2p-1}(-1).\qquad
\qquad
&\hbox{[d'apr\`es le
th\'eor\`eme 5.6]}\cr}
$$
Compte-tenu de la relation (2) le second membre de la 
formule~(9) du chapitre~III s'\'ecrit donc
$$\displaylines{
\exp\Bigl(\sum_{1\le p}{u^{2p}\over
(2p)!}\,(-1)^{p-1}A_{2p-1}(-1)\Bigr).\cr
\noalign{\hbox{En utilisant le fait que $A_n(-1)=0$ si
$n$ est pair, on en d\'eduit :}}
\sum_{0\le n}{u^n\over n!}\,\card {\cal B}_n=\exp\Bigl(
\sum_{2\le n}{(iu)^n\over n!}\,(-1)A_{n-1}(-1)\Bigr),\cr}
$$

\goodbreak\noindent
o\`u $i$ est le nombre complexe de module~1 et d'argument
$\pi/2$. Or le membre de droite de cette derni\`ere
\'equation est la valeur pour $t=-1$ de l'expression
$$\displaylines{
\exp\Bigl(\sum_{2\le n}{(iu)^n\over n!}\,t\,
A_{n-1}(t)\Bigr),\cr
\noalign{\hbox{qui d'apr\`es le th\'eor\`eme~4.1 est \'egale \`a}}
\sum_{0\le n}{(iu)^n\over n!}\,B_n(t).\cr}
$$
En identifiant terme \`a terme, il en r\'esulte que l'on a
$B_n(-1)=0$ si~$n$ est impair et que pour $n=2p$, on
a $(-1)^p\,B_{2p}(-1)=\card {\cal B}_{2p}=\card {\cal
T}_{2p}$ d'apr\`es la proposition~1.14.\cqfd

\sectiona 4. Les d\'eveloppements de $\tg u$ et de $1/\cos
u$|
De l'identit\'e (5) du th\'eor\`eme 4.2, on tire
$$
\displaylines{
\sum_{0\le n}{(iu)^n\over n!}\,A_n(-1)={2\over
1+e^{-2iu}},\cr
\noalign{\hbox{qu'on peut r\'ecrire}}
\sum_{1\le n}{u(iu)^{n-1}\over n!}\,A_n(-1)
={1-e^{-2iu}\over i(1+e^{-2iu})}=\tg u,\cr
\noalign{\hbox{soit, en utilisant le th\'eor\`eme 5.6,}}
\hfill
\tg u=\sum_{1\le p}{u^{2p-1}\over (2p-1)!}\,
\card {\cal T}_{2p-1}.\hfill\llap{\rm (4)}\cr
\noalign{\hbox{De m\^eme, d'apr\`es l'identit\'e (6) du
th\'eor\`eme 4.2, on a}}
\sum_{0\le n}{(iu)^n\over n!}\,B_n(-1)={2\over
e^{-iu}+e^{iu}}={1\over \cos u}.\cr
\noalign{\hbox{D'apr\`es le th\'eor\`eme 5.9, on d\'eduit donc :}}
\hfill{1\over \cos u}=1+\sum_{1\le p}{u^{2p}\over
(2p)!}\,\card {\cal T}_{2p}.\hfill\llap{(5)}\cr
}
$$
Les identit\'es (4) et (5) sont dues \`a D\'esir\'e Andr\'e [1]. Nous
avons pu les \'etablir ici sans recourir aux m\'ethodes
traditionnelles du calcul diff\'erentiel et int\'egral, en
n'utilisant que l'identit\'e de Cauchy et des constructions
sur la cat\'egorie des ensembles totalement ordonn\'es finis.

Nous laissons au lecteur l'amusement de v\'erifier par
les m\^emes techniques la formule \'el\'ementaire
$$
{1\over \cos u}=\exp\Bigl(\int \tg u\,du\Bigr)
$$
en utilisant une d\'efinition appropri\'ee de l'int\'egrale.

\goodbreak

\titrecourant={\eightrm 5. TABLE DES NOMBRES D'EULER}

\sectiona 5. Table des nombres d'Euler|
On a souvent appel\'e {\it nombres d'Euler} les coefficients
du d\'eveloppement de $\tg u$ et de $1/\cos u$. Les
valeurs num\'eriques de ces premiers coefficients ont d\'ej\`a
\'et\'e obtenues par Euler lui-m\^eme (voir [16],
p.~299--301). Nous reproduisons ci-dessous ces
premi\`eres valeurs. Rappelons que pour $n\ge 1$ on note
${\cal T}_n$ le sous-ensemble de~${\goth S}_n$ form\'e
par les permutations {\it altern\'ees}, qu'on a ensuite les
identit\'es
$$
(-1)^{p-1}\,A_{2p-1}(-1)=\card {\cal T}_{2p-1};\quad
(-1)^{p-1}\,B_{2p}(-1)=\card {\cal T}_{2p}
\quad (p\ge 1).
$$
Le tableau des quantit\'es
$t_n=\card {\cal T}_n$ pour $n=1,2,\ldots, 14$ est alors
le suivant.

$$
\vbox{\halign{\vrule\thinspace \hfil$#$\ \vrule
&\strut\ $#$\hfil\  \vrule\cr
\noalign{\hrule}
n&t_n\cr
\noalign{\hrule}
1&1\cr
2&1\cr
3&2\cr
4&5\cr
5&16\cr
6&61\cr
7&272\cr
8&1385\cr
9&7936\cr
10&50521\cr
11&353792\cr
12&2702765\cr
13&22368256\cr
14&199360981\cr
\noalign{\hrule}
}}
$$

\vfill\eject
\pagetitretrue
\titrecourant={\eightrm BIBLIOGRAPHIE}
\auteurcourant={\eightrm BIBLIOGRAPHIE}
\vglue 44pt\vskip0pt plus 5pt\raggedbottom

{\eightpoint

\centerline{\bf BIBLIOGRAPHIE}

\bigskip

\article 1|D\'esir\'e Andr\'e|D\'evelopements de $\sec
x$ et de ${\rm tang}\,x$|C. R. Acad. Sc.
Paris|88|1879|965--967|

\article 2|D\'esir\'e Andr\'e|M\'emoire sur le nombre des
permutations altern\'ees|Journal de
Math.|7|1881|167|

\livre 3|D. E. Barton, F. N. David|Combinatorial
Chance|Griffin, London, {\oldstyle 1962}|

\article 4|L. Carlitz|Eulerian numbers and
polynomials|Math. Magazine|32|1959|247--260|

\article 5|L. Carlitz|Eulerian numbers and polynomials
of higher order|Duke Math. J.|27|1960|401--423|

\article 6|L. Carlitz|A note on Eulerian
numbers|Arch . Math.|14|1963|383--390|

\article 7|L. Carlitz, J. Riordan|Congruences for
Eulerian Numbers|Duke Math. J.|20|1953|339--343|

\article 8|L. Carlitz, D. P. Roselle and R. A.
Scoville|Permutations and Sequences with Repetitions by
Number of Increases|J. Combin. Theory|1|1966|350--374|

\livre 9|P. Cartier, D. Foata|Probl\`emes combinatoires
de commutation et r\'earrangements|Lecture Notes
in Math., no.~85, Springer-Verlag, Berlin, {\oldstyle 1969}|

\article 10|J. F. Dillon, D. P. Roselle|Eulerian numbers
of higher order|Duke Math. J.|35|1968|247--256|

\article 11|R. C. Entringer|A combinatorial
interpretation of the Euler and Bernoulli
numbers|Nieuw Arch. V. Wiskunde|14|1966|241--246|

\article 12|Dominique Foata|\'Etude alg\'ebrique de
certains probl\`emes d'analyse combinatoire et du
calcul des probabilit\'es|Publ. Inst.
Statist. Univ. Paris|14|1965|81--241|

\divers 13|G. Frobenius|\"Uber die Bernoullischen Zahlen
und die Eulerschen Polynome, {\sl Sitz. Berichte Preuss.
Akad. Wiss.}, {\oldstyle 1910}, p. 808--847|

\divers 14|R. Frucht|A combinatorial approach to the Bell
polynomials and their generalizations. {\sl Recent
Progress in Combinatorics} (W. T. Tutte, ed.), Academic
Press, London and New York, {\oldstyle 1964}, p. 69--74|

\article 15|R. Frucht, G.-C. Rota|Polynomios de Bell y
partitiones de conjuntos
finitos|Scientia|126|1965|5--10|

\livre 16|Ch. Jordan|Calculus of finite differences|R\"ottig
\& Romwalter, Budapest, {\oldstyle 1939}|

\article 17|W.O. Kermack, A.G. McKendrick|Some
distributions associated with a randomly arranged set
of numbers|Proc. Roy. Soc. Edinburgh
Sec. A.|57|1937|332--376|

\divers 18|B. Kittel|Communication priv\'ee|

\article 19|P. A. MacMahon|Second memoir on the
composition of numbers|Phil. Trans.
Royal Soc. London, A|207|1908|65--134|

\livre 20|P. A. MacMahon|Combinatory Analysis,
{\rm vol.~1 and~2}|Cambridge, Cambridge Univ. Press,
{\oldstyle 1915},  (Reprinted by Chelsea, New York,
{\oldstyle 1955})|

\livre 21|E. Netto|Lehrbuch der Combinatorik|B. G. Teubner,
Leipzig, {\oldstyle 1900}|

\article 22|F. Poussin|Sur une propri\'et\'e
arithm\'etique de certains polyn\^omes associ\'es aux
nombres d'Euler|C. R. Acad. Sc. Paris|266|1968|392--393|

\article 23|J. Riordan|Triangular permutation
numbers|Proc. Amer. Math. Soc.|2|1951|404--407|

\livre 24|J. Riordan|An Introduction to
Combinatorial Analysis|New York, J. Wiley, {1958}|

\article 25|D. P. Roselle|Permutations by number of
rises and successions|Proc. Amer. Math.
Soc.|19|1968|8--16|

\article 26|G.-C. Rota|On the foundations of
Combinatorial Theory|J.
Wahrscheinlichkeits\-theorie|2|1966|340--368|

\article 27|E. B. Shanks|Iterated sums of powers of
binomial coefficients|Amer. Math.
Monthly|58|1951|404--407|

\divers 28|R. Strosser|S\'eminaire de th\'eorie combinatoire,
I.R.M.A., Universit\'e de Strasbourg,  {\oldstyle
1969}--{\oldstyle 1970}|

\divers 29|G. E. Uhlenbeck, G. W. Ford|The theory
of graphs with applications to the virial development of
the properties of gases, {\sl Studies in Statistical
Mechanics}, vol.~1 (J.~de Boer and G. E. Uhlenbeck, eds.),
North-Holland, Amsterdam, {\oldstyle 1962}, p. 119--211|

\divers 30|Ph. Welschinger|S\'eminaire de th\'eorie
combinatoire, I.R.M.A., Universit\'e de Strasbourg,  {\oldstyle
1969}--{\oldstyle 1970}|

\article 31|J. Worpitzky|Studien \"uber die
Bernoullischen und Eulerschen Zahlen|J.
f\"ur die reine und angewandte Math.|94|1883|203--232|

}

\bye